\input ./style/arxiv-general.cfg
\documentclass[aap,MSNbibl,seceqn,nameyear,dvips]{arximspdf}
\makeatletter
   \@ifpackageloaded{graphicx}{}{\usepackage{graphicx}}
\makeatother

\doi{10.1214/14-AAP1092}
\volume{25}
\issue{6}
\pubyear{2015}
\firstpage{3624}
\lastpage{3683}
\docsubty{FLA}

\makeatletter
\renewcommand{\mid}{|}
\def\sfrac#1#2{#1/#2}
\def\vfrac#1#2{(#1)/#2}
\def\afrac#1#2{#1/(#2)}

\def\sklfrac#1#2{(#1/#2)}

\newcommand{\rrvert}{\vert}
\newcommand{\rrVert}{\Vert}
\newcommand{\llvert}{\vert}
\newcommand{\llVert}{\Vert}
\newcommand{\erE}{\mathrm{E}}
\newcommand{\overset}{\stackrel}
\newcommand{\bbB}{\mathbf{B}}
\newcommand{\bbC}{\mathbf{C}}
\newcommand{\bbX}{\mathbf{X}}
\newcommand{\bbI}{\mathbf{I}}
\newcommand{\bgL}{\bolds{\Lambda}}
\newcommand{\bbF}{\mathbf{F}}                                            
\newcommand{\bbe}{\mathbf{e}}
\newcommand{\bgma}{\bolds{\gamma}}
\newcommand{\bbM}{\mathbf{M}}
\newcommand{\bbbeta}{\bolds{\beta}}
\newcommand{\bbR}{\mathbf{R}}
\newcommand{\gd}{{\delta}}
\newcommand{\rE}{\mathrm{E}}
\newcommand{\rtr}{\operatorname{tr}}
\newcommand{\bbA}{\mathbf{A}}
\newcommand{\bga}{\bolds{\alpha}}
\newcommand{\ep}{\varepsilon}
\newcommand{\rP}{\mathrm{P}}
\newcommand{\ga}{\alpha}
\newcommand{\la}{\lambda}
\newtheorem{cor}{Corollary}[section]
\newproclaim{remark}{Remark}[section]
\newtheorem{lemma}{Lemma}[section]
\newtheorem{teo}{Theorem}[section]
\makeatother

\begin{document}
\begin{frontmatter}

\title{Strong limit of the extreme eigenvalues of a~symmetrized
auto-cross covariance matrix}
\runtitle{Strong limit of extreme eigenvalues}

\begin{aug}
\author[A]{\fnms{Chen}~\snm{Wang}\thanksref{M1}\ead[label=e1]{stawc@nus.edu.sg}},
\author[B]{\fnms{Baisuo}~\snm{Jin}\thanksref{M2,T1}\ead[label=e2]{jbs@ustc.edu.cn}},
\author[C]{\fnms{Z.~D.}~\snm{Bai}\corref{}\thanksref{T2,M3}\ead[label=e3]{baizd@nenu.edu.cn}},
\author[D]{\fnms{K.~Krishnan}~\snm{Nair}\thanksref{M4}\ead[label=e4]{kknair@stanford.edu}}
\and
\author[E]{\fnms{Matthew}~\snm{Harding}\thanksref{M5}\ead[label=e5]{matthew.harding@duke.edu}}
\runauthor{C. Wang et al.}
\affiliation{National University of Singapore\thanksmark{M1},
University of Science and Technology of China\thanksmark{M2},
Northeast Normal University\thanksmark{M3},\\
Stanford University\thanksmark{M4} and
Duke University\thanksmark{M5}}
\address[A]{C. Wang\\
Department of Statistics\\
\quad and Applied Probability\\
National University of Singapore\\
Singapore 117546\\
\printead{e1}}
\address[B]{B. Jin\\
Department of Statistics\\
\quad and Finance\\
University of Science\\
\quad and Technology of China\\
96, Jinzhai Road\\
Hefei 23 0026\\
P.~R. China\\
\printead{e2}}
\address[C]{Z.~D. Bai\\
KLASMOE and School of Math. and Stat.\\
Northeast Normal University\\
5268 Renmin Street, Changchun\\
Jilin Province 130024\\
P.~R.~China\\
\printead{e3}\\
\printead{u1}}
\address[D]{K.~K. Nair\\
Department of Civil  and\\
\quad Environmental Engineering\\
Stanford University\\
439 Panama Mall\\
Stanford, California 94305\\
USA\\
\printead{e4}}
\address[E]{M. Harding\\
Sanford School of Public Policy\\
Duke University\\
Durham North California 27708\\
USA\\
\printead{e5}}
\end{aug}
\thankstext{T1}{Supported by NSF China Young Scientist Grant 11101397.}
\thankstext{T2}{Supported by NSF China Grant 11171057, Program for
Changjiang Scholars and Innovative Research Team in University,
and the Fundamental Research Funds for the Central Universities.}

%
\received{\smonth{9} \syear{2014}}
%
\revised{\smonth{12} \syear{2014}}

%
\begin{abstract}
The auto-cross covariance matrix is defined as
\[
\mathbf{M}_n=\frac{1}{2T}\sum_{j=1}^{T}
\bigl(\mathbf{e}_{j}\mathbf{e}_{j+\tau}^{*}+\mathbf{e}
_{j+\tau}\mathbf{e}_{j}^{*}\bigr),
\]
where $\mathbf{e}_{j}$'s are $n$-dimensional vectors of independent standard
complex components
with a common mean 0, variance $\sigma^{2}$, and uniformly bounded
$2+\eta$th moments
and $\tau$ is the lag.
Jin et~al. [\textit{Ann. Appl. Probab.} \textbf{24} (2014) 1199--1225]
has proved that the LSD of
$\mathbf{M}_n$ exists
uniquely and nonrandomly, and independent of $\tau$ for all $\tau\ge
1$. And in addition
they gave an analytic expression of the LSD. As a continuation of Jin et~al. [\textit{Ann. Appl. Probab.} \textbf{24} (2014) \mbox{1199--1225}],
this paper proved that under the condition of uniformly bounded fourth
moments, in any
closed interval outside the support of the LSD, with probability 1
there will be no
eigenvalues of $\mathbf{M}_{n}$ for all large $n$. As a consequence of the
main theorem,
the limits of the largest and smallest eigenvalue of $\mathbf{M}_n$ are also
obtained.
\end{abstract}

%
\end{frontmatter}

\section{Introduction}\label{sec1}

For a $p\times p$ random Hermitian matrix $\bbA$ with eigenvalues $\la
_j,j=1,2,\ldots, p$, we define the empirical spectral distribution
(ESD) of $\bbA$ by
\[
F^{\bbA}(x)=\frac{1}{p}\sum_{j=1}^pI(
\la_j\leq x).
\]
The limit distribution $F$ of $\{F^{\bbA_n}\}$ for a given sequence of
random matrices $\{\bbA_n\}$ is called the limiting spectral
distribution (LSD).
Let $\{\ep_{it}\}$ be independent random variables with common mean 0
and variance 1. Define $\bbe_k=(\ep_{1k},\ldots,\ep_{nk})'$, $\bgma
_k=\frac{1}{\sqrt{2T}}\bbe_k$ and $\bbM_n(\tau)=\sum_{k=1}^T(\bgma_k\bgma_{k+\tau}^*+\bgma_{k+\tau}\bgma_k^*)$. Here,\vspace*{1pt}
$\tau\ge1$ is the number of lags.
Under the condition of bounded $2+\eta$th moments, \citet{Jinetal14} or under the weaker condition of second moments, \citet{BaiWan15}
derived the LSD of $\bbM_n(\tau)$, namely, $F^{\bbM_n(\tau
)}=:F_n\overset{w}{\rightarrow} F_{c}$ a.s. and $F_{c}$ has a density
function given by
%
\begin{eqnarray}
\label{Jin} \phi_c(x)=\frac{1}{2c\pi}\sqrt{
\frac{y_0^2}{1+y_0}- \biggl(\frac
{1-c}{\llvert  x\rrvert  }+\frac{1}{\sqrt{1+y_0}}
\biggr)^2},
\nonumber\\[-8pt]\\[-8pt]
\eqntext{-d(c)\le x\le d(c).}
\end{eqnarray}
Here, $c=\lim_{n\to\infty}c_n:=\lim_{n\to\infty}\frac{n}T$ and
$y_0$ is the largest real root of the equation
\[
y^3-\frac{(1-c)^2-x^2}{x^2}y^2-\frac{4}{x^2}y-
\frac{4}{x^2}=0
\]
and
\begin{eqnarray*}
d(c)=\cases{ \displaystyle\frac{(1-c)\sqrt{1+y_1}}{y_1-1}, &\quad$c\ne1$,
\vspace*{5pt}\cr
\displaystyle
\lim_{c\to1}\frac{(1-c)\sqrt{1+y_1}}{y_1-1}=\lim_{c\to1}\sqrt
{\frac{1+y_1}{y_1^3}}\sqrt{1+y_1}=2, &\quad$c=1$,}
\end{eqnarray*}
where $y_1$ is a real root of the equation:
\[
\bigl((1-c)^2-1 \bigr)y^3+y^2+y-1=0
\]
such that $y_1>1$ if $c<1$ and $y_1\in(0,1)$ if $c>1$.
Further, if $c>1$, then $F_{c}$ has a point mass $1-1/c$ at the origin.

The model of consideration comes from a high-dimensional dynamic
$k$-factor model with lag $q$, that is, $\bbR_{t}=\sum_{i=0}^{q}\bgL_{i}
\bbF_{t-i}+\bbe_{t},t=1,\ldots,T$. The factor $\bbF_{t-\tau}$
captures the structural part of the model at lag $\tau$, while $\bbe
_t$ corresponds to the noise component. Readers are referred to \citet{Jinetal14} for more details. An interesting problem to economists is
how to estimate $k$ and $q$. To solve this problem, for $\tau
=0,1,\ldots,$ define $\Phi_{n}(\tau)=\frac{1}{2T}\sum_{j=1}^{T}(\bbR_{j}\bbR_{j+\tau}^{*}+\bbR_{j+\tau}\bbR_{j}^{*})$.
Note that essentially, $\bbM_n(\tau)$ and $\Phi_{n}(\tau)$ are
symmetrized auto-cross covariance matrices at lag $\tau$ and
generalize the standard sample covariance matrices $\bbM_{n}(0)$ and
$\Phi_{n}(0)$, respectively.
The matrix $\bbM_{n}(0)$ has been intensively studied in the
literature and it is well known that the LSD has an MP law [\citet{MarPas67}].
Moreover, when $\tau=0$ and $\operatorname{Cov}(\bbF_t)=\Sigma_f$, the
population covariance matrix of $\bbR_{t}$ is a \emph{spiked
population model} [\citet{Joh01}, \citet{BaiSil06}, \citet{BaiYao08}].
In fact, under certain conditions, $k(q+1)$ can be estimated by
counting the number of eigenvalues of $\Phi(0)$ that are significantly
larger than $(1+\sqrt{c})^{2}$. What remains is to separate the
estimates of $k$ and $q$, which can be achieved using the LSD of $\bbM
_n=\bbM_{n}(\tau)$ for general $\tau\ge1$. A related work has been
found in \citet{LiWanYao} in which the number $k$ was detected by a
different symmetrized covariance matrix for a factor model without
lags. \citet{Jinetal14}\vadjust{\goodbreak} has proved that the LSD of $\bbM_n$ exists uniquely
and nonrandomly, and independent of $\tau$ for all $\tau\ge1$,
whose Stieltjes transform $m(z)$ satisfies the following equation:
\[
\bigl(1-c^2m^2(z)\bigr) \bigl(c+czm(z)-1
\bigr)^2=1,
\]
from which four roots are obtained, with $y_0$ defined as above:
\begin{eqnarray*}
m_1(z)&=&\frac{(\vfrac{1-c}{z}+\sqrt{1+y_0})+\sqrt{(\vfrac
{1-c}{z}-\sfrac{1}{\sqrt{1+y_0}})^2-\afrac{y_0^2}{1+y_0}}}{2c},
\\
m_2(z)&=&\frac{(\vfrac{1-c}{z}+\sqrt{1+y_0})-\sqrt{(\vfrac
{1-c}{z}-\sfrac{1}{\sqrt{1+y_0}})^2-\afrac{y_0^2}{1+y_0}}}{2c},
\\
m_3(z)&=&\frac{(\vfrac{1-c}{z}-\sqrt{1+y_0})+\sqrt{(\vfrac
{1-c}{z}+\sfrac{1}{\sqrt{1+y_0}})^2-\afrac{y_0^2}{1+y_0}}}{2c},
\\
m_4(z)&=&\frac{(\vfrac{1-c}{z}-\sqrt{1+y_0})-\sqrt{(\vfrac
{1-c}{z}+\sfrac{1}{\sqrt{1+y_0}})^2-\afrac{y_0^2}{1+y_0}}}{2c}.
\end{eqnarray*}
Here, as convention, we assume that the square root with a complex
number is the one whose imaginary part is positive and the Stieltjes
transform for a function of bounded variation $G$ is defined as
\[
m_G(z)=\int\frac{1}{x-z}\,dG(x)\qquad\mbox{for complex } \Im(z)>0.
\]
However, the number of eigenvalues of
$\Phi_{n}(\tau)$ that lie outside the support of the LSD of $\bbM_n$
at lags $1 \leq\tau\leq q $ is different from that at lags $\tau>q$.
Thus, the estimates of $k$ and $q$ can be separated by counting the
number of eigenvalues of
$\Phi_{n}(\tau)$ that lie outside the support of the LSD of $\bbM_n$
from $\tau=0,1,2,\ldots,q,q+1,\ldots.$

It is worth noting that for the above method to work, one should expect
no eigenvalues outside the support of the LSD of $\bbM_n$ so that
if an eigenvalue of $\Phi_{n}(\tau)$ goes out of the support of the
LSD of $\bbM_n$, it must come from the signal part. As a continuation
of \citet{Jinetal14}, this paper establishes limits of the largest and
smallest eigenvalues of $\bbM_n$, after showing that no eigenvalues
exist outside the support of the LSD of $\bbM_n$, along the similar
lines as in \citet{BaiSil98}.

In \citet{BaiSil98}, the authors considered the separation
problem of the general sample covariance matrices. Later, \citet{PauSil09} extended the result to a more general class of
matrices taking the form of $\frac{1}{n}\bbA_n^{1/2}\bbX_n\bbB
_n\bbX_n^*\bbA_n^{1/2}$ and \citet{BaiSil12} established
the result for the information-plus-noise matrices.

Compared with \citet{BaiSil98}, the model we considered here
is more complicated and some new techniques are employed. Besides the
recursive method to solve a disturbed difference equation as in \citet{Jinetal14}, a relationship between the convergence rates of polynomial
coefficients and those of the roots is established and applied.
Moreover, the results in this paper will pave the way for establishing
other results such as limit theorems for sample eigenvalues of the
spiked model.
The main results can now be stated.
%
\begin{teo} \label{teo1}Assume:
\begin{longlist}[(a)]
\item[(a)] $\tau\ge1$ is a fixed integer.
\item[(b)] $\bbe_{k}=(\ep_{1k},\ldots,\ep_{nk})'$, $k=1,2,\ldots,T+\tau$, are $n$-vectors of independent standard complex components
with $\sup_{i, t}\rE\llvert  \ep_{it}\rrvert  ^4\le M$ for some $M>0$.
\item[(c)] There exist $K>0$ and a random variable $X$ with finite
fourth-order moment such that, for any $x>0$, for all $n,T$
%
\begin{equation}
\label{eta} \frac{1}{nT}\sum_{i=1}^n
\sum_{t=1}^{T+\tau}\rP\bigl(\llvert
\ep_{it}\rrvert > x\bigr)\le K\rP\bigl(\llvert X\rrvert >x\bigr).
\end{equation}
\item[(d)] $\bbM_n=\sum_{k=1}^{T}(\bgma_{k}\bgma_{k+\tau}^{*}+\bgma
_{k+\tau}\bgma_{k}^{*})$, where $\bgma_{k}=\frac{1}{\sqrt{2T}}\bbe_{k}$.
\item[(e)] $c_n\equiv n/T\rightarrow c\in(0,1)\cup(1,\infty)$ as
$n\to\infty$.
\item[(f)] The interval $[a,b]$ lies outside the support of $F_{c}$.
\end{longlist}
Then $\rP (\mbox{no eigenvalues of $\bbM_n$ appear in $[a,b]$
for all large n} )=1$.
\end{teo}
By definition of $\bbe_k$ and the convergence of the largest
eigenvalue of the sample covariance matrix [\citet{YinBaiKri88}], we have, for any $\delta>0$ and all large $n$,
%
\begin{eqnarray}\label{upbd}
\llVert \bbM_n\rrVert &\le&\frac{1}{2T} \bigl(
\bigl\llVert \mathbf{E}\mathbf{E}_{\tau
}^*\bigr\rrVert +\bigl\llVert
\mathbf{E}_{\tau}\mathbf{E}^*\bigr\rrVert \bigr) \nonumber
\\
&\le&\frac{1}{T}s_{\max}(
\mathbf{E})s_{\max}(\mathbf{E}_{\tau}) =s_{\max}
\biggl(\frac{\mathbf{E}}{\sqrt{T}} \biggr)s_{\max} \biggl(\frac{\mathbf{E}_{\tau}}{\sqrt{T}}
\biggr)
\\
&\le&(1+\sqrt c)^2+\delta\qquad\mbox{a.s.}\nonumber
\end{eqnarray}
Here, $\mathbf{E}=(\bbe_1,\ldots,\bbe_T)$, $\mathbf{E}_{\tau
}=(\bbe_{1+\tau},\ldots,\bbe_{T+\tau})$ and $s_{\max}(\bbA)$
denotes the largest singular value of a matrix $\bbA$. This, together
with Theorem~\ref{teo1}, implies the following result.
%
\begin{teo}\label{teo2} Assuming conditions \textup{(a)--(e)} in Theorem~\ref{teo1} hold, we have
\begin{eqnarray*}
&&\lim_{n\to\infty}\lambda_{\min}(\bbM_n)=-d(c)
\qquad\mbox {a.s.}\quad\mbox{and}\quad\lim_{n\to\infty}
\lambda_{\max}(\bbM _n)=d(c)\qquad\mbox{a.s.}
\end{eqnarray*}
Here, $-d(c)$ and $d(c)$ are the left and right boundary points of the
support of the LSD of $\bbM_n$, as defined in (\ref{Jin}).
\end{teo}
\begin{pf} When $c\in(0,1)\cup(1,\infty)$, let $\ep>0$ be given and
consider the interval $[d(c)+\ep, b]$ with $b>(1+\sqrt c)^2+\delta$
for some $\delta>0$. By (\ref{upbd}), with probability one, there is
no eigenvalue in the interval $(b, \infty)$. This, together with
Theorem~\ref{teo1}, implies that with probability one, there is no
eigenvalue in the interval $[d(c)+\ep, \infty)$. Therefore, we have
\begin{eqnarray*}
&&\limsup_{n\to\infty} \lambda_{\max}(\bbM_n)\le
d(c)+\ep\qquad \mbox{a.s.}
\end{eqnarray*}

Next, we claim that, for all large $n$, there exists at least one
eigenvalue in $[d(c)-\ep, d(c)]$. Otherwise, we have
$F_n(d(c))-F_n(d(c)-\ep)=0$ for infinitely many $n$, which contradicts
the fact that $F_n\to F_c$, or equivalently that
$F_c(d(c))-F_c(d(c)-\ep)>0$. Hence, our claim is proved. Therefore, we have
\begin{eqnarray*}
&&\liminf_{n\to\infty} \lambda_{\max}(\bbM_n)\ge
d(c)-\ep\qquad \mbox{a.s.}
\end{eqnarray*}
Now, let $\ep\to0$, and we then have $\lim_{n\to\infty} \lambda
_{\max}(\bbM_n)=d(c)$, a.s. By symmetry, $\lim_{n\to\infty}
\lambda_{\min}(\bbM_n)=-d(c)$, a.s. This completes the proof of the theorem.
\end{pf}

One can extend Theorem~\ref{teo2} to the case $c=1$ as follows.
%
\begin{teo}\label{teo3} When $c=1$, Theorem~\ref{teo2} still holds,
that is,
\[
\lim_{n\to\infty}\lambda_{\min}(\bbM_n)=-d(1)=-2
\qquad\mbox {a.s.}
\]
and
\[
\lim_{n\to\infty}
\lambda_{\max}(\bbM _n)=d(1)=2\qquad\mbox{a.s.}
\]
\end{teo}
\begin{pf}
To prove this theorem, we need to enlarge the matrix $\bbM_n$ with a
larger dimension. To this end, denote $\bbM_n=\bbM_{n,T}=\bbM
_{n,T(n)}$. Fix $T$, we show that $\lambda_{\max}(\bbM_{n,T})$ is
nondecreasing and $\lambda_{\min}(\bbM_{n,T})$ is nonincreasing in
$n$, or more precisely,
$\lambda_{\max}(\bbM_{n,T(n)})\le\lambda_{\max}(\bbM
_{n+1,T(n)})$ and $\lambda_{\min}(\bbM_{n,T(n)})\ge\lambda_{\min
}(\bbM_{n+1,T(n)})$.

To prove these relations, we will employ the interlacing theorem (Lemma
\ref{interlacing}) by showing that
$\bbM_{n,T(n)}$ is a major sub-matrix of $\bbM_{n+1,T(n)}$. Rewrite
\begin{eqnarray*}
\bbM_{n,T(n)}&=&\sum_{k=1}^{T(n)}
\bigl(\bgma_{k}\bgma_{k+\tau
}^{*}+
\bgma_{k+\tau}\bgma_{k}^{*}\bigr) =\sum
_{k=1}^{T(n)}\bigl(\bgma_{k,n}
\bgma_{k+\tau,n}^{*}+\bgma_{k+\tau,n}\bgma_{k,n}^{*}
\bigr).
\end{eqnarray*}
By introducing, $x_{t,n+1}=\frac{1}{\sqrt{2T(n)}}\ep_{(n+1)t}$, we obtain
{\fontsize{10}{11}\selectfont{\begin{eqnarray*}
\hspace*{-3pt}&&\bbM_{n+1,T(n)}
\\
\hspace*{-3pt}&&\qquad =\sum_{k=1}^{T(n)}\bigl(
\bgma_{k,n+1}\bgma_{k+\tau,n+1}^{*}+\bgma _{k+\tau,n+1}
\bgma_{k,n+1}^{*}\bigr)
\\
\hspace*{-3pt}&&\qquad =\sum_{k=1}^{T(n)} \Biggl[\pmatrix{
\bgma_{k,n}
\cr
x_{k,n+1}}\bigl(\bgma_{k+\tau,n}^{*},x_{k+\tau,n+1}^{*}
\bigr)+ \pmatrix{ \bgma_{k+\tau,n}
\cr
x_{k+\tau,n+1}}\bigl(
\bgma_{k,n}^{*},x_{k,n+1}^{*}\bigr) \Biggr]
\\
\hspace*{-3pt}&&\qquad = \pmatrix{ \displaystyle\sum_{k=1}^{T(n)}\bigl(
\bgma_{k,n}\bgma_{k+\tau,n}^{*}+\bgma_{k+\tau,n}
\bgma_{k,n}^{*}\bigr) & \displaystyle\sum_{k=1}^{T(n)}
\bigl(\bgma_{k,n}x_{k+\tau,n+1}^{*}+\bgma_{k+\tau,n}x_{k,n+1}^{*}
\bigr)
\vspace*{5pt}\cr
\displaystyle\sum_{k=1}^{T(n)}
\bigl(x_{k,n+1}\bgma_{k+\tau,n}^{*}+x_{k+\tau,n+1}\bgma
_{k,n}^{*}\bigr) & \displaystyle\sum_{k=1}^{T(n)}
\bigl(x_{k,n+1}x_{k+\tau,n+1}^{*}+x_{k+\tau,n+1}x_{k,n+1}^{*}
\bigr)}
\\[5pt]
\hspace*{-3pt}&&\qquad =\pmatrix{ \bbM_{n,T(n)} & \displaystyle\sum_{k=1}^{T(n)}
\bigl(\bgma_{k,n}x_{k+\tau,n+1}^{*}+\bgma_{k+\tau,n}x_{k,n+1}^{*}
\bigr)
\vspace*{5pt}\cr
\displaystyle \sum_{k=1}^{T(n)}
\bigl(x_{k,n+1}\bgma_{k+\tau,n}^{*}+x_{k+\tau,n+1}\bgma
_{k,n}^{*}\bigr) & \displaystyle\sum_{k=1}^{T(n)}
\bigl(x_{k,n+1}x_{k+\tau,n+1}^{*}+x_{k+\tau,n+1}x_{k,n+1}^{*}
\bigr)}.
\end{eqnarray*}}}%
By Lemma~\ref{interlacing}, we have $\lambda_{\max}(\bbM
_{n+1,T(n)})\ge\lambda_{\max}(\bbM_{n,T(n)})$. By symmetry, we also
have $\lambda_{\min}(\bbM_{n+1,T(n)})\le\lambda_{\min}(\bbM_{n,T(n)})$.
This together with Theorem~\ref{teo2} implies that for any $\ep>0$,
we have a.s.
\begin{eqnarray*}
&&\mathop{\limsup_{n\to\infty}}_{n/T(n)\to1}\lambda_{\max}(
\bbM _{n,T(n)})\le\mathop{\lim_{n\to\infty}}_{n/T(n)\to1}
\lambda_{\max
}(\bbM_{[(1+\ep)n],T(n)})=d(1+\ep).
\end{eqnarray*}
Note that $d(c)$ is continuous in $c$. By letting $\ep\to0$, we have
a.s.
\[
\mathop{\limsup_{n\to\infty}}_{n/T(n)\to1}\lambda
_{\max}(\bbM_{n,T(n)})\le d(1)=2.
\]
Since the LSD of $\bbM_n$ exists
with right support boundary $d(1)=2$, we have proved that
\[
\mathop{\lim_{n\to\infty}}_{n/T(n)\to1}\lambda_{\max}(\bbM
_{n,T(n)})=2.
\]
By symmetry, we have a.s.
$\lim_{n\to\infty, n/T(n)\to1}\lambda_{\min
}(\bbM_{n,T(n)})=-d(1)=-2$. The proof of the theorem is complete.
\end{pf}

As an immediate consequence of Theorem~\ref{teo3}, Corollary~\ref
{teo4} complements Theorem~\ref{teo1} for $c=1$.

\begin{cor}\label{teo4} Theorem~\ref{teo1} still holds when $c=1$.
\end{cor}

\begin{figure}[t]

\includegraphics{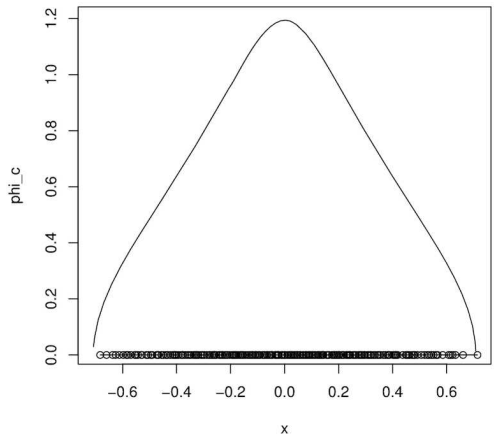}

\caption{Density function $\phi_c(x)$ of $F_c$ and distribution of
sample eigenvalues with $\tau=1, c=0.2$ $(n=200,T=1000)$.}
\label{fig1}
\end{figure}

\begin{figure}

\includegraphics{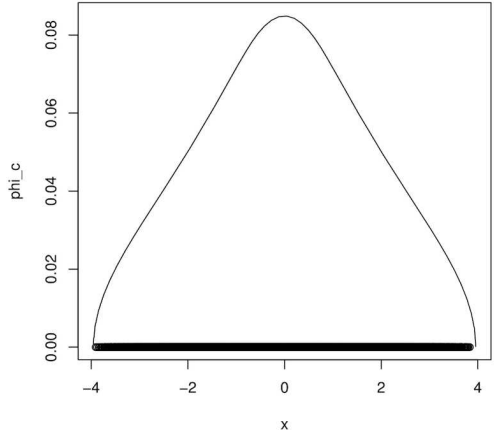}

\caption{Density function $\phi_c(x)$ of $F_c$ and distribution of
sample eigenvalues with $\tau=1, c=2.5$ $(n=2500,T=1000)$. Note that the
area under the density function curve is $1/c$.}\label{fig2}
\end{figure}

Figures~\ref{fig1} and \ref{fig2} display the density functions $\phi_c(x)$ and the
distributions of sample eigenvalues with $\tau=1, c=0.2$
$(n=200,T=1000)$ and $c=2.5$ $(n=2500,T=1000)$, respectively.
%

We will now focus on proving Theorem~\ref{teo1}. As in \citet{Jinetal14}, we denote the Stieltjes transform of $\bbM_n$ as $m_n(z)=\frac{1}n\rtr(\bbM_n-z\bbI_n)^{-1}$ where, and throughout the paper,
$z=u+iv_n$, $v_n>0$, and let $m_n^0(z)$ be the Stieltjes transform of
$\phi_{c_n}$ with limiting ratio
of $c_n=n/T$. Using the truncation technique employed in Section~3 of
\citet{BaiSil98}, we further assume that the $\ep_{ij}$'s
satisfy the conditions that
%
\begin{equation}
\label{trun} \llvert \ep_{ij}\rrvert \le C,\qquad \rE
\ep_{ij}= 0,\qquad \rE\llvert \ep_{ij}\rrvert
^2=1,\qquad \rE\llvert \ep_{ij}\rrvert ^4<M
\end{equation}
for some $C, M>0$. More detailed justifications are provided in the
\hyperref[append]{Appendix}.

The rest of the paper is structured as follows. Section~\ref{sec2} contains some
lemmas of known results. Section~\ref{s3} provides some technical lemmas.
Convergence rates of $\llVert  F_n-F_{c_n}\rrVert  $ and $m_{n}(z)-m_n^0(z)$ are
obtained in Sections~\ref{s4} and~\ref{s5}, respectively. Section~\ref{sec6} concludes the
proof of Theorem~\ref{teo1}. Justifications of variable truncation,
centralization and rescaling and proofs of lemmas presented in
Section~\ref{s3} are given in the \hyperref[append]{Appendix}.

\section{Mathematical tools}\label{sec2}

In this section, we provide some known results.

\begin{lemma}[{[\citet{Bur73}]}]\label{Burkholder} Let $\{X_k\}$ be a
complex martingale difference sequence with respect to the increasing
$\sigma$-fields $\{\mathcal{F}_n\}$. Then, for $p\ge2$, we have
\begin{eqnarray*}
&&\erE\Bigl\llvert \sum X_k\Bigr\rrvert ^p\le
K_p \Bigl(\erE \Bigl(\sum\erE\bigl(\llvert X_k
\rrvert ^2| \mathcal{F}_{k-1}\bigr)
\Bigr)^{p/2}+\erE\sum\llvert X_k\rrvert
^p \Bigr).
\end{eqnarray*}
\end{lemma}

\begin{lemma}[{[\citet{Bur73}]}]\label{Burkholder1} Let $\{X_k\}$ be
as above. Then, for $p\ge2$, we have
\begin{eqnarray*}
&&\erE\Bigl\llvert \sum X_k\Bigr\rrvert ^p\le
K_p\erE \Bigl(\sum\llvert X_k\rrvert
^2 \Bigr)^{p/2}.
\end{eqnarray*}
\end{lemma}

\begin{lemma}[{[Theorem A.43 of \citet{BaiSil10}]}]\label{A43} Let $\bbA$ and $\bbB$ be two $n\times n$ Hermitian matrices. Then
\begin{eqnarray*}
&&\bigl\llVert F^\bbA-F^\bbB\bigr\rrVert \le
\frac{1}{n}\mathrm{rank}(\bbA-\bbB),
\end{eqnarray*}
where $F^\bbA$ is the empirical spectral distribution of $\bbA$ and
$\llVert  f\rrVert  =\sup_x\llvert  f(x)\rrvert  $.
\end{lemma}

\begin{lemma}[{[\citet{Bai93} or Corollary B.15 of \citet{BaiSil10}]}]\label{B14}  Let~$F$ be a distribution function and let $G$ be a function of
bounded variation satisfying $\int\llvert  F(x)-G(x)\rrvert  \,dx <\infty$. Denote
their Stieltjes transforms by $f(z)$ and $g(z)$, respectively. Assume
that for some constant $B > 0$, $F([-B,B]) = 1$ and $\llvert  G\rrvert  ((-\infty,-B))
= \llvert  G\rrvert  ((B,\infty)) = 0$, where $\llvert  G\rrvert  ((a, b))$ denotes the total
variation of the signed measure $G$ on the interval $(a, b)$. Then we have
\begin{eqnarray*}
\llVert F-G\rrVert &:=&\sup_x\bigl\llvert F(x)-G(x)\bigr
\rrvert
\\
&\le& \frac{1}{\pi(1-\kappa)(2\gamma- 1)}
\\
&&{}\times \biggl[\int_{-A}^{A}
\bigl\llvert f(z)-g(z)\bigr\rrvert \,du+ v^{-1} \sup_x
\int_{\llvert  y\rrvert  \le2va} \bigl\llvert G(x + y)-G(x)\bigr\rrvert \,dy \biggr],
\end{eqnarray*}
where $z=u+iv$, $v>0$, $a$ and $\gamma$ are positive constants such
that $\gamma=\frac{1}\pi\int_{\llvert  u\rrvert  <a}\frac{1}{u^2+1}\,du>\frac{1}2$. $A$
is a positive constant such that $A>B$ and $\kappa=\frac{4B}{\pi
(A-B)(2\gamma-1)}<1$.
\end{lemma}

\begin{lemma}[{[Lemma B.26 of \citet{BaiSil10}]}]\label{B26}
Let $\bbA=(a_{ij})$ be an $n\times n$ nonrandom matrix and $\bbX=
(x_1,\ldots,x_n)'$ be a random vector of independent entries. Assume
that $\erE x_i=0$,
$\erE\llvert  x_i\rrvert  ^2=1$, and $\erE\llvert  x_j\rrvert  ^{\ell}\le v_{\ell}$.
Then, for any $p\ge1$,
\begin{eqnarray*}
&&\erE\bigl\llvert \bbX^*\bbA\bbX-\rtr \bbA\bigr\rrvert ^p\le
C_p \bigl( \bigl(v_4\rtr\bigl(\bbA\bbA^*\bigr)
\bigr)^{p/2}+v_{2p}\rtr\bigl(\bbA \bbA^*\bigr)^{p/2}
\bigr),
\end{eqnarray*}
where $C_p$ is a constant depending on $p$ only.
\end{lemma}

\begin{lemma} [{[The interlacing theorem, \citet{RaoRao98}]}]\label{interlacing} If $\bbC$ is an $(n-1) \times(n-1)$ major sub-matrix of the
$n \times n$ Hermitian matrix $\bbA$, then $\lambda_1(\bbA)\ge
\lambda_1(\bbC)\ge\lambda_2(\bbA)\ge\cdots\ge\lambda
_{n-1}(\bbC)\ge\lambda_n(\bbA)$. Here
$\lambda_i(\bbA)$ denotes the $i$th largest eigenvalue of the
Hermitian matrix $\bbA$.
\end{lemma}

\section{Some technical lemmas}\label{s3}\label{sec25}

Before proceeding, some technical lemmas are presented with proofs
postponed in the \hyperref[append]{Appendix}. The first three are about the convergence
rates of roots of a polynomial.

\begin{lemma}\label{P1} Let $\{r_n\}$ be a sequence of positive real
numbers converging to $0$ and $m$ be a fixed positive integer,
independent of $n$. Let $B(x_0,r_n)$ denote the open ball centered at
$x_0$ with radius $r_n$. Given $m$ points $x_1,\ldots,x_m$ in
$B(x_0,r_n)$, one can find $x\in B(x_0,r_n)$ and $d>0$ such that $\min_{i\in\{1,\ldots,m\}}\llvert  x-x_i\rrvert  \ge \,dr_n$.
\end{lemma}

\begin{lemma} \label{P2}For each $n\in\mathbb{N}$, let
$P_n(x)=x^k+a_{n,k-1}x^{k-1}+\cdots+a_{n,1}x+a_{n,0}$ be a polynomial
of degree k, with roots $x_{n1},\ldots,x_{nk}$. Moreover, for
$i=0,1,\ldots,k-1$, $\lim_{n\to\infty}a_{n,i}=a_i$. Let
$P(x)=x^k+a_{k-1}x^{k-1}+\cdots+a_{1}x+a_{0}$. Suppose $P(x)$ has
distinct roots $x_1,\ldots,x_m$, and each $x_j$ has multiplicity $\ell
_j$ with $\sum_{j=1}^{m}\ell_j=k$. Then for $n$ large enough, for
each $j\in{\{1,\ldots,m\}}$, there are exactly $\ell_j$ $x_{ni}$'s
in $B(x_j,r_n^{1/\ell_j})$, where $r_n=\max_{i\in\{0,1,\ldots,k-1\}
}\llvert  a_{n,i}-a_i\rrvert  $.
\end{lemma}

\begin{lemma}\label{P3} For each $n\in\mathbb{N}$, let
$P_n(x)=x^k+a_{n,k-1}x^{k-1}+\cdots+a_{n,1}x+a_{n,0}$ and
$Q_n(y)=y^k+b_{n,k-1}y^{k-1}+\cdots+b_{n,1}y+b_{n,0}$ be two
polynomials of degree k, with roots $x_{n1},\ldots,x_{nk}$ and
$y_{n1},\ldots,y_{nk}$, respectively. Moreover, for $i=0,1,\ldots,k-1$, $\lim_{n\to\infty}b_{n,i}=\lim_{n\to\infty}a_{n,i}=a_i$. Let
$P(x)=x^k+a_{k-1}x^{k-1}+\cdots+a_{1}x+a_{0}$. Suppose $P(x)$ has
distinct roots $x_1,\ldots,x_m$, and each $x_j$ has the multiplicity
$\ell_j$ with $\sum_{j=1}^{m}\ell_j=k$. Then\vspace*{-1pt} for $n$ large enough,
for each $j\in{\{1,\ldots,m\}}$, for any $x_{ni}\in B(x_j,r_n^{1/\ell
_j})$, there exists at least one $y_{nl}$ such that $\llvert  x_{ni}-y_{nl}\rrvert  \le
d\widetilde r_n^{1/\ell_j}$ for some $d>0$. Here, $r_n=\max_{i\in\{
0,1,\ldots,k-1\}}\llvert  a_{n,i}-a_i\rrvert  $ and $\widetilde r_n=\max_{i\in\{
0,1,\ldots,k-1\}}\llvert  a_{n,i}-b_{n,i}\rrvert  $.
\end{lemma}

To establish the following lemmas, we need some notation: let
$z=u+iv_n$, where $u\in[-A,A]$ and $v_n\in[n^{-1/52},n^{-1/212}]$ and
$A>0$ is a large constant.\vadjust{\goodbreak} Define
\begin{eqnarray*}
\bbA&=&\bbM_n-z\bbI_n,
\\
\bbA_k&=&\bbM_{n,k}-z\bbI_{n}=\bbA-
\bgma_{k}(\bgma_{k-\tau
}+\bgma_{k+\tau})^*-(
\bgma_{k-\tau}+\bgma_{k+\tau})\bgma_k^*,
\\
&\vdots&
\\
\bbA_{k,\ldots,k+s\tau}&=&\bbA-\sum_{t=0}^s
\bigl[\bgma_{k+t\tau
}(\bgma_{k+(t-1)\tau}+\bgma_{k+(t+1)\tau})^*
\\
&&\hspace*{37pt}{} +(
\bgma_{k+(t-1)\tau
}+\bgma_{k+(t+1)\tau})\bgma_{k+t\tau}^* \bigr],
\end{eqnarray*}
with the convention that $\bgma_{l}=0$ for $l\le0$ or $l >T +\tau$.

The following lemma will be frequently used.

\begin{lemma}\label{rs}Let $r,s$ be fixed positive integers. For $l\ne
k$, we have
\begin{eqnarray*}
&&\rE\bigl\llvert \bgma_{l}^*\bbA_k^{-s}
\bgma_k\bigr\rrvert ^{2r}\le\frac{K}{T^rv_n^{2rs}}
\end{eqnarray*}
for some $K>0$.
\end{lemma}

Define $a_n=\frac{c_n\rE m_n}{2}$ and let $x_{n1}, x_{n0}$ be two
roots of the equation $x^2=x-a_n^2$ with $\llvert  x_{n1}\rrvert  >\llvert  x_{n0}\rrvert  $. Some
properties regarding $x_{n1}$ and $x_{n0}$ are stated in the next lemma.

In the following, if a lemma contains two sets of results
simultaneously, then the results labelled by ``a'' hold for all
$z=u+iv_n$, and $u$ lies in a bounded interval $[-A,A] \subseteq
\mathbb{R}$, whereas results labelled by ``b'' hold for all $z=u+iv_n$
with $u\in[a,b]$ and are obtained under the additional condition that
$\rP(\llVert  F_{n}-F_{c_n}\rrVert  \ge n^{-1/104})=o(n^{-t})$ for any fixed $t>0$,
where $[a,b]$ is defined in Theorem~\ref{teo1}. Results ``a'' will be
used to establish a preliminary convergence rate of the ESD of $\bbM
_n$ in Section~\ref{sec3} and the results ``b'' will be applied to the
refinement of the convergence rate when $u\in[a, b]$ in Section~\ref
{sec4}. If a lemma contains only one set of results, the results will
be established for all $u\in[a,b]$ and under the additional assumption
that $\rP(\llVert  F_{n}-F_{c_n}\rrVert  \ge n^{-1/104})=o(n^{-t})$.

\begin{lemma}\label{lkj} When $u\in[a,b]$, let $\lambda_{kj}$ denote
the $j$th largest eigenvalue of $\bbM_n-\bgma_{k}(\bgma_{k+\tau
}+\bgma_{k-\tau})^*-(\bgma_{k+\tau}+\bgma_{k-\tau})\bgma_{k}^*$,
for $\Im(z)\ge n^{-\gd}$ with $\gd=1/106$, we have, for any $t>0$
\begin{eqnarray*}
&&\rP \biggl(\frac{1}{2T}\sum\frac{1}{\llvert  \lambda_{kj}-z\rrvert  ^2}> K \biggr)=o
\bigl(n^{-t}\bigr)
\end{eqnarray*}
for some $K>0$.
\end{lemma}

\begin{remark}\label{lkj1}When $u\in[a,b]$, with similar proofs, for
$\Im(z)\ge n^{-\gd}$ with $\gd=1/53$, we have, for any $t>0$,
\begin{eqnarray*}
&&\rP \biggl(\frac{1}{2T}\bigl\llvert \operatorname{tr}\bbA_k^{-1}
\bigr\rrvert > K \biggr)\le\rP \biggl(\frac{1}{2T}\sum
\frac{1}{\llvert  \lambda_{kj}-z\rrvert  }> K \biggr)=o\bigl(n^{-t}\bigr)
\end{eqnarray*}
and when $\Im(z)\ge n^{-\gd}$ with $\gd=1/212$,
\begin{eqnarray*}
\rP \biggl(\frac{1}{2T}\bigl\llvert \operatorname{tr}\bbA_k^{-4}
\bigr\rrvert > K \biggr)\le\rP \biggl(\frac{1}{2T}\sum
\frac{1}{\llvert  \lambda_{kj}-z\rrvert  ^4}> K \biggr)=o\bigl(n^{-t}\bigr)
\end{eqnarray*}
for some $K>0$.
\end{remark}

\begin{remark}\label{lkj2}When $u\in[a,b]$, and $\lambda_{kj}$'s are
eigenvalues of $\bbM_{n,k}=\bbM_n-\bgma_{k}(\bgma_{k+\tau}+\bgma
_{k-\tau})^*-(\bgma_{k+\tau}+\bgma_{k-\tau})\bgma_{k}^*$, for
$\Im(z)\ge n^{-\gd}$ with $\gd=1/212$, with a similar proof, we have
\begin{eqnarray*}
&&\rP \biggl(\frac{1}{2T}\sum\frac{1}{\llvert  \lambda_{kj}-z\rrvert  ^2}> K \biggr)=o
\bigl(n^{-t}\bigr)
\end{eqnarray*}
for some $K>0$.
\end{remark}

\begin{lemma}\label{x1} With $x_{n1}$ and $x_{n0}$ defined as above,
for any $v_n\ge n^{-1/52}$, we have:
\begin{enumerate}[(iii)]
\item[(i)] There exists some $\eta>0$ such that for all large $n$:
\begin{longlist}[(a)]
\item[(a)]
$\sup_{u\in[-A, A], \Im(z)=v_n}\llvert  \frac
{x_{n0}(z)}{x_{n1}(z)}\rrvert  <1-\eta v_n^3$.\vspace*{2pt}
\item[(b)]
$\sup_{u\in[a,b], \Im(z)=v_n}\llvert  \frac
{x_{n0}(z)}{x_{n1}(z)}\rrvert  <1-\eta$.
\end{longlist}
\item[(ii)]\mbox{}
\begin{longlist}[(a)]
\item[(a)] When $u\in[-A, A]$, we have $\llvert  x_{n1}\rrvert  \ge\frac{1}2$ and
$\llvert  x_{n1}\rrvert  \le Kv_n^{-1}$ for some constant $K$.
\item[(b)] When $u\in[a,b]$, we have $\llvert  x_{n1}\rrvert  \ge\frac{1}2$ and
$\llvert  x_{n1}\rrvert  \le K$ for some constant~$K$.
\end{longlist}
\item[(iii)]\mbox{}
\begin{longlist}[(a)]
\item[(a)] When $u\in[-A, A]$, we have $\llvert  x_{n1}-x_{n0}\rrvert  \ge\eta v_n$
for some constant \mbox{$\eta>0$.}
\item[(b)] When $u\in[a,b]$, we have $\llvert  x_{n1}-x_{n0}\rrvert  \ge\eta$ for
some constant $\eta>0$.
\end{longlist}
\item[(iv)]\mbox{}
\begin{longlist}[(a)]
\item[(a)] When $u\in[-A, A]$, we have $\frac
{\llvert  x_{n1}\rrvert  }{\llvert  x_{n1}-x_{n0}\rrvert  }\le Kv_n^{-1}$ for some constant $K$.\vspace*{2pt}
\item[(b)] When $u\in[a,b]$, we have $\frac
{\llvert  x_{n1}\rrvert  }{\llvert  x_{n1}-x_{n0}\rrvert  }\le K$ for some constant $K$.
\end{longlist}
\item[(v)]
When $u\in[a,b]$, we have $\llvert  a_n\rrvert  <\frac{1}2-\eta$ for some
constant $\eta>0$.
\end{enumerate}
\end{lemma}

\begin{lemma}\label{Ek+tau}\label{kkbar0}\label{ep41}\label{ep43}
For any $v_n\ge n^{-1/52}$ and $t>0$:
\begin{longlist}[(b1)]
\item[(a)] for any $u\in[-A,A]$ and $k\le T-v_{n}^{-4}$, we have
\begin{eqnarray*}
\rP\biggl(\biggl\llvert \bgma_{k+\tau}^*\bbA_k^{-1}
\bgma_{k+\tau}-\frac{c_n\rE
m_n}{2x_{n1}}\biggr\rrvert \ge v_n^6
\biggr)=o\bigl(n^{-t}\bigr)
\end{eqnarray*}
and for any $k\ge v_n^{-4}$,
\begin{eqnarray*}
\rP\biggl(\biggl\llvert \bgma_{k-\tau}^*\bbA_k^{-1}
\bgma_{k-\tau}-\frac{c_n\rE
m_n}{2x_{n1}}\biggr\rrvert \ge v_n^6
\biggr)=o\bigl(n^{-t}\bigr),
\end{eqnarray*}
\item[(b1)]  for any $u\in[a,b]$, there is a constant
$\eta\in(0,\frac{1}2)$ such that
$\rP(\llvert  \bgma_{k+\tau}^*\bbA_k^{-1}\*\bgma_{k+\tau}\rrvert  \ge1-\eta)=o(n^{-t})$,
\item[(b2)]  for any $u\in[a,b]$, when $k\le T-\log
^{2}n$, we have
$\llvert  \rE\bgma_{k+\tau}^*\bbA^{-1}_{k}\bgma_{k+\tau}-\frac
{a_n}{x_{n1}}\rrvert  =o(1/(nv_n))$,
and when $k\ge\log^{2}n$, we have
$\llvert  \rE\bgma_{k-\tau}^*\bbA^{-1}_{k}\bgma_{k-\tau}-\frac
{a_n}{x_{n1}}\rrvert  =o(1/(nv_n))$,
\item[(b3)]  for any $u\in[a,b]$, when $k\le T-\log
^{2}n$, we have
$\rE\llvert  \bgma_{k+\tau}^*\bbA^{-1}_{k}\bgma_{k+\tau}-\frac
{a_n}{x_{n1}}\rrvert  ^2=o(1/(nv_n))$,\vspace*{1pt}
and when $k\ge\log^{2}n$, we have
$\rE\llvert  \bgma_{k-\tau}^*\bbA^{-1}_{k}\bgma_{k-\tau}-\frac
{a_n}{x_{n1}}\rrvert  ^2=o(1/(nv_n))$.
\end{longlist}
\end{lemma}



\begin{lemma}\label{cross}\label{ep42}\label{ep44}
For any $v_n\ge n^{-1/52}$ and $t>0$:
\begin{longlist}[(b2)]
\item[(a)] for any $u\in[-A, A]$, we have
\begin{eqnarray*}
\rP\bigl(\bigl\llvert \bgma_{k-\tau}^*\bbA_k^{-1}
\bgma_{k+\tau}\bigr\rrvert >v_n^6\bigr)=o
\bigl(n^{-t}\bigr);
\end{eqnarray*}
\item[(b1)] for any $u\in[a,b]$, we have
$\llvert  \rE\bgma_{k-\tau}^*\bbA^{-1}_{k}\bgma_{k+\tau}\rrvert  =o(1/(nv_n))$;
\item[(b2)] for any $u\in[a,b]$, we have
$\rE\llvert  \bgma_{k-\tau}^*\bbA^{-1}_{k}\bgma_{k+\tau}\rrvert  ^2=o(1/(nv_n))$.
\end{longlist}
\end{lemma}


\begin{lemma}\label{kkbar}For any $v_n\ge n^{-1/212}$, $u\in[a,b]$
and $t>0$, there exists a constant $K>0$ such that
\begin{eqnarray*}
\rP\bigl(\bigl\llvert \bgma_{k+\tau}^*\bbA_k^{-1}
\bigl(\bbA^*_k\bigr)^{-1}\bgma_{k+\tau
}\bigr\rrvert
\ge K\bigr)=o\bigl(n^{-t}\bigr).
\end{eqnarray*}
\end{lemma}

\begin{lemma}\label{kkbar2}For any $v_n\ge n^{-1/212}$, $u\in[a,b]$
and $t>0$, we have
\begin{eqnarray*}
&&\rP\bigl(\bigl\llvert \bgma_{k+\tau}^*\bbA_k^{-2}
\bigl(\bbA^*_k\bigr)^{-2}\bgma_{k+\tau
}\bigr\rrvert
\ge K\bigr)=o\bigl(n^{-t}\bigr)
\end{eqnarray*}
for some $K>0$.
\end{lemma}

\begin{lemma}\label{ep40} Let $u\in[a,b]$, then for any $v_n\ge
n^{-1/212}$, we have
\begin{eqnarray*}
\bigl\llvert \rE\rtr\bbA^{-1}-\rE\rtr\bbA_k^{-1}
\bigr\rrvert &=&O(1)\quad\mbox{and}
\\
\bigl\llvert \rE\rtr\bbA_{k,\ldots,k+(s-1)\tau}^{-1}-\rE\rtr
\bbA_{k,\ldots,k+s\tau}^{-1}\bigr\rrvert &=&O(1).
\end{eqnarray*}
\end{lemma}
%

\section{A convergence rate of the empirical spectral distribution}\label{s4}\label{sec3}

In this section, we give a convergence rate of $\llVert  F_n-F_{c_n}\rrVert  $.
\subsection{A preliminary convergence rate of $m_{n}(z)-\mathrm{E}m_{n}(z)$}\label{sec4.1}\label{sec31}
Let $\rE_k$ denote the conditional expectation given $\bgma
_{k+1},\ldots,\bgma_{T+\tau}$. With this notation, we have
$m_{n}(z)=\rE_{0}(m_{n}(z))$ and $\rE m_{n}(z)=\rE_{T}(m_{n}(z))$.
Therefore, we obtain
\begin{eqnarray*}
m_{n}(z)-\rE m_{n}(z)&=&\sum
_{k=1}^{T+\tau}\bigl(\rE_{k-1}m_{n}(z)-
\rE _km_{n}(z)\bigr)
\\
&=&\sum_{k=1}^{T+\tau}\frac{1}{n}(
\rE_{k-1}-\rE_k) \bigl(\rtr A^{-1}-\rtr
A_{k}^{-1}\bigr)
\\
&\equiv&\sum_{k=1}^{T+\tau}\frac{1}{n}(
\rE_{k-1}-\rE_k)\alpha_{k}.
\end{eqnarray*}
Write
\begin{eqnarray*}
\bbM_n&=&\bbM_{n,k}+(\bgma_{k+\tau},
\bgma_k, \bgma_{k-\tau
})\pmatrix{ 0& 1 & 0
\cr
1 & 0 & 1
\cr
0
& 1 & 0}\pmatrix{ \bgma_{k+\tau}^*
\cr
\bgma_k^*
\cr
\bgma_{k-\tau}^* }
\\
&\equiv&\bbM_{n,k}+\mathbf{C}_k.
\end{eqnarray*}
Let $\la_i(\bbB)$ denote the $i$th smallest eigenvalue for a
Hermitian matrix $\bbB$. Then, for any $i>3$, we have
%
\begin{eqnarray}\label{eqlacing}
\la_i(\bbM_n) &=&\sup_{\bolds \alpha_1,\ldots,\bolds \alpha_{i-1}}\mathop{\inf
_{\bolds \beta \bot\bolds \alpha_1,\ldots,\bolds \alpha_{i-1}}}_{\llVert
\bolds \beta \rrVert  =1}\bigl(\bolds \beta ^*\bbM _{n,k}
\bolds \beta +\bolds \beta ^*\mathbf{C}_k \bolds\beta\bigr)
\nonumber
\\
&\ge&\sup_{\bolds \alpha_1,\ldots,\bolds \alpha_{i-4}}\mathop{\inf_{\bolds \beta \bot\bolds \alpha_1,\ldots,\bolds \alpha_{i-4},\bgma_{k+\tau
}, \bgma_k, \bgma_{k-\tau}}}_{\llVert  \bolds \beta \rrVert  =1}\bolds \beta^*\bbM_{n,k}\bolds\beta
\nonumber\\[-8pt]\\[-8pt]\nonumber
&\ge&\sup_{\bolds \alpha_1,\ldots,\bolds \alpha_{i-4}}\mathop{\inf_{\bolds \beta \bot\bolds \alpha_1,\ldots,\bolds \alpha_{i-4}}}_{\llVert
\bolds \beta \rrVert  =1}\bolds \beta^*\bbM _{n,k}\bolds \beta
\nonumber
\\
&=&\la_{i-3}(\bbM_{n,k}). \nonumber
\end{eqnarray}
Similarly, we have $\la_i(\bbM_n)\le\la_{i+3}(\bbM_{n,k})$.
Therefore, with
\[
G(x):=\sum_{i=1}^nI_{\{\la_i(\bbM_n)\le x\}}\quad\mbox{and}\quad G_k(x):=\sum_{i=1}^nI_{\{\la_i(\bbM_{n,k})\le x\}},
\]
we have
%
\begin{eqnarray} \label{eqlacing1}
\llvert \alpha_{k}\rrvert &=&\bigl\llvert \rtr\bbA^{-1}-
\rtr\bbA_{k}^{-1}\bigr\rrvert
\nonumber
\\
&=&\biggl\llvert \int\frac{1}{x-z}\,d\bigl(G(x)-G_k(x)\bigr)
\biggr\rrvert
\nonumber
\\
&\le&\int\frac{\llvert  G(x)-G_k(x)\rrvert  }{\llvert  x-z\rrvert  ^2}\,dx
\\
&\le&3\int\frac{1}{(x-u)^2+v_n^2}\,dx
\nonumber
\\
&\le&\frac{3\pi}{v_n}.\nonumber
\end{eqnarray}
Here, the third equality follows from integration by parts.
Therefore, by Lem\-ma~\ref{Burkholder1},
%
\begin{eqnarray}\label{W1}
\rP\bigl(\bigl\llvert m_n(z)-\rE m_n(z)\bigr\rrvert
>v_n\bigr)&=&\rP \Biggl(\Biggl\llvert \sum
_{k=1}^{T+\tau
}(\rE_{k-1}-\rE_k)
\alpha_k\Biggr\rrvert >nv_n \Biggr)
\nonumber
\\
&\leq&\rE \Biggl(\frac{1}{(nv_n)^p}\Biggl\llvert \sum
_{k=1}^{T+\tau}(\rE _{k-1}-\rE_k)
\alpha_k\Biggr\rrvert ^p \Biggr)
\nonumber\\[-8pt]\\[-8pt]\nonumber
&\leq&\frac{K}{(nv_n)^p}\rE \Biggl(\sum_{k=1}^{T+\tau}
\bigl\llvert (\rE _{k-1}-\rE_k)\alpha_k\bigr
\rrvert ^2 \Biggr)^{p/2}
\nonumber
\\
&\leq&Kn^{-\sfrac{p}{2}}v_n^{-2p}.\nonumber
\end{eqnarray}
Hence, when $v_n\ge n^{-\alpha}$ for some $0<\alpha<\frac{1}{4}$, we
can choose $p>1$ such that $p(\frac{1}{2}-2\alpha)>t$, and thus
%
\begin{eqnarray}
\rP\bigl(\bigl\llvert m_n(z)-\rE m_n(z)\bigr\rrvert
>v_n\bigr)=o\bigl(n^{-t}\bigr),\label{eqw2}
\end{eqnarray}
for any fixed $t>0$. This implies $\llvert  m_{n}(z)-\rE m_{n}(z)\rrvert  =o(v_n)$, a.s.
\subsection{A preliminary convergence rate of $\mathrm{E}m_n(z)-m_n^0(z)$}\label{sec4.2}\label{prev}

Next, we want to show that when $v_n\ge n^{-1/52}$,
%
\begin{eqnarray}
\bigl\llvert \rE m_n(z)-m_n^0(z)\bigr\rrvert
=o(v_n). \label{eqw3}
\end{eqnarray}
By
\[
\bbA=\sum_{k=1}^T\bigl(
\bgma_k\bgma_{k+\tau}^*+\bgma_{k+\tau}\bgma
_k^*\bigr)-z\bbI_n %
\]
we have
\[
\bbI_n=\sum_{k=1}^T\bigl(
\bbA^{-1}\bgma_k\bgma_{k+\tau}^*+\bbA ^{-1}
\bgma_{k+\tau}\bgma_k^*\bigr)-z\bbA^{-1}. %
\]
Taking trace and dividing by $n$, we obtain
\begin{eqnarray*}
1+zm_n(z)=\frac{1}n\sum_{k=1}^T
\bigl(\bgma_{k+\tau}^*\bbA^{-1}\bgma _k+
\bgma_k^*\bbA^{-1}\bgma_{k+\tau}\bigr).
\end{eqnarray*}
Taking expectation on both sides, we obtain
\begin{eqnarray*}
1+z\rE m_n(z)=\frac{1}n\sum_{k=1}^T
\rE\bgma_k^*\bbA^{-1}(\bgma _{k+\tau}+
\bgma_{k-\tau}),
\end{eqnarray*}
or equivalently, by noticing $1-\frac{c_n^2}{2x_{n1}}\rE
^2m_n(z)=x_{n1}-x_{n0}$,
\begin{eqnarray}\label{Dn}
&& c_n+c_nz\rE m_n(z)\nonumber
\\
&&\qquad =\frac{1}T\sum_{k=1}^T\rE
\bgma_k^*\bbA ^{-1}(\bgma_{k+\tau}+
\bgma_{k-\tau})
\nonumber
\\
&&\qquad = \frac{1}T\sum_{k=1}^T
\biggl[1-\rE\frac{1}{1+\bgma_k^*\tilde
\bbA_k^{-1}(\bgma_{k+\tau}+\bgma_{k-\tau})} \biggr]
\nonumber\\[-8pt]\\[-8pt]\nonumber
&&\qquad = \frac{1}T\sum_{k=1}^T
\biggl[1-\rE
\biggl({1}\Big /\biggl(1+\bgma_k^*\bbA
_k^{-1}(\bgma_{k+\tau}+\bgma_{k-\tau})\nonumber
\\
&&\hspace*{119pt}{} -\frac{\bgma_k^*\bbA_k^{-1}
\bgma_k(\bgma_{k+\tau}^*+\bgma_{k-\tau}^*)\bbA_k^{-1}(\bgma
_{k+\tau}+\bgma_{k-\tau})}{1+(\bgma_{k+\tau}^*+\bgma_{k-\tau
}^*)\bbA_k^{-1}\bgma_{k}}\biggr)\biggr) \biggr]
\nonumber
\\
&&\qquad =1-\frac{1}{1-(\sfrac{c_n^2}{(2x_{n1})})\rE^2m_n(z)}+\delta_n,\nonumber
\end{eqnarray}
where
\begin{eqnarray*}
\tilde\bbA_k&=&\bbA-(\bgma_{k+\tau}+\bgma_{k-\tau})
\bgma _k^*=\bbA_k+\bgma_k\bigl(
\bgma_{k+\tau}^*+\bgma_{k-\tau}^*\bigr),
\\
\delta_n&=&-\frac{1}T\sum_{k=1}^T
\biggl(\rE \biggl( 1\Big/ \biggl(1+\bgma_k^*\bbA
_k^{-1}(\bgma_{k+\tau}+\bgma_{k-\tau})
\\
&&\hspace*{77pt}{} -\frac{\bgma_k^*\bbA_k^{-1}
\bgma_k(\bgma_{k+\tau}^*+\bgma_{k-\tau}^*)\bbA_k^{-1}(\bgma
_{k+\tau}+\bgma_{k-\tau})}{1+(\bgma_{k+\tau}^*+\bgma_{k-\tau
}^*)\bbA_k^{-1}\bgma_{k}}\biggr)\biggr)
\\
&&\hspace*{245pt}{} -\frac{1}{x_{n1}-x_{n0}} \biggr),
\end{eqnarray*}
$x_{n1}, x_{n0}$ are the roots of the equation $x^2=x-a_n^2$ with
$\llvert  x_{n1}\rrvert  >\llvert  x_{n0}\rrvert  $, and $a_n=\frac{c_n\rE m_n}{2}$, as defined below
the statement of Lemma~\ref{rs}. Substituting the expression of
$x_{n1}$, we have
%
\begin{eqnarray}
\label{Emn} \bigl(1-c_n^2\bigl(\rE m_n(z)
\bigr)^2\bigr) \bigl(c_n+c_nz\rE
m_n(z)-1-\delta_n\bigr)^2=1.
\end{eqnarray}
Meanwhile, by (3.8) of \citet{Jinetal14}, we have
%
\begin{eqnarray}\label{eq4.8}
\label{m} \bigl(1-c^2m^2(z)\bigr) \bigl(c+czm(z)-1
\bigr)^2=1.
\end{eqnarray}
Similarly, $m_n^0(z)$ satisfies
%
\begin{eqnarray}
\label{mn} \bigl(1-c_n^2\bigl(m_n^0(z)
\bigr)^2\bigr) \bigl(c_n+c_nzm_n^0(z)-1
\bigr)^2=1.
\end{eqnarray}

We can regard the three expressions above as polynomials of $\rE
m_n(u+iv_n)$, $m(u)$ and $m_n^0(u+iv_n)$, respectively. Compared with
(\ref{m}), coefficients in (\ref{Emn}) and~(\ref{mn}) are different
in terms of $\delta_n$ and $c_n$.

\subsubsection{Identification of the solution to equation \texorpdfstring{(\protect\ref{m})}{(4.8)}}\label{sec4.2.1}\label{uniqueness}
In this subsection, we show that for $c\ne1$ and every $A>0$, there is
a constant $\eta>0$ such that
for every $z$ with $\Im(z)\in(0,\eta)$ and $\llvert  \Re(z)\rrvert  \le A$,
equation (\ref{m})
\begin{eqnarray*}
\bigl(1-c^2m^2(z)\bigr) \bigl(1-c-czm(z)
\bigr)^2=1
\end{eqnarray*}
has only one solution satisfying $\Im(m(z))>\eta v$ and the other
three satisfying $\Im(m(z))<-\eta v$ when $c<1$; and one satisfying
$\Im(m(z)+\frac{c-1}{cz})>\eta v$ and the other three satisfying $\Im
(m(z)+\frac{c-1}{cz})<-\eta v$ when $c>1$.

At first, we claim that the statement is true when $\llvert  z\rrvert  <\delta$ for
some small positive~$\delta$. In \citet{Jinetal14}, it has been proved
that the four solutions for a $z$ with $\Im(z)>0$ are
\begin{eqnarray*}
m_1(z)&=&\frac{(\vfrac{1-c}{z}+\sqrt{1+y_0})+\sqrt{(\vfrac
{1-c}{z}-\sfrac{1}{\sqrt{1+y_0}})^2-\afrac{y_0^2}{1+y_0}}}{2c},
\\
m_2(z)&=&\frac{(\vfrac{1-c}{z}+\sqrt{1+y_0})-\sqrt{(\vfrac
{1-c}{z}-\sfrac{1}{\sqrt{1+y_0}})^2-\afrac{y_0^2}{1+y_0}}}{2c},
\\
m_3(z)&=&\frac{(\vfrac{1-c}{z}-\sqrt{1+y_0})+\sqrt{(\vfrac
{1-c}{z}+\sfrac{1}{\sqrt{1+y_0}})^2-\afrac{y_0^2}{1+y_0}}}{2c},
\\
m_4(z)&=&\frac{(\vfrac{1-c}{z}-\sqrt{1+y_0})-\sqrt{(\vfrac
{1-c}{z}+\sfrac{1}{\sqrt{1+y_0}})^2-\afrac{y_0^2}{1+y_0}}}{2c},
\end{eqnarray*}
where as convention, we assume that the square root of a complex number
is the one with positive imaginary part, and $y_0$ is the root of the
largest absolute value to the equation
\begin{eqnarray*}
y^3-\frac{(1-c)^2-z^2}{z^2}y^2-\frac{4}{z^2}y-
\frac{4}{z^2}=0
\end{eqnarray*}
or equivalently
%
\begin{eqnarray}
z^2y^3-\bigl((1-c)^2-z^2
\bigr)y^2-4y-4=0. \label{eqb1}
\end{eqnarray}

We first consider the case where $z\to0$.
At first, by Lemma 4.1 of \citet{BaiMiaRao91}, we see that $y_0\to\infty$ as $z\to0$. Dividing both sides
of (\ref{eqb1}) by~$y^2$, we obtain that
$y_0=\frac{(1-c)^2}{z^2}(1+o(1))$. Writing $y_0=\frac
{(1-c)^2}{z^2}+d$ and substituting it into (\ref{eqb1}), we obtain
%
\begin{eqnarray}\label{eqb2}
&&\frac{(1-c)^6}{z^4}+3d\frac
{(1-c)^4}{z^2}+3d^2(1-c)^2+d^3z^2
\nonumber
\\
&&\quad{} -\bigl((1-c)^2-z^2\bigr) \biggl(\frac{(1-c)^4}{z^4}+
\frac{2d(1-c)^2}{z^2}+d^2 \biggr)-\frac{4(1-c)^2}{z^2}-4d-4\hspace*{-10pt}
\nonumber\\[-8pt]\\[-8pt]\nonumber
&&\qquad =\frac{d(1-c)^4}{z^2}-\frac
{4(1-c)^2-(1-c)^4}{z^2}+2\bigl(d^2+d\bigr)
(1-c)^2-4(d+1)
\nonumber
\\
&&\quad\qquad{} +\bigl(d^3+d^2\bigr)z^2=0.\nonumber
\end{eqnarray}

By equation (\ref{eqb2}), we have
\begin{eqnarray*}
d=\frac{4}{(1-c)^2}-1+O\bigl(z^2\bigr).
\end{eqnarray*}
That is,
%
\begin{eqnarray}
y_0=\frac{(1-c)^2}{z^2}+\frac{4}{(1-c)^2}-1+O\bigl(z^2
\bigr). \label{eqb3}
\end{eqnarray}
Therefore, we have
%
\begin{eqnarray}
\sqrt{1+y_0}=-\frac{\llvert  1-c\rrvert  }{z} \biggl(1+\frac{2z^2}{(1-c)^4}+O
\bigl(z^4\bigr) \biggr). \label{eqb4}
\end{eqnarray}
Consequently,
%
\begin{eqnarray}
\frac{1-c}{z}+\sqrt{1+y_0}=\frac{1-c-\llvert  1-c\rrvert  }{z}-
\frac
{2z}{\llvert  1-c\rrvert  ^3}+O\bigl(z^3\bigr),\label{eqb5}
\\
\frac{1-c}{z}-\sqrt{1+y_0}=\frac{1-c+\llvert  1-c\rrvert  }{z}+
\frac
{2z}{\llvert  1-c\rrvert  ^3}+O\bigl(z^3\bigr).\label{eqb6}
\end{eqnarray}
%

Because
\begin{eqnarray*}
&& \biggl(\frac{1-c}{z}\mp\frac{1}{\sqrt{1+y_0}} \biggr)^2-
\frac
{y_0^2}{1+y_0}=\frac{(1-c)^2}{z^2}\mp2\frac{1-c}{z\sqrt
{1+y_0}}+1-y_0
\\
&&\qquad =-\frac{4}{(1-c)^2}\pm2\frac{1-c}{\llvert  1-c\rrvert  +O(z^2)}+2+O\bigl(z^2\bigr)
\\
&&\qquad =-\frac{4}{(1-c)^2}\pm2\frac{1-c}{\llvert  1-c\rrvert  }+2+O\bigl(z^2\bigr),
\end{eqnarray*}
we obtain
%
\begin{eqnarray}\label{eqb7}
&&\sqrt{ \biggl(\frac{1-c}{z}\mp\frac{1}{\sqrt{1+y_0}}
\biggr)^2-\frac
{y_0^2}{1+y_0}}
\nonumber\\[-8pt]\\[-8pt]\nonumber
&&\qquad =i\sqrt{\frac{4}{(1-c)^2}\mp2
\frac{1-c}{\llvert  1-c\rrvert  }-2}+O\bigl(z^2\bigr) .
\end{eqnarray}
When $c<1$, from (\ref{eqb5}) and (\ref{eqb7}), as $z\to0$, we obtain
%
\begin{eqnarray}\label{eqbb4}
\Im(2cm_1)&=&\Im \biggl(O(z)+i\sqrt{\frac{4}{(1-c)^2}-4} \biggr)>
\frac{\sqrt{c(2-c)}}{(1-c)},
\nonumber
\\
\Im(2cm_2)&=&\Im \biggl(O(z)-i\sqrt{\frac{4}{(1-c)^2}-4} \biggr)<-
\frac{\sqrt{c(2-c)}}{(1-c)},
\nonumber\\[-8pt]\\[-8pt]\nonumber
\Im(2cm_3)&=&\Im \biggl(\frac{2(1-c)}{z}+i\sqrt{
\frac
{4}{(1-c)^2}}+O(z) \biggr)<-\frac{1-c}{\llvert  z\rrvert  ^2}v,
\nonumber
\\
\Im(2cm_4)&=&\Im \biggl(\frac{2(1-c)}{z}-i\sqrt{
\frac
{4}{(1-c)^2}}+O(z) \biggr)<-\frac{2}{(1-c)}. \nonumber
\end{eqnarray}

When $c\in(1,2]$, as $z\to0$, we have
%
\begin{eqnarray}\label{eqbb5418}
\Im\biggl(2c\biggl(m_1+\frac{c-1}{cz}\biggr)\biggr)&=&\Im
\biggl(i\sqrt{\frac
{4}{(1-c)^2}}+O(z) \biggr)>\frac{1}{c-1},
\nonumber
\\
\Im\biggl(2c\biggl(m_2+\frac{c-1}{cz}\biggr)\biggr)&=&\Im
\biggl(-i\sqrt{\frac
{4}{(1-c)^2}}+O(z) \biggr)<-\frac{1}{c-1},\nonumber
\\
\Im\biggl(2c\biggl(m_3+\frac{c-1}{cz}\biggr)\biggr)&=&\Im
\biggl(\frac{2(c-1)}{z}+i\sqrt {\frac{4c(2-c)}{(1-c)^2}}+O(z) \biggr)
\nonumber\\[-8pt]\\[-8pt]\nonumber
&<& - \frac{c-1}{\llvert  z\rrvert  ^2}v,
\nonumber
\\
\Im\biggl(2c\biggl(m_4+\frac{c-1}{cz}\biggr)\biggr)&=&\Im
\biggl(\frac{2(c-1)}{z}-i\sqrt {\frac{4c(2-c)}{(1-c)^2}}+O(z) \biggr)\nonumber
\\
&<& -\frac{c-1}{\llvert  z\rrvert  ^2}v.\nonumber
\end{eqnarray}

When $c>2$, as $z\to0$, we have
%
\begin{eqnarray}\label{eqbb519}
\Im\biggl(2c\biggl(m_1+\frac{c-1}{cz}\biggr)\biggr)&=&\Im
\biggl(i\sqrt{\frac
{4}{(1-c)^2}}+O(z) \biggr)>\frac{1}{c-1},
\nonumber
\\
\Im\biggl(2c\biggl(m_2+\frac{c-1}{cz}\biggr)\biggr)&=&\Im
\biggl(-i\sqrt{\frac{4}{(1-c)^2}}+O(z) \biggr)<-\frac{1}{c-1},
\nonumber
\\
\Im\biggl(2c\biggl(m_3+\frac{c-1}{cz}\biggr)\biggr)&=&\Im
\biggl(\frac{2(c-1)}{z}-\sqrt {\frac{4c(c-2)}{(1-c)^2}}+O(z) \biggr)\nonumber
\\
&<& -\frac{c-1}{\llvert  z\rrvert  ^2}v,
\\
\Im\biggl(2c\biggl(m_4+\frac{c-1}{cz}\biggr)\biggr)&=&\Im
\biggl(\frac{2(c-1)}{z}+\sqrt {\frac{4c(c-2)}{(1-c)^2}}+O(z) \biggr)\nonumber
\\
&<& - \frac{c-1}{\llvert  z\rrvert  ^2}v. \nonumber
\end{eqnarray}
This proves the result when $\llvert  z\rrvert  <\gd$ for some $\gd>0$.

For $\llvert  z\rrvert  \ge\gd$, we first consider the case where $c<1$. Suppose that
$m(z)$ is one of the four continuous branches of the solutions of the
equation (\ref{eq4.8}).
If the conclusion is incorrect for $m(z)$, then there exist a sequence
of constants $\zeta_n\downarrow0$ and a sequence of complex numbers
$z_n=u_n+iv_n$ satisfying $\llvert  z_n\rrvert  \ge\gd$, $\llvert  u_n\rrvert  \le A$, $v_n\in
(0,\eta)$ with $\eta=\gd^2/2$ and $\llvert  \Im(m(z_n))\rrvert  \le\zeta_nv_n$.
Then there is a subsequence $\{n'\}$ such that $z_{n'}\to z_0=u_0+iv_0$
with $u_{n'}\to u_0\in[-A,A]$ and $v_{n'}\to v_0\in[0,\eta]$.

Write $m(z_n)=m_1(z_n)+im_2(z_n)$, where $m_1(z_n)$ and $m_2(z_n)$ are
real. Since $m(z_n)$ satisfies the equation (\ref{eq4.8}), we have
%
\begin{eqnarray}
\bigl(1-c^2m^2(z_n)\bigr)
\bigl(1-c-cz_nm(z_n)\bigr)^2=1. \label{eqzzdd0}
\end{eqnarray}
Comparing the imaginary parts of both sides of (\ref{eqzzdd0}), we obtain
\begin{eqnarray*}
&& c^2m_1(z_n)m_2(z_n)
\\
&&\quad{}\times
\bigl[\bigl(1-c-cu_nm_1(z_n)+cv_nm_2(z_n)
\bigr)^2
-\bigl(cu_nm_2(z_n)+cv_nm_1(z_n)
\bigr)^2\bigr]
\\
&&\qquad{} +\bigl(1-c^2m_1^2(z_n)+c^2m_2^2(z_n)
\bigr) \bigl(cu_nm_2(z_n)+cv_nm_1(z_n)
\bigr)
\\
&&\quad\qquad {}\times  \bigl(1-c-cu_nm_1(z_n)+cv_nm_2(z_n)
\bigr)=0.
\end{eqnarray*}
Dividing by $v_n$ both sides of the equation above, we obtain
%
\begin{eqnarray}
\bigl(1-c^2m_1^2(z_0)\bigr)
\bigl(cm_1(z_0)\bigr) \bigl(1-c-cu_0m_1(z_0)
\bigr)=0. \label{eqzzdd1}
\end{eqnarray}
By the condition that
$\llvert  \Im(m(z_n))\rrvert  \le\zeta_n v_n\to0$, we have that $m(z_0)=m_1(z_0)$
is real. The solutions $\pm1/c$ and $0$ of the equation (\ref
{eqzzdd1}) for $m(z_0)$ do not satisfy equation~(\ref{eq4.8}). Therefore, we have
$1-c-cu_0m(z_0)=0$, and hence by (\ref{eq4.8})
%
\begin{eqnarray}
-\bigl(1-c^2m^2(z_0)\bigr)c^2v_0^2m^2(z_0)=1.
\label{eqzzdd2}
\end{eqnarray}
Note that $v_0=0$ contradicts to the equation above. Thus, we have
$v_0\in(0,\gd^2/2]$. By (\ref{eqzzdd2}) and the fact that
$1-c-cu_0m(z_0)=0$, we obtain
\[
\frac{(1-c)^2}{u_0^2}=\frac{v_0^2+\sqrt{v_0^4+4v_0^2}}{2v_0^2}\quad \mbox{or}\quad u_0^2=
\frac{2v_0^2(1-c)^2}{v_0^2+\sqrt{v_0^4+4v_0^2}}. %
\]
The expression of $u_0^2$ implies that $u_0^2<v_0<\gd^2/2$. On the
other hand, by the assumption that
$\llvert  z_0\rrvert  >\gd$, we have $u_0^2+v_0^2>\gd^2$ and $v_0^2<v_0<\gd^2/2$
which implies that $u_0^2>\gd^2/2$, the contradiction proves our assertion.

Now, we consider the case $c>1$. Let $\underline m(z)=cm(z)+\frac
{c-1}{z}$. Then equation (\ref{eq4.8}) becomes
%
\begin{equation}
z^2\underline m^2(z) \biggl(1- \biggl(
\frac{1-c}{z}+\underline m(z) \biggr)^2 \biggr)=1. \label{eqbm1}
\end{equation}
If the conclusion is untrue, similar to the case where $c<1$, there
exist sequences $\zeta_n\downarrow0$ and
$z_n=u_n+iv_n\to z_0=u_0+i0$ such that $\llvert  \Im(\underline m(z_n))\rrvert  \le
\zeta_nv_n$, and
$\llvert  u_n\rrvert  \le A$. By the continuity of the solution $\underline m(z)$ for
$\llvert  z\rrvert  \ge\gd$, we may assume the inequality above is an equality, for
otherwise, one may shift $\Re(z_n)=u_n$ toward the origin.
Write $\underline m(z_n)=\underline m_1(z_n)+i\underline m_2(z_n)$,
where $\underline m_1(z_n)$ and $\underline m_2(z_n)$ are both real. By
the equality of imaginary parts of (\ref{eqbm1}), we have
%
\begin{eqnarray}\label{eqbm3}
&& \underline m_1(z_n)\underline m_2(z_n)\nonumber
\\
&&\quad{}\times
\bigl(u_n^2-v_n^2-
\bigl(1-c+u_n\underline m_1(z_n)-v_0
\underline m_2(z_n)\bigr)^2\nonumber
\\
&&\hspace*{98pt}{}
+\bigl(u_n\underline m_2(z_n)+v_n
\underline m_1(z_n)\bigr)^2 \bigr)
\nonumber\\[-8pt]\\[-8pt]\nonumber
&&\qquad{}
 -\bigl(
\underline m_1^2(z_n)-\underline
m_2^2(z_n)\bigr)
\\
&&\quad\qquad{}\times  \bigl(u_nv_n
-\bigl(1-c+u_n\underline m_1(z_n)-v_n
\underline m_2(z_n)\bigr) \bigl(u_n\underline
m_2(z_n)+v_n\underline m_1(z_n)
\bigr) \bigr)\nonumber
\\
&&\qquad =0 \nonumber
\end{eqnarray}
Dividing both sides by $v_n$ and making $n\to\infty$ on both sides of
the equation above, by assumption, we obtain
%
\begin{eqnarray}
\underline m_1^2(z_0) \bigl(u_0-
\bigl(1-c+u_0\underline m_1(z_0)\bigr)
\underline m_1(z_0) \bigr)=0. \label{eqbm4}
\end{eqnarray}
This implies that
%
\begin{eqnarray}
u_0=\frac{(1-c)\underline m_1(z_0)}{(1-\underline m^2_1(z_0))}. \label{eqnm2}
\end{eqnarray}
Similarly, we have $\underline m(u_0)=\underline m_1(u_0)$ which is
real. By the real part of (\ref{eqbm1}), we have
\begin{eqnarray*}
\underline m^2(u_0) \bigl(u_0^2-
\bigl(1-c+u_0\underline m(u_0)\bigr)^2
\bigr)=1.
\end{eqnarray*}
The solution to the equation above in $u_0$ is
%
\begin{eqnarray}
u_0=\frac{\underline m^3(u_0)(1-c)\pm\sqrt{\underline
m^2(u_0)-c(2-c)\underline m^4(u_0)}}{\underline m^2(u_0)(1-\underline
m^2(u_0))}. \label{eqbmm1}
\end{eqnarray}
If $\underline m^2(u_0)\ne\frac{1}{c(2-c)}$, then (\ref{eqbmm1})
contradicts (\ref{eqnm2}).

Now, we consider the case where $c\in(1,2)$ and $\underline
m^2(u_0)=\frac{1}{c(2-c)}$.
By differentiating (\ref{eqbm1}) with respect to $z$, we obtain
\begin{eqnarray*}
\frac{d\underline m(z)}{dz}&=&-\frac{\underline m(z-\underline
m(1-c+z\underline m))}{z^2-(1-c+z\underline m)^2-z\underline
m(1-c+z\underline m)}
\\
&=&-\frac{\underline m(z-\underline m(1-c+z\underline
m))}{z^2-(1-c)^2-z(1-c)\underline m}.
\end{eqnarray*}
Because
\begin{eqnarray*}
&& \Im\bigl(z_n-\underline m\bigl(1-c+z_n\underline
m(z_n)\bigr)\bigr) = v_n \bigl[\bigl(1-\underline
m^2_1(u_0)\bigr) +o(1) \bigr],
\\
&& \Re\bigl(z_n-\underline m\bigl(1-c+z_n\underline
m(z_n)\bigr)\bigr)
\\
&&\qquad = \bigl[u_n-\underline
m_1(z_n) \bigl(1-c+u_n\underline
m_1(z_n) \bigr) \bigr]+O\bigl(\underline
m_2(z_n)\bigr)
\\
&&\qquad = \bigl[u_n\bigl(1-\underline m_1^2(z_n)
\bigr)-(1-c)\underline m_1(z_n) \bigr]+O\bigl(\underline
m_2(z_n)\bigr)\qquad\bigl(\mbox{by (\ref{eqbm3})}
\bigr)
\\
&&\qquad =-\frac{\underline m_2(z_n)}{v_n\underline m_1(z_n)} \bigl[u_n^2-(1-c)^2-u_n(1-c)
\underline m_1(z_n)+o(1) \bigr]
\\
&&\qquad \simeq\zeta_n \frac{(1-c)^2[1-2\underline m^2(u_0)]}{\underline
m(u_0)(1-\underline m(u_0)^2)^2},
\\
&& \frac{z^2_n-(1-c)^2-z(1-c)\underline m(z_n)}{\underline m(z_n)} \simeq \frac{(1-c)^2[2\underline m^2(u_0)-1]}{\underline m(u_0)(1-\underline
m^2(u_0))^2}.
\end{eqnarray*}
Therefore,
\[
\frac{\partial\underline m_2(z_n)}{\partial u}\simeq v_n\frac
{\underline m(u_0)(1-\underline m^2(u_0))^3}{(1-c)^2(2\underline m^2(u_0)-1)}, %
\]
and
\[
\frac{\partial\underline m_1(z_n)}{\partial u}\simeq\zeta_n. %
\]
Hence,
%
\begin{eqnarray}
G_n&=&\underline m_1^2(z_n)
\bigl(u_n-\bigl(1-c+u_n\underline m_1(z_n)
\bigr)\underline m_1(z_n) \bigr)
\nonumber
\\
&&{}-\underline m_1^2(z_0)
\bigl(u_0-\bigl(1-c+u_0\underline m_1(z_0)
\bigr)\underline m_1(z_0) \bigr)
\\
&=&(u_n-u_0) \bigl(\underline m_1^2
\bigl(z_n^*\bigr) \bigl(1-\underline m_1^2(z_0)
\bigr)+O(\zeta_n) \bigr).\nonumber
\end{eqnarray}
On the other hand, we have
%
\begin{eqnarray}
\zeta_nv_n&=&\underline m_2(z_n)-
\underline m_2(z_0)
\nonumber
\\
&=&(u_n-u_0)\frac{\partial\underline m_2(z_n^*)}{\partial u}
\\
&\simeq&(u_n-u_0)\frac{v_n\underline m(z_0)(1-\underline
m^2(z_0))^3}{(1-c)^2(2\underline m^2(z_0)-1)}.\nonumber
\end{eqnarray}
%
Therefore,
%
\begin{eqnarray}
G_n&\simeq&\zeta_n\frac{(1-c)^2\underline m(u_0)(2\underline
m^2(z_0)-1)}{(1-\underline m^2(z_0))^2}.
\end{eqnarray}
Substituting the above into (\ref{eqbm3}) and dividing $\underline
m_2(z_n)=\zeta_nv_n$ on both sides and letting $n\to\infty$, we obtain
%
\begin{eqnarray}
\qquad 0&=&\underline m(u_0) \bigl(u_0^2-
\bigl(1-c+u_0\underline m(u_0)\bigr)^2\bigr)+
\underline m^2(u_0) \bigl(1-c+u_0\underline
m(u_0)\bigr)u_0
\nonumber
\\
&&{}+\frac{(1-c)^2\underline m(u_0)(2\underline
m^2(u_0)-1)}{(1-\underline m^2(u_0))^2}
\nonumber\\[-8pt]\\[-8pt]\nonumber
&=&\underline m(u_0) \bigl(u_0^2-(1-c)^2-u_0(1-c)
\underline m(u_0)\bigr)
\\
&&{} +\frac
{(1-c)^2\underline m(u_0)(2\underline m^2(u_0)-1)}{(1-\underline m^2(u_0))^2}.\nonumber
\end{eqnarray}
By substitution of (\ref{eqnm2}), the equation above becomes
\begin{eqnarray*}
\frac{2(1-c)^2\underline m(u_0)(2\underline m^2(u_0)-1)}{(1-\underline
m^2(u_0))^2}=0
\end{eqnarray*}
which also implies that $\underline m^2(u_0)=\frac{1}2$.
This contradicts to the assumption that $\underline m^2(u_0)=\frac
{1}{c(2-c)}$ and the assertion is finally proved.

Consequently, under the condition that $\llvert  \gd_n\rrvert  \le Kv_n^{\eta}$ with
$\eta>1$, we have
$\max_{j=2,3,4, z=u+iv_n}\llvert  m_j(z)-\rE m_n(z)\rrvert  \ge\eta v_n$ and thus\vspace*{1pt}
$\max_{z=u+iv_n}\llvert  m_1(z)- \rE m_n(z)\rrvert  \le Kv_n^{\eta}$ when $c<1$.
Similarly for $\underline m(z)$ when $c>1$.

Hence, to prove (\ref{eqw3}), it remains to show
%
\begin{eqnarray}
\label{Kc} \llvert \delta_n\rrvert \le K v_n^\eta
\end{eqnarray}
for some $K>0$, and $\eta>1$.
\subsubsection{Convergence rate of \texorpdfstring{$\delta_n$}{deltan}}\label{sec4.2.2}
Let $v_n\ge n^{-1/52}$. By (\ref{Dn}), we have
\begin{eqnarray*}
\delta_n &=&c_n+c_nz\rE m_n(z)-1+
\frac{1}{x_{n1}-x_{n0}}=:\frac{1}T\sum_{k=1}^T
\rE\eta_k,
\end{eqnarray*}
where
\begin{eqnarray*}
\eta_k&=&\bgma_k^*\bbA^{-1}(
\bgma_{k+\tau}+\bgma_{k-\tau
})-1+\frac{1}{x_{n1}-x_{n0}}.
\end{eqnarray*}
When $k\le v_n^{-4}$ or $\ge T-v_n^{-4}$, by (iii)(a) of Lemma~\ref
{x1}, we have
\begin{eqnarray*}
\llvert \rE\eta_k\rrvert &\le& v_n^{-1}\sqrt{
\rE\llvert \bgma_k\rrvert ^2\bigl(\rE\llvert
\bgma_{k-\tau
}\rrvert ^2+\rE\llvert \bgma_{k+\tau}
\rrvert ^2\bigr)}+1+\frac{1}{\llvert  x_{n1}-x_{n0}\rrvert  }
\\
&\le&Kv_n^{-1}.
\end{eqnarray*}
Therefore, for all large $n$,
\begin{eqnarray}
\frac{1}T \Biggl(\sum_{k=1}^{[v_ n^{-4}]}+
\sum_{k=[T-v_n^{-4}
]}^T \Biggr)\llvert \rE
\eta_k\rrvert \le\frac{K}{Tv_n^5}\le Kv_n^{47}.
\label{eqw4}
\end{eqnarray}
When $k\in([v^{-4}_n],[T-v_n^{-4}])$,
denote
%
\begin{eqnarray} \label{defep}
\ep_1&=&\bigl(\bgma_{k+\tau}^*+\bgma_{k-\tau}^*\bigr)
\bbA_k^{-1}\bgma _{k},
\nonumber
\\
\ep_2&=&\bgma_k^*\bbA_k^{-1}(
\bgma_{k+\tau}+\bgma_{k-\tau
}),
\nonumber
\\
\ep_3&=&\bgma_k^*\bbA_k^{-1}
\bgma_k-\frac{1}{2T}\rtr\bbA _k^{-1},
\\
\ep_4&=&\frac{1}{2T}\rtr\bbA_k^{-1}-
\frac{c_n}{2}\rE m_n(z),
\nonumber
\\
\ep_5&=&\bigl(\bgma_{k+\tau}^*+\bgma_{k-\tau}^*\bigr)
\bbA_k^{-1}(\bgma _{k+\tau}+\bgma_{k-\tau})-
\frac{c_n}{x_{n1}}\rE m_n(z).\nonumber
\end{eqnarray}
Then, by the fact that $x_{n1}-x_{n0}=1-2a_n^2/x_{n1}$, we have
\begin{eqnarray*}
-\rE\eta_k
&=&\rE \biggl(1\Big/ \biggl(1+\bgma_k^*\bbA_k^{-1}(\bgma_{k+\tau}+\bgma
_{k-\tau})
\\
&&\hspace*{33pt}{}-\frac{\bgma_k^*\bbA_k^{-1}
\bgma_k(\bgma_{k+\tau}^*+\bgma_{k-\tau}^*)\bbA_k^{-1}(\bgma
_{k+\tau}+\bgma_{k-\tau})}{1+(\bgma_{k+\tau}^*+\bgma_{k-\tau}^*)
\bbA_k^{-1}\bgma_{k}}\biggr)\biggr)-\frac{1}{x_{n1}-x_{n0}}
\\
&=&\frac{1}{x_{n1}-x_{n0}} \rE\beta_{k} \biggl(-2\ep_1
\frac{a_n^2}{x_{n1}}-\ep_2-\ep_1\ep_2
\\
&&\hspace*{68pt}{} +\bigl(\bgma_{k+\tau}^*+\bgma_{k-\tau}^*\bigr)
\bbA_k^{-1} (\bgma_{k+\tau}+\bgma_{k-\tau}) (
\ep_3+\ep_4)+a_n\ep_5 \biggr),
\end{eqnarray*}
where
\begin{eqnarray*}
\beta_{k}&=&\frac{1}{1+\ep_1+\ep_2+\ep_1\ep_2-\bgma_k^*\bbA_k^{-1}
\bgma_k(\bgma_{k+\tau}^*+\bgma_{k-\tau}^*)\bbA_k^{-1}(\bgma
_{k+\tau}+\bgma_{k-\tau})}
\\
&=&\frac{1}{1+\ep_1+\ep_2+\ep_1\ep_2-(a_n+\ep_3+\ep_4)(\sfrac
{2a_n}{x_{n1}}+\ep_5)}.
\end{eqnarray*}

Define a random set ${\mathcal E}_n=\{| \ep_i| \le v_n^6,i=1,2,3,4,5\}$.
When ${\mathcal E}_n$ happens, by the facts $| a_n| \le Kv_n^{-1}$, $
| \frac{2a_n}{x_{n1}}| \le2$ and Lemma~\ref{x1}(iii)(a), we have
\begin{eqnarray*}
\llvert \beta_k\rrvert &\le&\frac{1}{\llvert  1-\sfrac
{2a_n^2}{x_{n1}}-9v_n^6-Kv_n^5\rrvert  }
\\
&=&\frac{1}{
\llvert  1-2x_{n0}-9v_n^6-Kv_n^5\rrvert  }
\\
&=&
\frac{1}{
\llvert  x_{n1}-x_{n0}-9v_n^6-Kv_n^5\rrvert  }
\\
&\le& Kv_n^{-1}.
\end{eqnarray*}
Together with Lemma~\ref{x1}(ii)(a) and (iii)(a), we obtain that
\begin{eqnarray*}
\llvert \eta_k\rrvert &\le&\frac{1}{\llvert  x_{n1}-x_{n0}\rrvert  }
\\[-1pt]
&&{}\times Kv_n^{-1}
\bigl(v_n^6\bigl(2\llvert x_{n0}\rrvert
\bigr)+v_n^6+v_n^{12}+v_n^{-1}
\llVert \bgma_{k+\tau}+\bgma _{k-\tau}\rrVert ^2
\bigl(2v_n^6\bigr)+Kv_n^{5} \bigr)
\\[-1pt]
&\le& Kv_n^3.
\end{eqnarray*}
Therefore, by Lemmas~\ref{rs},~\ref{Ek+tau}(a) and~\ref{cross}(a),
when $v_n\ge n^{-1/52}$, we have
%
\begin{eqnarray}\label{eqx5}
\rE\llvert \eta_{k}\rrvert &\le& Kv_n^3+Kv_n^{-1}
\Biggl(\sum_{i=1}^5\rP\bigl(\llvert \ep
_i\rrvert \ge v_n^6\bigr) \Biggr)
\nonumber
\nonumber\\[-10pt]\\[-10pt]\nonumber
&\le&Kv_n^3.
\end{eqnarray}
Then the conclusion (\ref{Kc}) follows from (\ref{eqw4}) and (\ref{eqx5}).

\subsection{Convergence rate of \texorpdfstring{$\|F_n-F_{c_n}\|$}{||Fn-Fcn||}}\label{sec4.3}
Choose $v_n=n^{-1/52}$. Let $F_n$ be the empirical distribution
function of $\bbM_n$ and $F_{c_n}$ be the LSD with the ratio parameter
$c_n=n/T$ whose Stieltjes transform is denoted by $m_n^0$.
By (\ref{upbd}), let $B=(1+\sqrt c)^2+\delta$, and we have
$F_{c_n}([-B,B]) = 1$. By Lemma~\ref{B14}
we have, for some $A>B$ and $a>0$,
\begin{eqnarray*}
&&\rP \bigl(\llVert F_n-F_{c_n}\rrVert >c'
\sqrt{v_n} \bigr)
\\[-1pt]
&&\qquad  \le\rP \Bigl(\sup_{u\in[-A,A]}\bigl\llvert m_n(z)-m_n^0(z)
\bigr\rrvert >K_0\sqrt {v_n} \Bigr)
\\[-1pt]
&&\quad\qquad{}  +\rP \biggl(\sup_x\int_{\llvert  y\rrvert  \le2v_na} \bigl
\llvert F_{c_n}(x + y)-F_{c_n}(x)\bigr\rrvert
\,dy>K_0\bigl(c'-1\bigr)v_n^{3/2}
\biggr)
\\[-1pt]
&&\qquad \le\rP \biggl(\sup_{u\in[-A,A]}\bigl\llvert m_n(z)-\rE
m_n(z)\bigr\rrvert >\frac{K_0\sqrt
{v_n}}2 \biggr)
\\[-1pt]
&&\quad\qquad{} +\rP \biggl(\sup
_{u\in[-A,A]}\bigl\llvert \rE m_n(z)-m_n^0(z)
\bigr\rrvert >\frac{K_0\sqrt{v_n}}2 \biggr)
\\[-1pt]
&&\quad\qquad{} +\rP \biggl(\sup_x\int_{\llvert  y\rrvert  \le2v_na} \bigl
\llvert F_{c_n}(x + y)-F_{c_n}(x)\bigr\rrvert
\,dy>K_0\bigl(c'-1\bigr)v_n^{3/2}
\biggr),
\end{eqnarray*}
where $K_0=\pi(1-\kappa)(2\gamma-1)$, and $a$ is a constant defined
in Lemma~\ref{B14}.
By $\llvert  \rE m_n(z)-m_n^0(z)\rrvert  =o(v_n)$, the second probability is $0$ for
all large $n$.

By the analysis of Section~3 of \citet{Jinetal14}, we see that $\phi
_{c_n}(x):=\frac{d}{dx}F_{c_n}(x)\le K\llvert  x\rrvert  ^{-1/2}$, which implies that
$F_{c_{n}}$ satisfies the Lipschitz condition with index $\frac{1}2$.
Hence, for some large $c'$, we have
\begin{eqnarray*}
&&\sup_x\int_{\llvert  y\rrvert  \le2v_na} \bigl\llvert
F_{c_n}(x + y)-F_{c_n}(x)\bigr\rrvert \,dy
\\[-1pt]
&&\qquad \le K\int_{\llvert  y\rrvert  \le2v_na}\llvert y\rrvert ^{{1/2}}\,dy=4Ka^2v_n^{3/2}<K_0
\bigl(c'-1\bigr)v_n^{3/2}.
\end{eqnarray*}
Therefore, the third probability is 0.

For the first probability, let ${\mathcal S}_n$ be the set containing
$n^2$ points that are equally spaced between $-n$ and $n$ and note that
$[-A,A]\subseteq[-n,n]$ for all large $n$. When $\llvert  u_1-u_2\rrvert  \le\frac{2}n$, we have
\begin{eqnarray*}
\bigl\llvert m_n(u_1+iv_n)-m_n(u_2+iv_n)
\bigr\rrvert \le\llvert u_1-u_2\rrvert
v_n^{-2}<\frac{K_0\sqrt
{v_n}}2,
\\
\bigl\llvert m_n^0(u_1+iv_n)-m_n^0(u_2+iv_n)
\bigr\rrvert \le\llvert u_1-u_2\rrvert
v_n^{-2}<\frac{K_0\sqrt{v_n}}2.
\end{eqnarray*}
Therefore, by (\ref{W1}), for any $t>0$, we have
\begin{eqnarray*}
&&\rP \biggl(\sup_{u\in[-A,A]}\bigl\llvert m_n(z)-\rE
m_n(z)\bigr\rrvert >\frac{K_0\sqrt
{v_n}}2 \biggr)
\\
&&\qquad =\rP \biggl(\sup_{u\in{\mathcal S}_n}\bigl\llvert m_n(z)-\rE
m_n(z)\bigr\rrvert >\frac
{K_0\sqrt{v_n}}2 \biggr)
\\
&&\qquad \le n^2\rP \biggl(\bigl\llvert m_n(z)-\rE
m_n(z)\bigr\rrvert >\frac{K_0\sqrt{v_n}}2 \biggr)
\\
&&\qquad \le Kn^{2-\sfrac{p}{2}}v_n^{-p}
\\
&&\qquad =o\bigl(n^{-t}\bigr)
\end{eqnarray*}
by selecting $p$ large enough. Thus, we have proved, for any fixed $t>0$
%
\begin{eqnarray}
\label{Fn} \rP\bigl(\llVert F_n-F_{c_n}\rrVert
>c'n^{-1/104}\bigr)=o\bigl(n^{-t}\bigr).
\end{eqnarray}
Next, let $a'=a-\underline{\ep}$ and $b'=b+\underline{\ep}$ for
some $\underline{\ep}>0$ such that $(a',b')\supseteq[a,b]$ is an
open interval outside the support of $F_{c_n}$ for all $n$ large enough.
By $\llvert  d(c_n)-d(c)\rrvert  \to0$, and hence $[a',b']$ is also outside the
support of $F_{c_n}$. We conclude that $F_{c_n}(b')-F_{c_n}(a')=0$ for
all large $n$. Hence, we have
\begin{eqnarray*}
F_n\bigl\{\bigl[a',b'\bigr]\bigr\}
&=&F_n\bigl(b'\bigr)-F_n
\bigl(a'\bigr)-\bigl(F_{c_n}\bigl(b'
\bigr)-F_{c_n}\bigl(a'\bigr)\bigr)
\\
&\le&2\llVert F_n-F_{c_n}\rrVert .
\end{eqnarray*}
Therefore,
%
\begin{eqnarray}
\label{Fnk} &&\rP \Bigl(\max_{k\le n}\rE_k
\bigl(F_n\bigl\{\bigl[a',b'\bigr]\bigr\}
\bigr)\ge4c' n^{-1/104} \Bigr)
\nonumber
\\
&&\qquad \le\rP \Bigl(\max_{k\le n}\rE_k
\bigl(F_n\bigl\{\bigl[a',b'\bigr]\bigr
\}I_{\{
\llVert  F_n-F_{c_n}\rrVert  < c'n^{-1/104}\}}\bigr)\ge2c' n^{-1/104} \Bigr)
\nonumber
\\
&&\quad\qquad{}+\rP \Bigl(\max_{k\le n}\rE_k\bigl(F_n
\bigl\{\bigl[a',b'\bigr]\bigr\}I_{\{\llVert  F_n-F_{c_n}\rrVert  \ge
c'n^{-1/104}\}}\bigr)
\ge2c'n^{-1/104} \Bigr)
\\
&&\qquad \le0+\rP\Bigl(\max_{k\le n}\rE_k I_{\{\llVert  F_n-F_{c_n}\rrVert  \ge c'n^{-1/104}\}
}\ne0\Bigr)
\nonumber
\\
&&\qquad \le n\rP\bigl(\llVert F_n-F_{c_n}\rrVert \ge
c'n^{-1/104}\bigr) =o\bigl(n^{-t}\bigr)\nonumber
\end{eqnarray}
for any $t>0$.

\section{A refined convergence rate of Stieltjes transform when \texorpdfstring{$u\in[a,b]$}{uin[a,b]}}\label{s5}\label{sec4}

In this section, we are to prove that for $v_n=n^{-1/212}$,
%
\begin{eqnarray}
m_{n}-m_{n}^0=o\bigl(1/(nv_{n})
\bigr)\qquad\mbox{a.s.} \label{eqw41}
\end{eqnarray}
by refining the convergence rates obtained in the last section.

\subsection{A refined convergence rate of $m_{n}-\mathrm{E}m_{n}$}\label{sec5.1}\label{sec42}
In this subsection, we want to show that
%
\begin{eqnarray}
\sup_{u\in[a,b]}\bigl\llvert m_{n}(z)-\rE
m_{n}(z)\bigr\rrvert =o\bigl(1/(nv_{n})\bigr),\qquad
\mbox{a.s.} \label{eqw42}
\end{eqnarray}
First, by recalling that $\tilde\bbA_k=\bbA-(\bgma_{k+\tau}+\bgma
_{k-\tau})\bgma_k^*$ and $\bbA_k=\tilde\bbA_k-\bgma_k(\bgma
_{k+\tau}+\bgma_{k-\tau})^*$, we have
\begin{eqnarray*}
\hspace*{-3pt}&&m_{n}(z)-\rE m_{n}(z)
\\
\hspace*{-3pt}&&\qquad =\sum_{k=1}^{T}\bigl(
\rE_{k-1}m_{n}(z)-\rE_km_{n}(z)\bigr)
\\
\hspace*{-3pt}&&\qquad =\sum_{k=1}^{T}
\frac{1}{n}(\rE_k-\rE_{k-1}) \bigl(\bigl(\rtr\bbA
_{k}^{-1}-\rtr\tilde\bbA_{k}^{-1}
\bigr)+\bigl(\rtr\tilde\bbA_{k}^{-1}-\rtr \bbA
^{-1}\bigr) \bigr)
\\
\hspace*{-3pt}&&\qquad =\sum_{k=1}^{T}\frac{1}{n}(
\rE_k-\rE_{k-1})
\\
\hspace*{-3pt}&&\hspace*{47pt}{}\times  \biggl(\frac{(\bgma
_{k+\tau}+\bgma_{k-\tau})^{*}\bbA_{k}^{-2}\bgma_{k}}{1+(\bgma
_{k+\tau}+\bgma_{k-\tau})^{*}\bbA_{k}^{-1}\bgma_{k}}
+\frac{\bgma_{k}^{*}\tilde\bbA_{k}^{-2}(\bgma_{k+\tau}+\bgma
_{k-\tau})}{1+\bgma_{k}^{*}\tilde\bbA_{k}^{-1}(\bgma_{k+\tau
}+\bgma_{k-\tau})} \biggr)
\\
\hspace*{-3pt}&&\qquad =\sum_{k=1}^{T}\frac{1}{n}(
\rE_k-\rE_{k-1})\frac{d}{ d z} \bigl(\log \bigl(1+(
\bgma_{k+\tau}+\bgma_{k-\tau})^{*}\bbA
_{k}^{-1}\bgma_{k} \bigr)
\\
\hspace*{-3pt}&&\hspace*{128pt}{}+ \log \bigl(1+
\bgma_{k}^{*}\tilde\bbA_{k}^{-1}(
\bgma_{k+\tau}+\bgma _{k-\tau}) \bigr) \bigr)
\\
\hspace*{-3pt}&&\qquad =\sum_{k=1}^{T}\frac{1}{n}(
\rE_k-\rE_{k-1})\frac{d}{ d z}
\\
\hspace*{-3pt}&&\hspace*{47pt}{}\times  \bigl(\log \bigl(
\bigl(1+(\bgma_{k+\tau}+\bgma_{k-\tau})^{*}\bbA
_{k}^{-1}\bgma_{k} \bigr)
\bigl(1+
\bgma_k^*\bbA_k^{-1}(\bgma_{k+\tau
}+
\bgma_{k-\tau}) \bigr)
\\
\hspace*{-3pt}&&\hspace*{127pt}{} -\bgma_k^*\bbA_k^{-1}\bgma_k
\bigl(\bgma_{k+\tau}^*+\bgma_{k-\tau
}^*\bigr)\bbA_k^{-1}(
\bgma_{k+\tau}+\bgma_{k-\tau}) \bigr)
\\
\hspace*{-3pt}&&\hspace*{264pt}{}-\log (x_{n1}-x_{n0}
) \bigr)
\\
\hspace*{-3pt}&&\qquad =\sum_{k=1}^{T}\frac{1}{n}(
\rE_k-\rE_{k-1})
\\
\hspace*{-3pt}&&\hspace*{47pt}{}\times  \biggl(\frac{d}{ d
z}\log \biggl(1+
\frac{\ep_1}{
x_{n1}-x_{n0}}+\frac{\ep_2}{
x_{n1}-x_{n0}}
\\
\hspace*{-3pt}&&\hspace*{68pt}\quad\qquad{} -\frac{\ep_3(\bgma_{k+\tau}^*+\bgma_{k-\tau
}^*)\bbA_k^{-1}(\bgma_{k+\tau}+\bgma_{k-\tau})}{
x_{n1}-x_{n0}}
\\
\hspace*{-3pt}&&\hspace*{68pt}\quad\qquad{} +\frac{\ep_1\ep_2-\ep_4(\bgma_{k+\tau}^*+\bgma
_{k-\tau}^*)\bbA_k^{-1}(\bgma_{k+\tau}+\bgma_{k-\tau})-a_n\ep_5}{
x_{n1}-x_{n0}} \biggr) \biggr)
\\
\hspace*{-3pt}&&\qquad :=\sum_{k=1}^T\frac{1}{n}(
\rE_k-\rE_{k-1})\frac{d}{dz}\log \bigl(1+
\alpha_{k1}(z)+ \alpha_{k2}(z)+\alpha_{k3}(z)+r_k(z)
\bigr)
\\
\hspace*{-3pt}&&\qquad :=\sum_{k=1}^T\frac{1}{n}(
\rE_k-\rE_{k-1})\frac{d}{dz}f_k(z),
\end{eqnarray*}
where $\ep_i$'s, $i=1,\ldots,5$, are defined in (\ref{defep}).

Let $\alpha_{k4}(z):=f_k(z)-\alpha_{k1}(z)-\alpha_{k2}(z)-\alpha
_{k3}(z)-r_k(z)$. It is easy to derive that
%
\begin{eqnarray}
\frac{d}{dz}\alpha_{k1}(z)
&=&\frac{1}{x_{n1}-x_{n0}} \bigl(\bgma_{k+\tau}^*+\bgma_{k-\tau}^*
\bigr)\bbA _k^{-2}\bgma_k
\nonumber\\[-8pt]\label{da1}  \\[-8pt]\nonumber
&&\quad\qquad{} -\frac{x_{n1}'-x_{n0}'}{(x_{n1}-x_{n0})^2}
\bigl(\bgma _{k+\tau}^*+\bgma_{k-\tau}^*\bigr)\bbA_k^{-1}
\bgma_k,\nonumber
\\
\frac{d}{dz}\alpha_{k2}(z)
&=& \frac{1}{x_{n1}-x_{n0}} \bgma_k^*\bbA_k^{-2}(
\bgma_{k+\tau
}+\bgma_{k-\tau})
\nonumber\\[-8pt]\label{da2} \\[-8pt]\nonumber
&&{}-\frac{x_{n1}'-x_{n0}'}{(x_{n1}-x_{n0})^2}
\bgma_{k+}^* \bbA_k^{-1}(\bgma_{k+\tau}+
\bgma_{k-\tau})\nonumber
\end{eqnarray}
and
%
\begin{eqnarray}
\label{da3} &&\frac{d}{dz}\alpha_{k3}(z)
\nonumber
\\
&&\qquad =\frac{1}{x_{n1}-x_{n0}}\nonumber
\\
&&\quad\qquad{}\times  \biggl(\biggl(\bgma_k^*\bbA_k^{-2}
\bgma_k-\frac
{1}{2T}\rtr\bbA_k^{-2}
\biggr) \bigl(\bgma_{k+\tau}^*+\bgma_{k-\tau}^*\bigr)\bbA
_k^{-1}(\bgma_{k+\tau}+\bgma_{k-\tau})
\nonumber\\[-8pt]\\[-8pt]\nonumber
&&\hspace*{20pt}\quad\qquad{} +\biggl(\bgma_k^*\bbA_k^{-1}
\bgma_k-\frac{1}{2T}\rtr\bbA _k^{-1}\biggr)
\bigl(\bgma_{k+\tau}^*+\bgma_{k-\tau}^*\bigr)\bbA_k^{-2}(
\bgma _{k+\tau}+\bgma_{k-\tau}) \biggr)
\nonumber
\\
&&\quad\qquad{}  +\frac{x_{n1}'-x_{n0}'}{(x_{n1}-x_{n0})^2}\nonumber
\\
&&\qquad\qquad{}\times \biggl(\bgma_k^*\bbA _k^{-1}
\bgma_k-\frac{1}{2T}\rtr\bbA_k^{-1}\biggr)
\bigl(\bgma_{k+\tau
}^*+\bgma_{k-\tau}^*\bigr)\bbA_k^{-1}(
\bgma_{k+\tau}+\bgma_{k-\tau}).\nonumber
\end{eqnarray}

Note that by (iii)(b) of Lemma~\ref{x1}, we have $\frac
{1}{\llvert  x_{n1}-x_{n0}\rrvert  }\le K$. Also, by Remarks~\ref{lkj1} and~\ref
{lkj2}, we have $\llvert  x_{n1}'-x_{n0}'\rrvert  =\llvert  {-}\frac
{4a_{n}a'_{n}}{x_{n1}-x_{n0}}\rrvert  \le K$. Together\vspace*{1pt} with Cauchy's formula
and the fact that $\llvert  \ln(1+x)-x\rrvert  \le\llvert  x\rrvert  ^2$ for any complex $x$ with
absolute value smaller than $\frac{1}2$, we have
%
\begin{eqnarray}
\label{da4} &&\biggl\llvert \frac{d}{dz}\alpha_{k4}(z)\biggr
\rrvert\nonumber
\\
&&\qquad =\biggl\llvert \frac{d}{dz} \bigl(\log \bigl(1+\alpha_{k1}(z)+
\alpha _{k2}(z)+\alpha_{k3}(z)+r_k(z) \bigr)\nonumber
\\
&&\hspace*{81pt}{}-
\alpha_{k1}(z)-\alpha _{k2}(z)-\alpha_{k3}(z)-r_k(z)
\bigr)\biggr\rrvert
\nonumber\\[-8pt]\\[-8pt]\nonumber
&&\qquad = \biggl\llvert \frac{1}{2\pi i}\oint_{\llvert  \xi-z\rrvert  =v_n/2}
\bigl(\bigl(\log
\bigl(1+\alpha_{k1}(\xi)+\alpha_{k2}(\xi)+\alpha_{k3}(z)+r_k(\xi)
\bigr)\nonumber
\\
&&\hspace*{140pt}{} -\alpha_{k1}(\xi)-\alpha_{k2}(\xi)-\alpha_{k3}(\xi)-r_k(\xi
)\bigr)\nonumber
\\
&&\hspace*{265pt}{} /
(\xi-z)^2\bigr)\,d
\xi\biggr\rrvert .
\nonumber
\end{eqnarray}
Therefore, for each $u\in[a,b]$, $\ell\ge1$, we have
%
\begin{eqnarray}
\label{dfk} &&\rE\bigl\llvert nv_n\bigl(m_{n}(z)-\rE
m_{n}(z)\bigr)\bigr\rrvert ^{2\ell}
\nonumber
\\
&&\qquad =\rE\Biggl\llvert v_n\sum_{k=1}^{T}(
\rE_k-\rE_{k-1})\frac
{d}{dz}f_k(z)\Biggr
\rrvert ^{2\ell}
\\
&&\qquad \le K\sum_{i=1}^4\rE\Biggl\llvert
v_n\sum_{k=1}^{T}(
\rE_k-\rE _{k-1})\frac{d}{dz}\alpha_{ki}
\Biggr\rrvert ^{2\ell}
+K\rE\Biggl\llvert v_n\sum
_{k=1}^{T}(\rE_k-\rE_{k-1})
\frac{d}{dz}r_{k}\Biggr\rrvert ^{2\ell}.\nonumber\hspace*{-10pt}
\end{eqnarray}
By Lemma~\ref{Burkholder}, for $i=1,2,3,4$, we have
\begin{eqnarray*}
&&\rE\Biggl\llvert v_n\sum_{k=1}^{T}(
\rE_k-\rE_{k-1})\frac{d}{dz}\alpha _{ki}
\Biggr\rrvert ^{2\ell}
\\
&&\qquad \le K_\ell v_n^{2\ell} \Biggl[\rE \Biggl(\sum
_{k=1}^{T}\rE _{k-1}\biggl\llvert
(\rE_k-\rE_{k-1})\frac{d}{dz}\alpha_{ki}
\biggr\rrvert ^2 \Biggr)^\ell
\\
&&\hspace*{89pt}{} +\sum
_{k=1}^{T}\rE\biggl\llvert (\rE_k-
\rE_{k-1})\frac{d}{dz}\alpha _{ki}\biggr\rrvert
^{2\ell} \Biggr]
\\
&&\qquad \le K'_\ell v_n^{2\ell} \Biggl[\rE
\Biggl(\sum_{k=1}^{T}\rE _{k-1}
\biggl\llvert \frac{d}{dz}\alpha_{ki}\biggr\rrvert
^2 \Biggr)^\ell+\sum_{k=1}^{T}
\rE \biggl\llvert \frac{d}{dz}\alpha_{ki}\biggr\rrvert
^{2\ell} \Biggr]. 
\end{eqnarray*}
Now we are ready to estimate the terms above.
By elementary calculation, we have 
%
\begin{eqnarray}
\label{Ak12} \rE_{k}\bigl\llvert \bgma_{k+\tau}^*
\bbA_k^{-1}\bgma_k\bigr\rrvert ^2
&=&\frac{1}{2T}\rE_{k}\bgma_{k+\tau}^*
\bbA_k^{-1}\bigl(\bbA _k^*\bigr)^{-1}
\bgma_{k+\tau}
\nonumber\\[-8pt]\\[-8pt]\nonumber
&\le& \frac{K}T+\frac{1}{2Tv_n^2}\rE_{k}I\bigl(\bigl
\llvert \bgma_{k+\tau}^*\bbA _k^{-1}\bigl(
\bbA_k^*\bigr)^{-1}\bgma_{k+\tau}\bigr\rrvert \ge K
\bigr)
\end{eqnarray}
%
and
%
\begin{eqnarray}
\label{Ak22} \rE_{k}\bigl\llvert \bgma_{k+\tau}^*
\bbA_k^{-2}\bgma_k\bigr\rrvert ^2
&=&\frac{1}{2T}
\rE_k\bgma_{k+\tau}^*\bbA_k^{-2}\bigl(
\bbA_k^*\bigr)^{-2}\bgma _{k+\tau}
\nonumber\\[-8pt]\\[-8pt]\nonumber
&\le&\frac{K}T+\frac{1}{2Tv_n^4}\rE_{k}I\bigl(\bigl
\llvert \bgma_{k+\tau}^*\bbA _k^{-2}\bigl(
\bbA_k^*\bigr)^{-2}\bgma_{k+\tau}\bigr\rrvert \ge K
\bigr),
\end{eqnarray}
for the constant $K>0$ such that Lemmas~\ref{kkbar} and~\ref{kkbar2}
hold.

Come back to the expressions of (\ref{da1}), (\ref{da2}) and (\ref
{da3}). By definition of $x_{ni}$ one can verify that
$x_{n1}'-x_{n0}'=-\frac{4a_{n}a'_{n}}{x_{n1}-x_{n0}}$ which is
bounded.\vspace*{2pt} By Remarks~\ref{lkj1},~\ref{lkj2}, Lemma~\ref{rs} and
estimates (\ref{Ak12}), (\ref{Ak22}), we have
\begin{eqnarray*}
\hspace*{-5pt}&& v_n^{2\ell} \Biggl[\rE \Biggl(\sum
_{k=1}^{T}\rE_{k}\biggl\llvert
\frac
{d}{dz}\alpha_{k1}\biggr\rrvert ^2
\Biggr)^\ell+\sum_{k=1}^{T}\rE
\biggl\llvert \frac
{d}{dz}\alpha_{k1}\biggr\rrvert
^{2\ell} \Biggr]
\\
\hspace*{-5pt}&&\!\qquad \le Kv_n^{2\ell} \Biggl[\rE \Biggl(\sum
_{k=1}^{T}\rE_{k}\bigl\llvert \bigl(\bgma
_{k+\tau}^*+\bgma_{k-\tau}^*\bigr)\bbA_k^{-2}
\bgma_k\bigr\rrvert ^2
\\
\hspace*{-5pt}&&\!\hspace*{75pt}{}+\sum_{k=1}^{T}
\rE_{k}\bigl\llvert \bigl(\bgma_{k+\tau}^*+\bgma_{k-\tau}^*
\bigr)\bbA _k^{-1}\bgma_k\bigr\rrvert
^2 \Biggr)^\ell
\\
\hspace*{-5pt}&&\!\hspace*{28pt}\quad\qquad{} +\sum_{k=1}^{T}\rE\bigl\llvert \bigl(
\bgma_{k+\tau}^*+\bgma_{k-\tau}^*\bigr)\bbA _k^{-2}
\bgma_k\bigr\rrvert ^{2\ell}
\\
\hspace*{-5pt}&&\!\hspace*{26pt}\hspace*{28pt}\quad\qquad{}+\sum_{k=1}^{T}
\rE\bigl\llvert \bigl(\bgma_{k+\tau}^*+\bgma _{k-\tau}^*\bigr)
\bbA_k^{-1}\bgma_k\bigr\rrvert ^{2\ell}
\Biggr]
\\[-2pt]
\hspace*{-5pt}&&\!\qquad \le Kv_n^{2\ell}
\\
\hspace*{-5pt}&&\!\!\quad\qquad{}+Kv_n^{-2\ell}\rE \Bigl(
\max_{k}\rE_kI\bigl(\bigl\llvert (\bgma
_{k+\tau}+\bgma_{k-\tau})^*\bbA_{k}^{-2}\bigl(
\bbA_k^*\bigr)^{-2}(\bgma _{k+\tau}+
\bgma_{k-\tau})\bigr\rrvert \ge K\bigr) \Bigr)^\ell
\nonumber
\\
\hspace*{-5pt}&&\!\!\quad\qquad {} +Kv_n^{2\ell}
\\
\hspace*{-5pt}&&\!\!\quad\qquad{}+K\rE \Bigl(\max_{k}
\rE_kI\bigl(\bigl\llvert (\bgma_{k+\tau}+\bgma
_{k-\tau})^*\bbA_{k}^{-1}\bigl(\bbA_k^*
\bigr)^{-1}(\bgma_{k+\tau}+\bgma _{k-\tau})\bigr\rrvert
\ge K\bigr) \Bigr)^\ell
\nonumber
\\
\hspace*{-5pt}&&\!\!\quad\qquad{} +Kv_n^{2\ell}\bigl(T^{1-\ell}v_n^{-4\ell}+T^{1-\ell}v_n^{-2\ell}
\bigr)
\\
\hspace*{-5pt}&&\!\!\qquad \le Kv_n^{2\ell}
\\
\hspace*{-5pt}&&\!\!\quad\qquad{}+Kv_n^{-2\ell}\sum
_{k=1}^T\rE \bigl(\rE _kI
\bigl(\bigl\llvert (\bgma_{k+\tau}+\bgma_{k-\tau})^*
\bbA_{k}^{-2}\bigl(\bbA _k^*\bigr)^{-2}(
\bgma_{k+\tau}+\bgma_{k-\tau})\bigr\rrvert \ge K\bigr) \bigr)
\nonumber
\\
\hspace*{-5pt}&&\!\!\quad\qquad{} +K\sum_{k=1}^T\rE \bigl(
\rE_kI\bigl(\bigl\llvert (\bgma_{k+\tau}+\bgma_{k-\tau
})^*
\bbA_{k}^{-1}\bigl(\bbA_k^*\bigr)^{-1}(
\bgma_{k+\tau}+\bgma_{k-\tau
})\bigr\rrvert \ge K\bigr) \bigr)
\nonumber
\\
\hspace*{-5pt}&&\!\qquad \le Kv_n^{2\ell},
\end{eqnarray*}
where Lemmas~\ref{kkbar} and~\ref{kkbar2} are used in the last
estimation. By similar arguments, one can show that
\begin{eqnarray*}
&& v_n^{2\ell} \Biggl[\rE \Biggl(\sum
_{k=1}^{T}\rE_{k}\biggl\llvert
\frac
{d}{dz}\alpha_{k2}\biggr\rrvert ^2
\Biggr)^\ell+\sum_{k=1}^{T}\rE
\biggl\llvert \frac
{d}{dz}\alpha_{k2}\biggr\rrvert
^{2\ell} \Biggr]\le Kv_n^{2\ell}.
\end{eqnarray*}
By Remarks~\ref{lkj1},~\ref{lkj2}, (\ref{Ak12}), (\ref{Ak22}) and
Lemmas~\ref{B26} and~\ref{lkj} we have
\begin{eqnarray*}
\hspace*{-4pt}&& v_n^{2\ell} \Biggl[\rE \Biggl(\sum
_{k=1}^{T}\rE_{k}\biggl\llvert
\frac
{d}{dz}\alpha_{k3}\biggr\rrvert ^2
\Biggr)^\ell+\sum_{k=1}^{T}\rE
\biggl\llvert \frac
{d}{dz}\alpha_{k3}\biggr\rrvert
^{2\ell} \Biggr]
\\
\hspace*{-4pt}&&\qquad \le Kv_n^{2\ell} \Biggl[\rE \Biggl(\sum
_{k=1}^{T}\rE_{k}\biggl\llvert \bgma
_k^*\bbA_k^{-2}\bgma_k-
\frac{1}{2T}\rtr\bbA_k^{-2}\biggr\rrvert ^2
\\
\hspace*{-4pt}&&\hspace*{90pt}{}\times
\bigl\llvert (\bgma_{k+\tau}+\bgma_{k-\tau})^*\bbA_k^{-1}(
\bgma_{k+\tau
}+\bgma_{k-\tau})\bigr\rrvert ^2
\\
\hspace*{-4pt}&&\hspace*{41pt}\quad\qquad{} +\sum_{k=1}^{T}\rE_{k}\biggl
\llvert \bgma_k^*\bbA_k^{-1}\bgma_k-
\frac
{1}{2T}\rtr\bbA_k^{-1}\biggr\rrvert
^2
\\
\hspace*{-4pt}&&\hspace*{100pt}{}\times \bigl\llvert (\bgma_{k+\tau}+\bgma_{k-\tau})^*\bbA
_k^{-2}(\bgma_{k+\tau}+\bgma_{k-\tau})\bigr
\rrvert ^2 \Biggr)^\ell
\\
\hspace*{-4pt}&&\hspace*{28pt}\quad\qquad{} +\sum_{k=1}^{T}\rE\biggl\llvert
\bgma_k^*\bbA_k^{-2}\bgma_k-
\frac{1}{2T}\rtr \bbA_k^{-2}\biggr\rrvert
^{2\ell}
\\
\hspace*{-4pt}&&\hspace*{87pt}{}\times \bigl\llvert (\bgma_{k+\tau}+\bgma_{k-\tau})^*\bbA
_k^{-1}(\bgma_{k+\tau}+\bgma_{k-\tau})\bigr
\rrvert ^{2\ell}
\\
\hspace*{-4pt}&&\hspace*{28pt}\quad\qquad{} +\sum_{k=1}^{T}\rE\biggl\llvert
\bgma_k^*\bbA_k^{-1}\bgma_k-
\frac{1}{2T}\rtr \bbA_k^{-1}\biggr\rrvert
^{2\ell}
\\
\hspace*{-4pt}&&\hspace*{107pt}{}\times \bigl\llvert (\bgma_{k+\tau}+\bgma_{k-\tau})^*\bbA
_k^{-2}(\bgma_{k+\tau}+\bgma_{k-\tau})\bigr
\rrvert ^{2\ell} \Biggr]
\\
\hspace*{-4pt}&&\qquad \le Kv_n^{2\ell}\rE \Biggl( \Biggl(\sum
_{k=1}^T\frac{1}{4T^2}\rE _k\rtr
\bbA_k^{-2}\bar\bbA_k^{-2}\bigl\llvert
(\bgma_{k+\tau}+\bgma_{k-\tau
})^*\bbA_k^{-1}(
\bgma_{k+\tau}+\bgma_{k-\tau})\bigr\rrvert ^2
\Biggr)^\ell
\\
\hspace*{-4pt}&&\hspace*{42pt}\quad\qquad{} + \Biggl(\sum_{k=1}^T
\frac{1}{4T^2}\rE_k\rtr\bbA_k^{-1}\bar\bbA
_k^{-1}
\\
\hspace*{-4pt}&&\hspace*{109pt}{}\times \bigl\llvert (\bgma_{k+\tau}+
\bgma_{k-\tau})^*\bbA_k^{-2}(\bgma _{k+\tau}+
\bgma_{k-\tau})\bigr\rrvert ^2 \Biggr)^\ell
\\
\hspace*{-4pt}&&\hspace*{42pt}\quad\qquad{} +\sum_{k=1}^T\rE \biggl(
\frac{1}{4T^2}\rtr\bbA_k^{-2}\bar\bbA
_k^{-2}
\\
\hspace*{-4pt}&&\hspace*{116pt}{}\times \bigl\llvert (\bgma_{k+\tau}+
\bgma_{k-\tau})^*\bbA_k^{-1}(\bgma _{k+\tau}+
\bgma_{k-\tau})\bigr\rrvert ^2\biggr)^\ell
\\
\hspace*{-4pt}&&\hspace*{42pt}\quad\qquad{} +\sum_{k=1}^T\rE \biggl(
\frac{1}{4T^2}\rtr\bbA_k^{-1}\bar\bbA
_k^{-1}
\\
\hspace*{-4pt}&&\hspace*{159pt}{}\times \bigl\llvert (\bgma_{k+\tau}+
\bgma_{k-\tau})^*\bbA_k^{-2}(\bgma _{k+\tau}+
\bgma_{k-\tau})\bigr\rrvert ^2 \biggr)^\ell \Biggr)
\\
\hspace*{-4pt}&&\qquad \le Kv_n^{2\ell}.
\end{eqnarray*}
By (\ref{da4}) and similar arguments, we have
\begin{eqnarray*}
&& v_n^{2\ell} \Biggl[\rE \Biggl(\sum
_{k=1}^{T}\rE_{k}\biggl\llvert
\frac
{d}{dz}\alpha_{k4}\biggr\rrvert ^2
\Biggr)^\ell+\sum_{k=1}^{T}\rE
\biggl\llvert \frac
{d}{dz}\alpha_{k4}\biggr\rrvert
^{2\ell} \Biggr]
\\
&&\qquad \le Kv_n^{2\ell} \Biggl[\rE \Biggl(\frac{1}{v_n^2}\sup
_{\llvert  \xi
-z\rrvert  =v_n/2}\sum_{k=1}^{T}
\rE_{k} \bigl(\bigl\llvert \alpha_{k1}(\xi)\bigr\rrvert
^4+\bigl\llvert \alpha _{k2}(\xi)\bigr\rrvert
^4
\\
&&\hspace*{90pt}\hspace*{76pt}{}+\bigl\llvert \alpha_{k3}(\xi)\bigr\rrvert
^4+\bigl\llvert r_k(\xi)\bigr\rrvert ^4
\bigr) \Biggr)^\ell
\\
&&\hspace*{28pt}\quad\qquad{} +\frac{1}{v_n^{2\ell}}\sup_{\llvert  \xi-z\rrvert  =v_n/2}\sum
_{k=1}^{T}\rE \bigl(\bigl\llvert
\alpha_{k1}(\xi)\bigr\rrvert ^{4\ell}+\bigl\llvert
\alpha_{k2}(\xi)\bigr\rrvert ^{4\ell}
\\
&&\hspace*{169pt}{}+\bigl\llvert \alpha
_{k3}(\xi)\bigr\rrvert ^{4\ell}+\bigl\llvert r_k(
\xi)\bigr\rrvert ^{4\ell} \bigr) \Biggr]
\\
&&\qquad \le KT^{-\ell}v_n^{-4\ell}.
\end{eqnarray*}
Finally, by measurable properties of some terms of $r_k$, we have
\begin{eqnarray*}
(\rE_{k-1}-\rE_k)r_k=(\rE_{k-1}-
\rE_k)\frac{\ep_1\ep_2}{x_{n1}-x_{n0}}
\end{eqnarray*}
from which and similar argument for $\alpha_{k1}$ and $\alpha_{k2}$,
we conclude that
\begin{eqnarray*}
v_n^{2\ell}\rE\Biggl\llvert \frac{d}{dz}\sum
_{k=1}^T(\rE_{k-1}-\rE
_k)r_k\Biggr\rrvert ^{2\ell}=KT^{-\ell}v_n^{-4\ell}.
\end{eqnarray*}

Substituting the five upper-bounds into (\ref{dfk}), we have
\begin{eqnarray*}
&&\rP \Bigl(\max_{u\in S_n}\bigl\llvert nv_n
\bigl(m_n(z)-\rE m_n(z)\bigr)\bigr\rrvert >\ep \Bigr)
\\
&&\qquad = Kn^2\rE\bigl\llvert nv_n\bigl(m_{n}(z)-
\rE m_{n}(z)\bigr)\bigr\rrvert ^{2\ell}
\\
&&\qquad \le Kn^2\bigl(v_n^{2\ell}+v_n^{-4\ell}T^{-\ell}
\bigr)
\end{eqnarray*}
which is summable when $\ell>318$ and $v_n\ge n^{-\alpha}$ for
$\alpha= 1/212$. Therefore, we have
proved that $\max_{u\in[a,b]}\llvert  m_n(z)-\rE m_n(z)\rrvert  =o(\frac{1}{nv_n})$ a.s.

\subsection{A refined convergence rate of $\mathrm{E}m_{n}(z)-m_n^0(z)$}\label{sec5.2}
To show
\[
\sup_{u\in[a,b]}\bigl\llvert  \rE m_{n}(z)-m_n^0(z)\bigr\rrvert  =o\biggl(\frac{1}{nv_n}\biggr),
\]
we follow the notation and expressions in Section~\ref{prev}.
Recall
%
\begin{eqnarray}
&&c_n+c_nz\rE m_n(z)
\nonumber
\\
&&\qquad =\frac{1}T\sum_{k=1}^T
\biggl[1-\rE\frac{1}{1+\bgma_k^*\widetilde
\bbA_k^{-1}(\bgma_{k+\tau}+\bgma_{k-\tau})} \biggr]
\nonumber
\\
&&\qquad =\frac{1}T\sum_{k=1}^T
\biggl[1-\rE \biggl(1\Big/\biggl(1+\bgma_k^*\bbA
_k^{-1}(\bgma_{k+\tau}+\bgma_{k-\tau})
\nonumber\\[-8pt]\\[-8pt]\nonumber
&&\hspace*{117pt}{}-\frac{\bgma_k^*\bbA_k^{-1}
\bgma_k(\bgma_{k+\tau}^*+\bgma_{k-\tau}^*)\bbA_k^{-1}(\bgma
_{k+\tau}+\bgma_{k-\tau})}{1+(\bgma_{k+\tau}^*+\bgma_{k-\tau
}^*)\bbA_k^{-1}\bgma_{k}}\biggr)\biggr) \biggr]
\nonumber
\\
&&\qquad =1-\frac{1}{x_{n1}-x_{n0}}+\delta_n,\nonumber
\end{eqnarray}
where
\begin{eqnarray*}
&&\delta_n=\frac{1}T\sum_{k=1}^T
\rE\eta_{k}
\end{eqnarray*}
with
\begin{eqnarray*}
\eta_k&=&- \biggl(1\Big/ \biggl(1+\bgma_k^*\bbA_k^{-1}(\bgma_{k+\tau
}+\bgma_{k-\tau})
\\
&&\hspace*{34pt}{} -\frac{\bgma_k^*\bbA_k^{-1}
\bgma_k(\bgma_{k+\tau}^*+\bgma_{k-\tau}^*)\bbA_k^{-1}(\bgma
_{k+\tau}+\bgma_{k-\tau})}{1+(\bgma_{k+\tau}^*+\bgma_{k-\tau
}^*)\bbA_k^{-1}\bgma_{k}}\biggr)-\frac{1}{x_{n1}-x_{n0}} \biggr).
\end{eqnarray*}

Consider expressions of (\ref{Emn}) and (\ref{m}). To apply Lemma
\ref{P2}, we only need to show $\llvert  \delta_n\rrvert  =o(\frac{1}{nv_n})$, which
can be reduced to showing $\llvert  \rE\eta_k\rrvert  =o(\frac{1}{nv_n})$ for $\log
^{2}n<k<T-\log^{2}n$ and $\llvert  \rE\eta_{k}\rrvert  =O(1)$ for $k\le\log^{2}n$
or $\ge T-\log^{2}n$.

When $\log^{2}n<k<T-\log^{2}n$, rewrite $\eta_k$ as
\begin{eqnarray*}
-\eta_k
&=& 1\Big/ \biggl(1+\bgma_k^*\bbA_k^{-1}(\bgma_{k+\tau}+\bgma_{k-\tau
})
\\
&&\hspace*{18pt}{} -\frac{\bgma_k^*\bbA_k^{-1}
\bgma_k(\bgma_{k+\tau}^*+\bgma_{k-\tau}^*)\bbA_k^{-1}(\bgma
_{k+\tau}+\bgma_{k-\tau})}{1+(\bgma_{k+\tau}^*+\bgma_{k-\tau
}^*)\bbA_k^{-1}\bgma_{k}}\biggr)-\frac{1}{1-(\sfrac{2a_n^2}{x_{n1}})}
\\
&=&\bigl(1+\bigl(\bgma_{k+\tau}^*+\bgma_{k-\tau}^*\bigr)\bbA_k^{-1}\bgma
_{k}\bigr)
\\
&&{} /
\bigl(\bigl(1+\bgma_k^*\bbA_k^{-1}(\bgma_{k+\tau}+\bgma_{k-\tau
})\bigr)\bigl(1+\bigl(\bgma_{k+\tau}^*+\bgma_{k-\tau}^*\bigr)\bbA_k^{-1}\bgma
_{k}\bigr)
\\
&&\hspace*{49pt}{}-\bgma_k^*\bbA_k^{-1}
\bgma_k\bigl(\bgma_{k+\tau}^*+\bgma_{k-\tau}^*\bigr)\bbA_k^{-1}(\bgma
_{k+\tau}+\bgma_{k-\tau})\bigr)
\\
&&{}-\frac{1}{1-\sklfrac{2a_n^2}{x_{n1}}}
\\
&=&(1+\ep_1)
\\
&&{} \Big/
\biggl(1+\ep_1+\ep_2+\ep_1\ep_2
\\
&&\hspace*{13pt}{} -\bigl(\bgma_{k+\tau
}^*+\bgma_{k-\tau}^*\bigr)\bbA_k^{-1}(\bgma_{k+\tau}+\bgma_{k-\tau
})(\ep_3+\ep_4)-a_n\ep_5-\frac{2a_n^2}{x_{n1}}\biggr)
\\
&&{} -\frac{1}{1-\sklfrac{2a_n^2}{x_{n1}}}
\\
&=&\frac{1}{1-\sklfrac{2a_n^2}{x_{n1}}}
\\
&&{} \times
\biggl(-\ep_1\frac{2a_n^2}{x_{n1}}-\ep_2-\ep_1\ep_2
\\
&&\hspace*{18pt}{}+\bigl(\bgma
_{k+\tau}^*+\bgma_{k-\tau}^*\bigr)\bbA_k^{-1}(\bgma_{k+\tau}+\bgma
_{k-\tau})(\ep_3+\ep_4)+a_n\ep_5\biggr)
\\
&&\hspace*{14pt}{}\Big/
\biggl(1+\ep_1+\ep_2+\ep_1\ep
_2
\\
&&\hspace*{28pt}{}-\bigl(\bgma_{k+\tau}^*+\bgma_{k-\tau}^*\bigr)\bbA_k^{-1}(\bgma_{k+\tau
}+\bgma_{k-\tau})(\ep_3+\ep_4)-a_n\ep_5-\frac{2a_n^2}{x_{n1}}\biggr),
\end{eqnarray*}
where $\ep_i$'s are defined as in Section~\ref{prev}.

For simplicity, denote $\tilde\ep=\ep_2+\ep_1\ep_2-(\bgma_{k+\tau
}^*+\bgma_{k-\tau}^*)\bbA_k^{-1}(\bgma_{k+\tau}+\bgma_{k-\tau
})(\ep_3+\ep_4)-a_n\ep_5$.
Applying the identity $\frac{x}{1+x+y}=\frac{x}{1+y}-\frac
{x^2}{(1+x+y)(1+y)}$ repeatedly, we have
\begin{eqnarray*}
-\eta_k&=&\frac{1}{1-\sklfrac{2a_n^2}{x_{n1}}}\times\frac{-\ep_1\sklfrac
{2a_n^2}{x_{n1}}-\tilde\ep}{1+\ep_1+\tilde\ep-\sklfrac
{2a_n^2}{x_{n1}}}
\\[-1pt]
&=&-\frac{\sfrac{2a_n^2}{x_{n1}}}{1-\sklfrac{2a_n^2}{x_{n1}}}\times \frac{\ep_1+\tilde\ep}{1+\ep_1+\tilde\ep-\sklfrac{2a_n^2}{x_{n1}}} -\frac{\tilde\ep}{1+\ep_1+\tilde\ep-\sklfrac{2a_n^2}{x_{n1}}}
\\[-1pt]
&=&-\frac{\sfrac{2a_n^2}{x_{n1}}}{1-\sklfrac{2a_n^2}{x_{n1}}}
\\[-1pt]
&&\hspace*{7pt}{}\times  \biggl(\frac{\ep_1+\tilde\ep}{1-\sklfrac{2a_n^2}{x_{n1}}} -\frac{(\ep_1+\tilde\ep)^2}{ (1-\sklfrac{2a_n^2}{x_{n1}}
) (1+\ep_1+\tilde\ep-\sklfrac{2a_n^2}{x_{n1}} )} \biggr)
\\[-1pt]
&&{}- \biggl(\frac{\tilde\ep}{1+\ep_1-\sklfrac{2a_n^2}{x_{n1}}}
\\[-1pt]
&&\hspace*{18pt}{}-\frac{\tilde\ep^2}{ (1+\ep_1-\sklfrac{2a_n^2}{x_{n1}} )
(1+\ep_1+\tilde\ep-\sklfrac{2a_n^2}{x_{n1}} )} \biggr)
\\[-1pt]
&=&-\frac{\sfrac{2a_n^2}{x_{n1}}}{1-\sklfrac{2a_n^2}{x_{n1}}}
\\[-1pt]
&&\hspace*{7pt}{}\times\biggl(\frac{\ep_1+\tilde\ep}{1-\sklfrac{2a_n^2}{x_{n1}}} -\frac{(\ep_1+\tilde\ep)^2}{ (1-\sklfrac{2a_n^2}{x_{n1}}
) (1+\ep_1+\tilde\ep-\sklfrac{2a_n^2}{x_{n1}} )} \biggr)
\\[-1pt]
&&{}- \biggl(\frac{\tilde\ep}{1-\sklfrac{2a_n^2}{x_{n1}}}-\frac{\tilde
\ep\ep_1}{ (1+\ep_1-\sklfrac{2a_n^2}{x_{n1}} ) (1-\sklfrac
{2a_n^2}{x_{n1}} )} \biggr)
\\[-1pt]
&&{} +\frac{\tilde\ep^2}{ (1+\ep_1-\sklfrac{2a_n^2}{x_{n1}} )
(1+\ep_1+\tilde\ep-\sklfrac{2a_n^2}{x_{n1}} )}.
\end{eqnarray*}

Therefore, by Lemma~\ref{x1}(iv)(b), we have $\llvert  {-}\frac{\sfrac
{2a_n^2}{x_{n1}}}{1-\sklfrac{2a_n^2}{x_{n1}}}\rrvert  =\llvert  \frac
{2x_{n0}}{x_{n1}-x_{n0}}\rrvert  \le\llvert  \frac
{2x_{n1}}{x_{n1}-x_{n0}}\rrvert  $ is bounded. Together with the fact
that all the denominators being bounded below and the Cauchy--Schwarz
inequality, to show $\llvert  \rE\eta_k\rrvert  =o(\frac{1}{nv_n})$, it suffices to
show $\llvert  \rE\ep_1\rrvert , \llvert  \rE\tilde\ep\rrvert , \llvert  \rE\ep_1^2\rrvert , \llvert  \rE\tilde\ep
^2\rrvert  $ are of $o(\frac{1}{nv_n})$. As $\llvert  \rE\ep_i\rrvert  =0$ for $i=1,2,3$,
it is clear that the above convergence rates achieve $o(\frac
{1}{nv_n})$ provided that so do $\rE\llvert  \ep_i\rrvert  ^2, i=1,2,3,4,5$, $\llvert  \rE
\ep_4\rrvert  $ and $\llvert  \rE\ep_5\rrvert  $ for $\log^{2}n<k<T-\log^{2}n$.

When $\log^{2}n<k<T-\log^{2}n$, for $i=1$, by Lemma~\ref{kkbar}, we
have, for any $t>0$,
\begin{eqnarray*}
\rE\bigl\llvert (\bgma_{k+\tau}+\bgma_{k-\tau})^*
\bbA_k^{-1}\bgma_k\bigr\rrvert ^2
&=&\frac{1}{2T}\rE(\bgma_{k+\tau}+\bgma_{k-\tau})^*\bbA
_k^{-1}\bigl(\bbA_k^*\bigr)^{-1}(
\bgma_{k+\tau}+\bgma_{k-\tau})
\\[-1pt]
&=&\frac{K}T+v_{n}^{{-2}}o\bigl(n^{-t}
\bigr)=O(1/n)=o\biggl(\frac{1}{nv_n}\biggr).
\end{eqnarray*}
Similarly, for $i=2$, $\rE\llvert  \ep_2\rrvert  ^2=O(1/n)=o(\frac{1}{nv_n})$.

For $i=3$, by Lemmas~\ref{B26} and~\ref{lkj}, we have
\begin{eqnarray*}
\rE\llvert \ep_3\rrvert ^2&=&\rE\biggl\llvert
\bgma_k^*\bbA_k^{-1}\bgma_k-
\frac{1}{2T}\rtr \bbA_k^{-1}\biggr\rrvert
^2\le\frac{K}{4T^2}\rE\bigl\llvert \rtr\bbA_k^{-1}
\bigl(\bbA ^*_k\bigr)^{-1}\bigr\rrvert
\\
&=&\frac{K}{4T^2}
\rE\sum\frac{1}{\llvert  \lambda_{kj}-z\rrvert  ^2}
\\
&\le& \frac{K}{2T} +\frac{K}{Tv_n^2}F_n\bigl(
\bigl[a',b'\bigr]\bigr)\le\frac{K}{T}+o
\bigl(T^{-1}\bigr) =O(1/n)=o\biggl(\frac{1}{nv_n}\biggr).
\end{eqnarray*}
For $\llvert  \rE\ep_4\rrvert  $, by Lemma~\ref{ep40} we have
\begin{eqnarray*}
\llvert \rE\ep_4\rrvert &=&\biggl\llvert \frac{1}{2T}\rE\rtr
\bbA_k^{-1}-a_n\biggr\rrvert =
\frac
{1}{2T}\bigl\llvert \rE\bigl(\rtr\bbA_k^{-1}-
\rtr\bbA^{-1}\bigr)\bigr\rrvert =O\bigl(T^{-1}\bigr)=o\biggl(
\frac{1}{nv_n}\biggr).
\end{eqnarray*}
For $\rE\llvert  \ep_4\rrvert  ^2$, by (\ref{eqlacing1}) and the convergence rate
obtained in Section~\ref{sec42}, we have
\begin{eqnarray*}
&&\rE\biggl\llvert \frac{1}{2T}\rtr\bbA_k^{-1}-a_n
\biggr\rrvert ^2
\\
&&\qquad \le 2\rE\biggl\llvert \frac{1}{2T}\rtr\bbA_k^{-1}-
\rE\frac{1}{2T}\rtr\bbA _{k}^{-1}\biggr\rrvert
^2+2\biggl\llvert \frac{1}{2T}\rE\rtr\bbA_{k}^{-1}-a_n
\biggr\rrvert ^2
\\
&&\qquad \le \frac{K}{n^2v_n^2}+O\bigl(n^{-1}\bigr)=o
\biggl(\frac{1}{nv_n} \biggr).
\end{eqnarray*}
Bounds of $\llvert  \rE\ep_5\rrvert  $ and $\rE\llvert  \ep_5\rrvert  ^2$ will follow Lemmas~\ref{ep41}(b2), (b3) and~\ref{ep44}(b1), (b2).

To show $\llvert  \rE\eta_{k}\rrvert  =O(1)$ when $k\le\log^{2}n$ or $\ge T-\log
^{2}n$, we just prove the case for $k\ge T-\log^{2}n$, as the case for
$k\le\log^{2}n$ follows by symmetry.

When $k\ge T-\log^{2}n$, by Lemma~\ref{Ek+tau}(b1), we have
$\rP(\llvert  \bgma_{k+\tau}^*\bbA_k^{-1}\bgma_{k+\tau}\rrvert  \ge1-\eta)=o(n^{-t})$.
By Lemma~\ref{Ek+tau}(a), we have $\rP(\llvert  \bgma_{k-\tau}^*\bbA
_k^{-1}\bgma_{k-\tau}-\frac{c_n\rE m_n}{2x_{n1}}\rrvert  \ge
v_n^6)=o(n^{-t})$, by Lemma~\ref{rs},
$\rP(\llvert  \bgma_k^*\bbA_k^{-1}\bgma_{k\pm\tau}\rrvert  \ge v_n^3)=o(n^{-t})$,
and by\vspace*{1pt} Lemmas~\ref{B26} and inequalities (\ref{eqlacing1}) and (\ref{W1}),
$\rP(\llvert  \bgma_k^*\bbA_k^{-1}\bgma_k-a_n\rrvert  \ge v_n^3)=o(n^{-t})$. By
Lemma~\ref{cross}(a),
$\rP(\llvert  \bgma_{k\pm\tau}^*\bbA_{k}^{-1}\bgma_{k\mp\tau}\rrvert  \ge
v_n^6)=o(n^{-t})$.
By Lemma~\ref{x1}(ii)(b) and (iv)(b), we have
$\llvert  \frac{1}{x_{n1}-x_{n0}}\rrvert  \le K$ and $\llvert  \rE\eta_k\rrvert  \le Kv_n^{-1}$.
Substitute\vspace*{2pt} the above results into the definition of $\eta_k$, and we
finally have
\begin{eqnarray*}
\llvert \rE\eta_k\rrvert &\le&\biggl\llvert \rE
\biggl(1\Big/ \biggl(1+\bgma_k^*\bbA_k^{-1}(\bgma
_{k+\tau}+\bgma_{k-\tau})
\\
&&\hspace*{36pt}{}-\frac{\bgma_k^*\bbA_k^{-1}
\bgma_k(\bgma_{k+\tau}^*+\bgma_{k-\tau}^*)\bbA_k^{-1}(\bgma
_{k+\tau}+\bgma_{k-\tau})}{1+(\bgma_{k+\tau}^*+\bgma_{k-\tau
}^*)\bbA_k^{-1}\bgma_{k}}\biggr)\biggr)
\biggr\rrvert
\\
&&{}+\biggl\llvert \frac{1}{x_{n1}-x_{n0}} \biggr\rrvert
\\
&\le&\biggl\llvert \frac{1+v_n^3}{(1-2v_n^3)-(\sfrac{1}2-\eta+v_n^3)(1-\eta
+3v_n^3+\sfrac{\llvert  a_{n}\rrvert  }{\llvert  x_{n1}\rrvert  })}\biggr\rrvert
\\
&&{} +K+Kv_n^{-1}o
\bigl(n^{-t}\bigr)
=O(1).
\end{eqnarray*}

\section{Completing the proof}\label{sec6}\label{sec5}

In this section, we follow the idea of \citet{BaiSil98} and
give the main steps here. From what has been obtained in the last two
sections, we have, with $v_n=n^{-\sfrac{1}{212}}$,
%
\begin{eqnarray}
\label{conclu} \sup_{u\in[a,b]}\bigl\llvert m_n(z)-m_n^0(z)
\bigr\rrvert =o\biggl(\frac{1}{nv_n}\biggr) \qquad\mbox{a.s.}
\end{eqnarray}
It is clear from the last two sections that (\ref{conclu}) is true
when $\Im(z)$ is replaced by a constant multiple of $v_n$. In fact, we have
\begin{eqnarray*}
\max_{k\in\{1,2,\ldots,106\}}\sup_{u\in[a,b]}\bigl\llvert
m_n(u+i\sqrt {k}v_n)-m_n^0(u+i
\sqrt{k}v_n)\bigr\rrvert =o\bigl(v_n^{211}\bigr)
\qquad\mbox{a.s.}
\end{eqnarray*}
Taking the imaginary part, we get
\begin{eqnarray*}
\max_{k\in\{1,2,\ldots,106\}}\sup_{u\in[a,b]}\biggl\llvert \int
\frac
{d(F_n(\lambda)-F_n^0(\lambda))}{(u-\lambda)^2+kv_n^2} \biggr\rrvert =o\bigl(v_n^{210}\bigr)
\qquad\mbox{a.s.}
\end{eqnarray*}
After taking difference, we obtain
\begin{eqnarray}
\max_{k_1\ne k_2}\sup_{u\in[a,b]}\biggl\llvert \int
\frac
{v_n^2d(F_n(\lambda)-F_n^0(\lambda))}{((u-\lambda
)^2+k_1v_n^2)((u-\lambda)^2+k_2v_n^2)}\biggr\rrvert &=&o\bigl(v_n^{210}\bigr)\nonumber
\\[-2pt]
\eqntext{\mbox {a.s.}}
\\
&\vdots&\nonumber
\\
\sup_{u\in[a,b]}\biggl\llvert \int\frac{(v_n^2)^{105}d(F_n(\lambda
)-F_n^0(\lambda))}{((u-\lambda)^2+v_n^2)((u-\lambda)^2+2v_n^2)\cdots
((u-\lambda)^2+106v_n^2)}\biggr\rrvert
&=&o\bigl(v_n^{210}\bigr)\nonumber
\\
\eqntext{\mbox{a.s.}}
\end{eqnarray}
Therefore,
\begin{eqnarray}
\sup_{u\in[a,b]}\biggl\llvert \int\frac{d(F_n(\lambda)-F_n^0(\lambda
))}{((u-\lambda)^2+v_n^2)((u-\lambda)^2+2v_n^2)\cdots((u-\lambda
)^2+106v_n^2)}\biggr\rrvert
=o(1)\nonumber
\\[-2pt]
\eqntext{\mbox{a.s.}}
\end{eqnarray}
After splitting the integral, we get
\begin{eqnarray*}
&&\sup_{u\in[a,b]}\biggl\llvert \int\frac{I_{[a',b']^c}(\lambda
)d(F_n(\lambda)-F_n^0(\lambda))}{((u-\lambda)^2+v_n^2)((u-\lambda
)^2+2v_n^2)\cdots((u-\lambda)^2+106v_n^2)}
\\
&&\hspace*{6pt}\qquad{} +\sum
_{\lambda_j\in[a',b']}\frac{v_n^{212}}{((u-\lambda
_j)^2+v_n^2)((u-\lambda_j)^2+2v_n^2)\cdots((u-\lambda
_j)^2+106v_n^2)}\biggr\rrvert
\\
&&\qquad  =o(1) \qquad
\mbox{a.s.}
\end{eqnarray*}
Note that the first term tends to 0 by dominated convergence theorem.
Now, if there is at least one eigenvalue contained in $[a,b]$, then the
second sum will be away from zero when $u$ takes one of such
eigenvalues. This contradicts the right-hand side. Therefore, with
probability 1, there are no eigenvalues of $\bbM_n$ in $[a,b]$ for all
$n$ large and the proof is complete.

\begin{appendix}\label{append}
\section{Justification of truncation, centralization and rescaling}\label{sec7}

Here, we give some justifications of (\ref{trun}), which will be
divided into two parts.

\subsection{Truncation and centralization}\label{sec7.1}
Fix some $C>0$, define $\hat{\ep}_{it}=\break  \ep_{it}I_{\{\llvert  x_{it}\rrvert  \le C\}
}-\rE\ep_{it}I_{\{\llvert  x_{it}\rrvert  \le C\}}$, $\hat{\bgma}_k=\frac{1}{\sqrt
{2T}}(\hat{\ep}_{1k},\ldots,\hat{\ep}_{nk})' \equiv\frac
{1}{\sqrt{2T}}\hat\bbe_k$,
$\hat{\mathbf{E}}=(\hat\bbe_1,\ldots,\hat\bbe_T)$, $\hat
{\mathbf{E}}_{\tau}=(\hat\bbe_{1+\tau},\ldots,\hat\bbe_{T+\tau})$
and $\hat\bbM_n=\sum_{k=1}^T(\hat{\bgma}_k\hat{\bgma}_{k+\tau
}^*+\hat{\bgma}_{k+\tau}\hat{\bgma}_k^*)
=\frac{1}{2T}(\hat{\mathbf{E}}\hat{\mathbf{E}}_{\tau}^*+\hat
{\mathbf{E}}_{\tau}\hat{\mathbf{E}}^*)$. By Theorem A.46 of \citet{BaiSil10},
\begin{eqnarray*}
&&\max_k\bigl\llvert \la_k(\hat{
\bbM}_n)-\la_k(\bbM_n)\bigr\rrvert
\\[-1.5pt]
&&\qquad \leq\llVert
\hat{\bbM}_n-\bbM_n\rrVert
\\[-1.5pt]
&&\qquad =\frac{1}{2T}\bigl\llVert (\mathbf{E}-\hat{\mathbf{E}})\hat{\mathbf
{E}}_{\tau}^*+\hat{\mathbf{E}}_{\tau}(\mathbf{E}-\hat{
\mathbf{E}})^*+ \mathbf{E}(\mathbf{E}_{\tau}-\hat{\mathbf{E}}_{\tau})^*+(
\mathbf{E}_{\tau}-\hat{\mathbf{E}}_{\tau})\mathbf{E}^*\bigr
\rrVert
\\[-1.5pt]
&&\qquad \le\frac{1}{T} \bigl(\llVert \mathbf{E}-\hat{\mathbf{E}}\rrVert \llVert
\hat {\mathbf{E}}_{\tau}\rrVert +\llVert \mathbf{E}-\hat{\mathbf{E}}
\rrVert \llVert \mathbf{E}\rrVert \bigr).
\end{eqnarray*}
By a similar approach as in \citet{YinBaiKri88}, one can show that almost surely
\begin{eqnarray*}
\limsup_n\frac{1}{\sqrt T}\llVert \mathbf{E}\rrVert
&\le&(1+\sqrt c)^2,
\\[-1.5pt]
\limsup_n\frac{1}{\sqrt T}\llVert \hat{
\mathbf{E}}_{\tau}\rrVert &\le&(1+\sqrt c)^2
\end{eqnarray*}
and
\begin{eqnarray*}
&&\limsup_n\frac{1}{\sqrt T}\llVert \mathbf{E}-\hat{
\mathbf{E}}\rrVert
\\[-1.5pt]
&&\qquad \le(1+\sqrt c)^2\max_{i,t}\operatorname{var}(
\ep_{it}-\hat{\ep }_{it})
\\[-1.5pt]
&&\qquad =(1+\sqrt c)^2\max_{i,t}\operatorname{var}(
\ep_{it}I_{\{\llvert  x_{it}\rrvert  \ge C\}})
\\[-1.5pt]
&&\qquad \le(1+\sqrt c)^2\max_{i,t}\rE(
\ep_{it}I_{\{\llvert  x_{it}\rrvert  \ge C\}})^2
\\[-1.5pt]
&&\qquad \le\frac{(1+\sqrt c)^2}{C^2}\max_{i,t}\rE\ep_{it}^4
\\[-1.5pt]
&&\qquad \le\frac{(1+\sqrt c)^2M}{C^2},
\end{eqnarray*}
which can be arbitrarily small by choosing $C$ large enough.
This verifies the truncation at a fixed point and centralization.
\subsection{Rescaling}\label{sec7.2}
Define $\sigma_{it}^2=\rE\llvert  \hat\ep_{it}\rrvert  ^2$, $\check\ep_{it}=\hat
\ep_{it}/\sigma_{it}$, $\check{\bgma}_k=\frac{1}{\sqrt{2T}}(\check
{\ep}_{1k},\ldots, \break \check{\ep}_{nk})' \equiv\frac{1}{\sqrt
{2T}}\check\bbe_k$,
$\check{\mathbf{E}}=(\check\bbe_1,\ldots,\check\bbe_T)$, $\check
{\mathbf{E}}_{\tau}=(\check\bbe_{1+\tau},\ldots,\check\bbe
_{T+\tau})$,
$\mathbf{D}= (\sigma_{it}^{-1} )_{n\times T}$, $\mathbf{D}_{\tau}= \break (\sigma_{i(t+\tau)}^{-1} )_{n\times T}$ and
$\check\bbM_n=\sum_{k=1}^T(\check{\bgma}_k\check{\bgma}_{k+\tau
}^*+\check{\bgma}_{k+\tau}\check{\bgma}_k^*)
=\frac{1}{2T}(\check{\mathbf{E}}\check{\mathbf{E}}_{\tau
}^*+\check{\mathbf{E}}_{\tau}\check{\mathbf{E}}^*)$. By Theorem
A.46 and Corollary A.21 of \citet{BaiSil10},
\begin{eqnarray*}
&&\max_k\bigl\llvert \la_k(\check{
\bbM}_\tau)-\la_k(\hat{\bbM}_\tau)\bigr\rrvert
\\
&&\qquad \le \llVert \check{\bbM}_\tau-\hat{\bbM}_\tau\rrVert
\\
&&\qquad \le\frac{1}T\bigl\llVert \hat{\mathbf{E}}\circ(\mathbf{D}-\mathbf{J})\bigr\rrVert \bigl\llVert \hat{\mathbf{E}}_{\tau}\circ(
\mathbf{D}_\tau-\mathbf{J})\bigr\rrVert
\\
&&\qquad \le\frac{1}T\llVert \hat{\mathbf{E}}\rrVert \llVert \hat{
\mathbf{E}}_{\tau}\rrVert \max_{i,t}\bigl(
\sigma_{it}^{-1}-1\bigr)^2.
\end{eqnarray*}
Here, $\circ$ denotes the Hadamard product and $\mathbf{J}$ is the
$n\times T$ matrix of all entries 1.

From \citet{YinBaiKri88}, we have, with probability 1 that
$\limsup_n\frac{1}T\llVert  \hat{\mathbf{E}}\rrVert  \llVert  \hat{\mathbf{E}}_{\tau
}\rrVert  \le(1+\sqrt c)^4$.

Also, we have
\begin{eqnarray*}
\max_{i,t}\bigl\llvert 1-\sigma_{it}^2
\bigr\rrvert &\le&\max_{i,t} \bigl(\rE\llvert
\ep_{it}\rrvert ^2I\bigl(\llvert \ep_{it}\rrvert
>C\bigr)+ \bigl(\rE\llvert \ep _{it}\rrvert I\bigl(\llvert
\ep_{it}\rrvert >C\bigr) \bigr)^2 \bigr)
\\
&\le&\max_{i,t}\frac{2}{C^2}\rE\llvert
\ep_{it}\rrvert ^4\le\frac{2M}{C^2}\to 0\qquad\mbox{as
} C\to\infty.
\end{eqnarray*}
Since $\min_{i,t}\sigma_{it}\to1$ as $n\to\infty$ and thus $\sigma
_{it}(1+\sigma_{it})\ge1$ for all large $n$. Therefore, we have
\begin{eqnarray*}
\sigma_{it}^{-1}-1=\frac{1-\sigma_{it}^2}{\sigma_{it}(1+\sigma
_{it})}\le1-
\sigma_{it}^2,
\end{eqnarray*}
which implies
$\max_k\llvert  \la_k(\check{\bbM}_\tau)-\la_k(\hat{\bbM}_\tau)\rrvert  \to0$
as $n\to\infty$.
\section{Proofs of lemmas in 
Section~\texorpdfstring{\lowercase{\protect\ref{sec25}}}{3}}\label{sec8}

\subsection{Proofs of Lemmas \texorpdfstring{\protect\ref{P1}}{3.1}, 
\texorpdfstring{\protect\ref{P2}}{3.2} 
and \texorpdfstring{\protect\ref{P3}}{3.3}}\label{sec8.1}
To show Lemma~\ref{P1}, take $d=\sqrt{\frac{1}{2m}}$ and denote $S$
the total area covered by the $m$ balls $B(x_i,dr_n)$, $i=1,\ldots,m$.
Then we have $S\le m\pi(dr_n)^2<\pi r_n^2$, which is the total area of
$B(x_0,r_n)$. Therefore, such $x$ must exist.

For Lemma~\ref{P2}, write $P_n(x)=\prod_{j=1}^{k}(x-x_{nj})$ and
$P(x)=\prod_{j=1}^{m}(x-x_{j})^{\ell_j}$.
Let
\[
\delta=\frac{1}{3}\mathop{\min_{i,j\in\{1,\ldots,m\}}}_{i\ne j}
\llvert x_i-x_j\rrvert >0.
\]
First, we claim that for any $i\in{\{1,\ldots,k\}}$, there exists
$j\in{\{1,\ldots,m\}}$ such that $x_{ni}\in B(x_j,\delta)$. Suppose
not, that is, there is some $x_{ni}$ with $\llvert  x_{ni}-x_{j}\rrvert  \ge\delta$
for any $j\in{\{1,\ldots,m\}}$. Then it follows that
$\llvert  P(x_{ni})\rrvert  =\prod_{j=1}^{m}\llvert  x_{ni}-x_{j}\rrvert  ^{\ell_j}\ge\delta^k$. On
the other hand, as $P_n(x_{ni})=0$, we have $L r_n\ge
\llvert  P_n(x_{ni})-P(x_{ni})\rrvert  =\llvert  P(x_{ni})\rrvert  $. This is a contradiction.

Also, by our construction of $\delta$, it follows that all the
$B(x_j,\delta)$'s are disjoint.

Suppose the lemma is not true, then as the sum of $\ell_j$'s is fixed,
there is at least one $j$ such that, there are $\ell_0$ $x_{ni}$'s in
$B(x_j,r_n^{1/\ell_j})$, with $0\le\ell_0<\ell_j$. WLOG, we can
assume $j=1$ and denote these $\ell_0$ $x_{ni}$'s by $x_{n1}^1,\ldots,x_{n\ell_0}^1$. By Lemma~\ref{P1}, we can choose $x^* \in
B(x_1,r_n^{1/\ell_1})$ such that $\min_{i\in\{1,\ldots,\ell_0\}
}\llvert  x^*-x_{ni}^1\rrvert  \ge dr_n^{1/\ell_1}$ for some $d>0$. By the
construction of $\delta$, we have $\llvert  x^*-x\rrvert  >\delta$ for all $x\in
B(x_j,r_n^{1/\ell_j})$, $j=2,\ldots,m$. Therefore, we have
$\llvert  P(x^*)\rrvert  =\prod_{j=1}^{m}\llvert  x^*-x_j\rrvert  ^{\ell_j}=\llvert  x^*-x_1\rrvert  ^{\ell_1}\prod_{j=2}^{m}\llvert  x^*-x_j\rrvert  ^{\ell_j}=O(r_n)$. On the other hand, we have
$\llvert  P_n(x^*)\rrvert  =\prod_{j=1}^{k}\llvert  x^*-x_{nj}\rrvert  =\prod_{i=1}^{\ell
_0}\llvert  x^*-x_{ni}^1\rrvert  \prod_{x_{nj}\notin B(x_1,r_n^{1/\ell
_1})}\llvert  x^*-x_{nj}\rrvert  >\break \delta^{k-\ell_0}r_n^{\ell_0/\ell_1}$,
contradicting $\llvert  P(x^*)-P(x_n^*)\rrvert  =O(r_n)$. Therefore, the lemma is
proved.

For Lemma~\ref{P3}, write $P_n(x)=\prod_{j=1}^{k}(x-x_{nj})$,
$Q_n(y)=\prod_{j=1}^{k}(y-y_{nj})$ and $P(x)=\prod_{j=1}^{m}(x-x_{j})^{\ell_j}$.
Let\vspace*{2pt} $\delta=\frac{1}{3} \min_{i,j\in\{1,\ldots,m\}, i\ne
j}\llvert  x_i-x_j\rrvert  >0$. By the definition of $\widetilde r_n$, there exists
some $L>0$ such that $L\widetilde r_n\ge\llvert  P_n(x_{ni})-Q_n(x_{ni})\rrvert  $
for all $x_{ni}$. Let $j\in{\{1,\ldots,m\}}$ be given, and let
$d:= (\frac{L}{\delta^{k-\ell_j}} )^{1/\ell_j}>0$. By
Lemma~\ref{P2}, we have exactly $\ell_j$ $x_{ni}$'s and exactly $\ell
_j$ $y_{ni}$'s in $B(x_j,r_n^{1/\ell_j})$. Let $x_{ni}\in
B(x_j,r_n^{1/\ell_j})$ be fixed. By our construction in the proof of
Lemma~\ref{P2}, if $y_{nl}\notin B(x_j,r_n^{1/\ell_j})$, one has
$d(x_{ni},y_{nl})>\delta$. Therefore, for the lemma to be true, we
only need to look at those $y_{nl}\in B(x_j,r_n^{1/\ell_j})$ and show
that at least one such $y_{nl}$ satisfies the desired distance. Suppose
not, that is, for this $x_{ni}\in B(x_j,r_n^{1/\ell_j})$, for any
$y_{nl}\in B(x_j,r_n^{1/\ell_j})$, one has
$d(x_{ni},y_{nl})>\widetilde r_n^{1/\ell_j}$. Note that when
$y_{nl}\notin B(x_j,r_n^{1/\ell_j})$, we have $d(x_{ni},y_{nl})>\delta
$. Hence, we have
$\llvert  Q_n(x_{ni})\rrvert  =\prod_{l=1}^{k}\llvert  x_{ni}-y_{nl}\rrvert  >\delta^{k-\ell
_j}(d\widetilde r_n^{1/\ell_j})^{\ell_j}=L\widetilde r_n$.
However, we also have $L\widetilde r_n\ge
\llvert  Q_n(x_{ni})-P_n(x_{ni})\rrvert  =\llvert  Q_n(x_{ni})\rrvert  $, which is a contradiction.

\subsection{Proof of Lemma \texorpdfstring{\protect\ref{rs}}{3.4}}\label{sec8.2}
Let $\bgma_{l}^*\bbA_k^{-s}={\mathbf b}=(b_1,\ldots,b_n)$. Noting $\llvert  \ep
_{it}\rrvert  <C$,
we have
\begin{eqnarray*}
\hspace*{-2pt}&&\rE\bigl\llvert \bgma_l^*\bbA_k^{-s}
\bgma_k\bigr\rrvert ^{2r}
\\
\hspace*{-2pt}&&\qquad =\frac{1}{2^rT^r}\rE\Biggl(
\Biggl\llvert \sum_{i=1}^n
\ep_{ki}b_i\Biggr\rrvert ^{2r}\Biggr)
\\
\hspace*{-2pt}&&\qquad =\frac{1}{2^rT^r}\rE\mathop{\sum_{i_1+\cdots+i_n=r}}_{j_1+\cdots
+j_n=r}
\frac{(r!)^2}{i_1!j_1!\cdots i_n!j_n!}(\ep_{k1}b_1)^{i_1}(\bar
\ep_{k1}\bar b_1)^{j_1}\cdots(
\ep_{kn}b_n)^{i_n}(\bar\ep_{kn}\bar
b_n)^{j_n}
\\
\hspace*{-2pt}&&\qquad =\frac{1}{2^rT^r}\rE\mathop{\mathop{\sum_{i_1+\cdots+i_n=r}}_{j_1+\cdots+j_n=r}}_{i_1+j_1\ne1
}
\frac{(r!)^2}{i_1!j_1!\cdots i_n!j_n!}(\ep_{k1}b_1)^{i_1}(\bar\ep
_{k1}\bar b_1)^{j_1}\cdots(\ep_{kn}b_n)^{i_n}(
\bar\ep_{kn}\bar b_n)^{j_n}.
\end{eqnarray*}
Let $l$ denote the number $k\le n$ such that $i_k+j_k\ge2$. By the
fact that $\frac{(r!)^2}{(2r)!}\le\frac{r}{2r}\frac
{r-1}{2r-1}\cdots\frac{1}{r+1}\le\frac{1}{2^r}$, we have
\begin{eqnarray*}
&&\rE\bigl\llvert \bgma_l^*\bbA_k^{-s}
\bgma_k\bigr\rrvert ^{2r}
\\
&&\qquad \leq\frac{1}{2^{2r}T^r}\sum_{l=1}^r
\sum_{1\leq j_1< \cdots<
j_l\leq n}\mathop{\sum_{i_1+\cdots+ i_l=2r}}_{i_1\ge2,\ldots,i_l\ge2}
\frac{(2r)!}{i_1 !\cdots i_l !l!} \rE\bigl\llvert \ep _{kj_1}^{i_1}b_{j_1}^{i_1}
\cdots\ep_{kj_l}^{i_l}b_{j_l}^{i_l}\bigr
\rrvert
\\
&&\qquad \leq\frac{1}{2^{2r}T^r}\rE\sum_{l=1}^{r}C^{2r}
\sum_{1\leq j_1<
\cdots<j_l\leq n} \mathop{\sum_{i_1+\cdots+ i_l=2r}}_{i_1\geq2,\ldots,i_l\geq2}
\frac{(2r)!}{i_1 !\cdots i_l !l!}\llvert b_{j_1}\rrvert ^{i_1}\cdots\llvert
b_{j_l}\rrvert ^{i_l}
\\
&&\qquad \le\frac{K_r}{T^r}\sum_{l=1}^r
\sum_{i_1+\cdots+i_l=2r}\rE\prod_{t=1}^l
\Biggl(\sum_{j=1}^n\llvert
b_j\rrvert ^{i_t} \Biggr)
\\
&&\qquad \le\frac{K_r}{T^r}\rE \Biggl(\sum_{j=1}^n
\bigl\llvert b_j^2\bigr\rrvert \Biggr)^r
\\
&&\qquad \le \frac{K_r}{T^r}\rE \bigl(\bgma_l^*\bbA^{-s}_{k}
\bigl(\bbA _k^*\bigr)^{-s}\bgma_l
\bigr)^r.
\end{eqnarray*}
Note that $\llVert  \bgma_l\rrVert  \le K$ and $\llVert  \bbA_k^{-1}\rrVert  \le v_n^{-1}$, we
finally obtain that
\begin{eqnarray*}
\rE\bigl\llvert \bgma_l^*\bbA_k^{-s}
\bgma_k\bigr\rrvert ^{2r}\le\frac{K}{T^rv_n^{2rs}}
\end{eqnarray*}
for some $K>0$. The proof of the lemma is complete.

\subsection{Proof of Lemma \texorpdfstring{\protect\ref{lkj}}{3.5}}\label{sec8.3}
Recall that $a'=a-\underline{\ep}$ and $b'=b+\underline{\ep}$, as
defined at the end of Section~\ref{sec3}. Therefore, we have
\begin{eqnarray*}
&&\rP \biggl(\frac{1}{2T}\sum\frac{1}{\llvert  \lambda_{kj}-z\rrvert  ^2}>K \biggr)
\\
&&\qquad \le\rP \biggl(\sum_{\lambda_{kj}\notin[a', b']}\frac{1}{\llvert  \lambda
_{kj}-u\rrvert  ^2+v_n^2}>TK
\biggr)
\\
&&\quad\qquad{}+\rP \biggl(\sum_{\lambda_{kj}\in[a',
b']}\frac{1}{\llvert  \lambda_{kj}-u\rrvert  ^2+v_n^2}>TK
\biggr)
\\
&&\qquad \le\rP \bigl(n\underline{\ep}^{-2}>TK \bigr)+\rP
\bigl(nv_n^{-2}F_{nk}\bigl(\bigl[a',
b'\bigr]\bigr)>TK \bigr)
\\
&&\qquad \le 0+ \rP\biggl(\llVert F_n-F_{c_{n}}\rrVert \ge
\frac{K}{2c}n^{-1/53}\biggr)=o\bigl(n^{-t}\bigr).
\end{eqnarray*}
Here, we pick $K>c\underline{\ep}^{-2}$ so that the first probability
is 0. The second probability follows (\ref{Fn}). The proof is complete.

\subsection{Proof of Lemma \texorpdfstring{\protect\ref{x1}}{3.6}, part \textup{(a)}}\label{sec8.4}
For (i)(a),
by definition of $x_{nj}$, $j=0,1$, we have
\begin{eqnarray*}
x_{n0,1}&=&\tfrac{1}2 \Bigl(1\pm\sqrt{1-4a_{n}^2}
\Bigr):=\tfrac{1}2\bigl(1\pm (\tilde\alpha+i\tilde\beta)
\bigr).
\end{eqnarray*}
Therefore,
\begin{eqnarray*}
\biggl\llvert \frac{x_{n0}}{x_{n1}}\biggr\rrvert &=&\cases{ \displaystyle\sqrt{
\frac{(1-\tilde\alpha)^2+\tilde\beta^2}{(1+\tilde\alpha
)^2+\tilde\beta^2}} <1-\frac{2\tilde\alpha}{(1+\tilde\alpha)^2+\tilde\beta^2}, &\quad if $\tilde\alpha>0$,
\vspace*{5pt}\cr
\displaystyle \sqrt{\frac{(1+\tilde\alpha)^2+\tilde\beta^2}{(1-\tilde\alpha
)^2+\tilde\beta^2}} <1-\frac{2\llvert  \tilde\alpha\rrvert  }{(1-\tilde\alpha)^2+\tilde\beta^2}, &\quad if $\tilde\alpha<0$}
\\
&=&1-\frac{\llvert  \tilde\alpha\rrvert  }{2\llvert  x_{n1}^2\rrvert  }<1-\eta_1 v_n^2
\llvert \tilde \alpha\rrvert ,
\end{eqnarray*}
where the last inequality follows from the fact that
$x_{n1}^2=x_{n1}-a_n^2=O(v_n^{-2})$.

Thus, to complete the proof of (i)(a), it suffices to show that there
is a constant $\eta_2>0$ such that $\llvert  \tilde\alpha\rrvert  >\eta_2v_n$.

Write $c_n\rE m_n(z)=2a_n=\alpha+i\beta$ where $\alpha$ and $\beta$
are real.
Then, by the formula of square root of complex numbers [see (2.3.2) of
\citet{BaiSil10}] we have
\begin{eqnarray*}
\sqrt{1-4a_n^2}=\tilde\alpha+i\tilde\beta,
\end{eqnarray*}
where
\begin{eqnarray*}
\tilde\alpha&=&\frac{-\sqrt{2}\alpha\beta}{\sqrt{\sqrt
{(1-\alpha^2+\beta^2)^2+4\alpha^2\beta^2}-(1-\alpha^2+\beta^2)}}.
\end{eqnarray*}
Obviously, when $1-\alpha^2+\beta^2>0$, by $\sqrt{(1-\alpha^2+\beta
^2)^2+4\alpha^2\beta^2}-(1-\alpha^2+\beta^2)<2\llvert  \alpha\rrvert  \beta$ we have
\[
\llvert \tilde\alpha\rrvert >1/ \sqrt{\llvert \alpha\rrvert \beta}>1/\bigl
\llvert c_n\rE m_n(z)\bigr\rrvert >\eta_2v_n,
\]
for all large $n$ such that $c_n\eta_2<1$, where $\eta_2\in(0,c^{-1})$.

On the other hand, if $1-\alpha^2+\beta^2<0$, by $\alpha^2>1+\beta
^2$ we have
\begin{eqnarray*}
\llvert \tilde\alpha\rrvert &>&\frac{\llvert  \alpha\rrvert  \beta}{\sqrt[4]{(1-\alpha^2+\beta
^2)^2+4\alpha^2\beta^2}}=\frac{\llvert  \alpha\rrvert  \beta}{\sqrt[4]{(1-\alpha
^2-\beta^2)^2+4\beta^2}}>\beta/
\sqrt{2}.
\end{eqnarray*}
Then the assertion that $\llvert  \tilde\alpha\rrvert  >\eta_2 v_n$ is proved if one
can show that $\beta>\eta_3v_n$ for some $\eta_3>0$. This is trivial
if one notices
\begin{eqnarray*}
\beta=v\int\frac{1}{(x-u)^2+v^2}\,d\rE F_n(x)>v_n
\bigl(4A^2+1\bigr)^{-1}\rE F_n\bigl([-A,A]
\bigr),
\end{eqnarray*}
when $\llvert  z\rrvert  <A$ and $v\in(v_n,1)$. The conclusion (i) is proved.

For (ii)(a), by $x_{n1}+x_{n0}=1$ and $\llvert  x_{n1}\rrvert  >\llvert  x_{n0}\rrvert  $, we conclude
that $\llvert  x_{n1}\rrvert  \ge\frac{1}2$. Since $x_{n1}=\frac{1}2(1\pm\sqrt
{1-4a_n^2})$, we conclude that
\[
\llvert x_{n1}\rrvert \le\frac{1}2\Bigl(1+\Bigl\llvert \sqrt
{1-4a_n^2}\Bigr\rrvert \Bigr)\le
Kv_n^{-1}. %
\]
For (iii)(a), by noting that
\begin{eqnarray*}
\llvert x_{n1}-x_{n0}\rrvert ^2&=&\bigl(1-
\alpha^2+\beta^2\bigr)^2+4\alpha^2
\beta ^2=\bigl(1-\alpha^2-\beta^2
\bigr)^2+4\beta^2.
\end{eqnarray*}
Then the conclusion (iii)(a) follows from the fact $\llvert  \beta\rrvert  >\eta
_3v_n$ that is shown in the proof of part (i)(a) of the lemma.

The conclusion (iv)(a) follows from
\[
\frac{\llvert  x_{n0}\rrvert  }{\llvert  x_{n1}-x_{n0}\rrvert  }\le\frac{1}2 \biggl(\frac{1}{\llvert  \sqrt
{1-4a_n^2}\rrvert  }+1 \biggr)\le
Kv_n^{-1}, %
\]
where the last inequality follows from conclusion (iii)(a).

The proof of the lemma is complete.

\subsection{Proof of Lemma \texorpdfstring{\protect\ref{Ek+tau}}{3.7}(a)}\label{sec8.5}
Recall that $a_n=\frac{c_n\rE m_n}2$. Write $W_{k}=\bgma_{k+\tau
}^*\*\bbA_k^{-1}\bgma_{k+\tau}$ and $W_{k,k+\tau,\ldots,k+s\tau
}=\bgma_{k+(s+1)\tau}^*\bbA_{k,k+\tau,\ldots,k+s\tau}^{-1}\bgma
_{k+(s+1)\tau}$.
Denote\break $\widetilde\bbA_{k,\ldots,k+(s-1)\tau}=\bbA_{k,\ldots,k+s\tau}+\bgma_{k+(s+1)\tau}\bgma_{k+s\tau}^*$. Apply the identity
\begin{eqnarray*}
\bigl(\bbB+\bga\bgma^*\bigr)^{-1}=\bbB^{-1}-
\frac{\bbB^{-1}\bga\bgma
^*\bbB^{-1}}{1+\bgma^*\bbB^{-1}\bga},
\end{eqnarray*}
we have
\begin{eqnarray*}
\bbA_{k,\ldots,k+(s-1)\tau}^{-1}&=&\bigl(\widetilde\bbA_{k,\ldots,k+(s-1)\tau}+
\bgma_{k+s\tau}\bgma_{k+(s+1)\tau}^*\bigr)^{-1}
\\
&=&\widetilde\bbA_{k,\ldots,k+(s-1)\tau}^{-1}-\frac{\widetilde
\bbA_{k,\ldots,k+(s-1)\tau}^{-1}\bgma_{k+(s+1)\tau}\bgma_{k+s\tau
}^*\widetilde\bbA_{k,\ldots,k+(s-1)\tau}^{-1}}{1+\bgma_{k+s\tau
}^*\widetilde\bbA_{k,\ldots,k+(s-1)\tau}^{-1}\bgma_{k+(s+1)\tau
}},
\\
\widetilde\bbA_{k,\ldots,k+(s-1)\tau}&=&\bigl(\bbA_{k,\ldots,k+s\tau
}+
\bgma_{k+(s+1)\tau}\bgma_{k+s\tau}^*\bigr)^{-1}
\\
&=&\bbA_{k,\ldots,k+s\tau}^{-1}-\frac{\bbA_{k,\ldots,k+s\tau
}^{-1}\bgma_{k+(s+1)\tau}\bgma_{k+s\tau}^*\bbA_{k,\ldots,k+s\tau
}^{-1}}{1+\bgma_{k+s\tau}^*\bbA_{k,\ldots,k+s\tau}^{-1}\bgma
_{k+(s+1)\tau}}.
\end{eqnarray*}
Therefore, we have
\begin{eqnarray*}
&&\bgma_{k+s\tau}^*\bbA_{k,\ldots,k+(s-1)\tau}^{-1}
\\
&&\qquad =\bgma_{k+s\tau}^*\widetilde\bbA_{k,\ldots,k+(s-1)\tau
}^{-1}
\\
&&\quad\qquad{}-
\frac{\bgma_{k+s\tau}^*\widetilde\bbA_{k,\ldots,k+(s-1)\tau}^{-1}\bgma_{k+(s+1)\tau}\bgma_{k+s\tau}^*\widetilde
\bbA_{k,\ldots,k+(s-1)\tau}^{-1}}{1+\bgma_{k+s\tau}^*\widetilde
\bbA_{k,\ldots,k+(s-1)\tau}^{-1}\bgma_{k+(s+1)\tau}}
\\
&&\qquad =\frac{\bgma_{k+s\tau}^*\widetilde\bbA_{k,\ldots,k+(s-1)\tau
}^{-1}}{1+\bgma_{k+s\tau}^*\widetilde\bbA_{k,\ldots,k+(s-1)\tau
}^{-1}\bgma_{k+(s+1)\tau}}
\end{eqnarray*}
and
%
\begin{eqnarray}\label{eqwwk}
&&\bgma_{k+s\tau}^*\bbA_{k,\ldots,k+(s-1)\tau}^{-1}
\bgma_{k+s\tau}\nonumber
\\
&&\qquad  = \frac{\bgma_{k+s\tau}^*\widetilde\bbA_{k,\ldots,k+(s-1)\tau
}^{-1}\bgma_{k+s\tau}}{1+\bgma_{k+(s+1)\tau}^*\widetilde\bbA
_{k,\ldots,k+(s-1)\tau}^{-1}\bgma_{k+s\tau}}
\nonumber
\\
&&\qquad =
\biggl(\bgma_{k+s\tau}^*\bbA_{k,\ldots,k+s\tau}^{-1}\bgma
_{k+s\tau}\nonumber
\\
&&\hspace*{37pt}{}-\frac{\bgma^*_{k+s\tau}\bbA_{k,\ldots,k+s\tau
}^{-1}\bgma_{k+(s+1)\tau}
\bgma_{k+s\tau}^*\bbA_{k,\ldots,k+s\tau}
^{-1}\bgma_{k+s\tau}}{1+\bgma_{k+s\tau}^*\bbA_{k,\ldots,k+s\tau
}^{-1}\bgma_{k+(s+1)\tau}}\biggr)
\\
&&\quad\qquad{} \Big/
\biggl(
1+\bgma_{k+(s+1)\tau}^*\bbA_{k,\ldots,k+s\tau}^{-1}\bgma
_{k+s\tau}\nonumber
\\
&&\hspace*{47pt}{}-\frac{\bgma^*_{k+(s+1)\tau}\bbA_{k,\ldots,k+s\tau}^{-1}
\bgma_{k+(s+1)\tau}\bgma_{k+s\tau}^*\bbA_{k,\ldots,k+s\tau}
^{-1}\bgma_{k+s\tau}}{1+\bgma_{k+s\tau}^*\bbA_{k,\ldots,k+s\tau
}^{-1}\bgma_{k+(s+1)\tau}}\biggr)
\nonumber
\\
&&\qquad =\frac{\sklfrac{c_n}2\rE m_n(z)+r_1(k+s\tau)}{1-\sklfrac{c_n}2\rE
m_n(z)\bgma_{k+(s+1)\tau}^*
\bbA_{k,\ldots,k+s\tau}^{-1}\bgma_{k+(s+1)\tau}+r_2(k+s\tau
)},\nonumber
\end{eqnarray}
that is,
%
\begin{eqnarray}
W_{k,\ldots,k+(s-1)\tau}=\frac{a_n+r_1(k+s\tau)}{1-a_nW_{k,\ldots,k+s\tau}+r_2(k+s\tau)}, \label{eqwwk1}
\end{eqnarray}
where
\begin{eqnarray*}
r_1(k+s\tau) &=&\bgma_{k+s\tau}^*\bbA_{k,\ldots,k+s\tau}^{-1}
\bgma_{k+s\tau
}-a_n,
\\
r_2(k+s\tau)&=&-\bigl(
\bgma_{k+s\tau}^*\bbA_{k,\ldots,k+s\tau
}^{-1}\bgma_{k+s\tau}-a_n
\bigr)\bgma_{k+(s+1)\tau}^*\bbA_{k,\ldots,k+s\tau}^{-1}
\bgma_{k+(s+1)\tau}
\\
&&{}+\bgma_{k+(s+1)\tau}^*\bbA_{k,\ldots,k+s\tau}^{-1}\bgma
_{k+s\tau}+\bgma_{k+s\tau}^*\bbA_{k,\ldots,k+s\tau}^{-1}\bgma
_{k+(s+1)\tau}
\\
&&{}+\bgma_{k+(s+1)\tau}^*\bbA_{k,\ldots,k+s\tau}^{-1}\bgma
_{k+s\tau}\bgma_{k+s\tau}^*\bbA_{k,\ldots,k+s\tau}^{-1}\bgma
_{k+(s+1)\tau}.
\end{eqnarray*}
When $k\le T-v_n^{-4}$, applying this relation $\ell$ times ($\ell=[
v^{-4}_n]$), we may express $W_k$ in the following form:
\begin{eqnarray*}
W_k=\frac{(a_n+r_1(k+\tau))(\ga_{k+\tau,\ell}-a_n\gamma_{k+\tau,\ell}W_{k,k+\tau,\ldots,k+(\ell+1)\tau})}{\ga_{k,\ell
}-a_n\gamma_{k,\ell}
W_{k,k+\tau,\ldots,k+(\ell+1)\tau}},
\end{eqnarray*}
where the coefficients satisfy the recursive relation
%
\begin{eqnarray}
\label{ga} \ga_{k+s\tau,\ell}&=&\bigl(1+r_2(k+s\tau)\bigr)
\ga_{k+(s+1)\tau,\ell
}\nonumber
\\
&&{} -a_n\bigl(a_n+r_1(k+s\tau)
\bigr)\ga_{k+(s+2)\tau,\ell},
\nonumber
\\
\ga_{k+\ell\tau,\ell}&=&1+r_2(k+\ell\tau), \qquad \ga_{k+(\ell+1)\tau,\ell}=1,
\nonumber\\[-8pt]\\[-8pt]\nonumber
\gamma_{k+s\tau,\ell}&=&\bigl(1+r_2(k+s\tau)\bigr)
\gamma_{k+(s+1)\tau,\ell
}\nonumber
\\
&&{}-a_n\bigl(a_n+r_1(k+s
\tau)\bigr)\gamma_{k+(s+2)\tau,\ell},
\nonumber
\\
\gamma_{k+\ell\tau,\ell}&=&1,\qquad \gamma_{k+(\ell+1)\tau,\ell
}=0.
\nonumber
\end{eqnarray}
Notice that $v_n=n^{-1/52}$. Employing Lemma~\ref{B26} and an
estimation similar to (\ref{W1}), for any fixed $t$, one has
%
\begin{eqnarray}
\rP\bigl(\bigl\llvert r_i(k+\ell\tau)\bigr\rrvert \ge
v_n^{12}\bigr)=o\bigl(n^{-t}\bigr)\qquad\mbox{for
}i=1,2. \label{eqw6}
\end{eqnarray}
As in the proof of Lemma B.3 of \citet{Jinetal14}, by letting $\ell
=[v_n^{-4}]$, it follows by induction that
%
\begin{eqnarray}
\ga_{k+l\tau,\ell}&=&(1-\ga)\prod_{\mu=1}^{\ell-l+1}
\nu_{\mu,1}+\ga\prod_{\mu=1}^{\ell-l+1}
\nu_{\mu,0}, \label{eqww2}
\end{eqnarray}
where $\nu_{1,i}$, $i=1,0$ (with $\llvert  \nu_{1,1}\rrvert  >\llvert  \nu_{1,0}\rrvert  $) are
defined by the two roots of the quadratic equation
\[
x^2=\bigl(1+r_2(k+\ell\tau)\bigr)x-a_n
\bigl(a_n+r_1(k+\ell\tau)\bigr) %
\]
and $\ga$ is such that
\[
(1-\ga)\nu_{1,1}+\ga\nu_{1,0}=1+r_2(k+\ell\tau)=
\ga_{k+\ell\tau,\ell}. %
\]
Recall that $x_{ni}$, $i=1,0$ (with $\llvert  x_{n1}\rrvert  >\llvert  x_{n0}\rrvert  $) are two roots
of the quadratic equation
\[
x^2=x-a_n^2. %
\]
Applying Lemmas~\ref{P1}--\ref{P3} to the above two quadratic
equations and using (\ref{eqw6}), we have
%
\begin{eqnarray}
&&\rP\bigl(\llvert \nu_{1,i}-x_{ni}\rrvert
\ge2v_n^6\bigr)
\nonumber\\[-8pt]\label{eqcc1} \\[-8pt]\nonumber
&&\qquad \le \rP\bigl(\bigl\llvert r_1(k+\ell\tau)\bigr\rrvert \ge
v_n^{12}\bigr)+\rP\bigl(\bigl\llvert r_2(k+\ell
\tau )\bigr\rrvert \ge v_n^{12}\bigr)=o\bigl(n^{-t}
\bigr),
\\
&&\rP\biggl(\biggl\llvert \ga-\frac{x_{n0}}{x_{n0}-x_{n1}}\biggr\rrvert
\ge3v_n^6\biggr)
\nonumber
\\
&&\qquad \le\rP\bigl(\llvert \nu_{1,0}-x_{n0}\rrvert \ge
v_n^6\bigr)+\rP\bigl(\llvert \nu_{1,1}-x_{n1}
\rrvert \ge v_n^6\bigr)+\rP\bigl(\bigl\llvert
r_2(k+\ell\tau)\bigr\rrvert \ge v_n^6
\bigr)\hspace*{-20pt}\label{eqcc2}
\\
&&\qquad =o\bigl(n^{-t}\bigr).\nonumber
\end{eqnarray}
By induction, one has for $\mu\in[1,\ell]$
\begin{eqnarray*}
\nu_{\mu+1,i}&=&1+r_2\bigl(k+(\ell-\mu)\tau\bigr)-
\frac{a_n(a_n+r_1(k+(\ell
-\mu)\tau))}{\nu_{\mu,i}}
\end{eqnarray*}
and can similarly verify that
\[
\rP\bigl(\llvert \nu_{\mu,i}-x_{ni}\rrvert \ge2\mu
v_n^{6}\bigr)\le\sum_{l=1}^{\mu}
\sum_{j=1}^2\rP\bigl(\bigl\llvert
r_j(k+l\tau)\bigr\rrvert \ge v_n^{12}\bigr)=o
\bigl(n^{-t}\bigr). %
\]
Therefore, we have
\begin{eqnarray*}
\rP \bigl(\bigl\llvert \ga_{k+\tau,\ell}- \bigl((1-\ga)x_{n1}^{\ell}+
\ga x_{n0}^{\ell} \bigr)\bigr\rrvert \ge v_n^6
\bigr)&\le&\sum_{\mu=1}^{\ell}\sum
_{i=0}^1\rP\bigl(\llvert \nu_{\mu,i}-x_{ni}
\rrvert \ge2\mu v_n^{6}\bigr)
\\
&=&o\bigl(n^{-t}
\bigr),
\\
\rP \bigl(\bigl\llvert \ga_{k,\ell}- \bigl((1-\ga)x_{n1}^{\ell+1}+
\ga x_{n0}^{\ell+1} \bigr)\bigr\rrvert \ge v_n^6
\bigr)&\le&\sum_{\mu=1}^{\ell
+1}\sum
_{i=0}^1\rP\bigl(\llvert \nu_{\mu,i}-x_{ni}
\rrvert \ge2\mu v_n^{6}\bigr)
\\
&=& o\bigl(n^{-t}
\bigr),
\end{eqnarray*}
and
\begin{eqnarray*}
&&\rP\biggl(\biggl\llvert \frac{\ga_{k+\tau,\ell}}{\ga_{k,\ell}}-\frac{1}{x_{n1}}\biggr\rrvert
\ge v_n^6\biggr)
\\
&&\qquad \le\rP \bigl(\bigl\llvert \ga_{k+\tau,\ell}- \bigl((1-\ga)x_{n1}^{\ell}+
\ga x_{n0}^{\ell} \bigr)\bigr\rrvert \ge v_n^6
\bigr)
\\
&&\quad\qquad{} +\rP \bigl(\bigl\llvert \ga_{k,\ell}- \bigl((1-
\ga)x_{n1}^{\ell+1}+\ga x_{n0}^{\ell+1} \bigr)
\bigr\rrvert \ge v_n^6 \bigr)
\\
&&\quad\qquad{} +\rP\bigl(\llvert \nu_{\ell+1,1}-x_{n1}\rrvert \ge2(\ell+1)
v_n^{6}\bigr)
\\
&&\qquad =o\bigl(n^{-t}\bigr).
\end{eqnarray*}
Similarly, we have
\begin{eqnarray*}
\gamma_{k+l\tau,\ell}&=&(1-\tilde\ga)\prod_{\mu=1}^{\ell
-l+1}
\tilde\nu_{\mu,1}+\tilde \ga\prod_{\mu=1}^{\ell
-l+1}
\tilde\nu_{\mu,0},
\end{eqnarray*}
where $\tilde\nu_{\mu,i}$, $i=1,0$, are the two roots of the
quadratic equation
\[
x^2=\bigl(1+r_2\bigl(k+(\ell-1)\tau\bigr)
\bigr)x-a_n\bigl(a_n+r_1\bigl(k+(\ell-1)\tau
\bigr)\bigr), %
\]
and $\tilde\ga$ satisfies
\[
(1-\tilde\ga)\tilde\nu_{1,1}+\tilde\ga\tilde\nu _{1,0}=1+r_2
\bigl(k+(\ell-1)\tau\bigr)=\gamma_{k+(\ell-1)\tau,\ell}. %
\]
One can similarly prove that $\tilde\nu_{\mu,i}$, $i=0,1$, satisfy
\[
\rP\bigl(\llvert \tilde \nu_{\mu,i}-x_{ni}\rrvert \ge2\mu
v_n^{6}\bigr)\le\sum_{l=0}^{\mu}
\sum_{j=1}^2\rP\bigl(\bigl\llvert
r_j(k+l\tau)\bigr\rrvert \ge v_n^{12}\bigr)=o
\bigl(n^{-t}\bigr), %
\]
and
\[
\rP\biggl(\biggl\llvert \tilde\ga-\frac{x_{n0}}{x_{n0}-x_{n1}}\biggr\rrvert
\ge3v_n^6\biggr)=o\bigl(n^{-t}\bigr).
\]
Therefore, we have
\begin{eqnarray*}
\rP \bigl(\bigl\llvert \gamma_{k+\tau,\ell}- \bigl((1-\tilde
\ga)x_{n1}^{\ell
}+\tilde\ga x_{n0}^{\ell}
\bigr)\bigr\rrvert \ge v_n^6 \bigr)&\le&\sum
_{\mu
=1}^{\ell}\sum_{i=0}^1
\rP\bigl(\llvert \tilde\nu_{\mu,i}-x_{ni}\rrvert \ge2\mu
v_n^{6}\bigr)
\\
&=&o\bigl(n^{-t}\bigr),
\\
\rP \bigl(\bigl\llvert \gamma_{k,\ell}- \bigl((1-\tilde
\ga)x_{n1}^{\ell
+1}+\tilde\ga x_{n0}^{\ell+1}
\bigr)\bigr\rrvert \ge v_n^6 \bigr)&\le&\sum
_{\mu
=1}^{\ell+1}\sum_{i=0}^1
\rP\bigl(\llvert \tilde\nu_{\mu,i}-x_{ni}\rrvert \ge2\mu
v_n^{6}\bigr)
\\
&=&o\bigl(n^{-t}\bigr),
\end{eqnarray*}
and
\begin{eqnarray*}
&&\rP\biggl(\biggl\llvert \frac{\gamma_{k+\tau,\ell}}{\gamma_{k,\ell}}-\frac{1}{x_{n1}}\biggr\rrvert
\ge v_n^6\biggr)
\\
&&\qquad \le\rP \bigl(\bigl\llvert \gamma_{k+\tau,\ell}- \bigl((1-\tilde\ga
)x_{n1}^{\ell}+\tilde\ga x_{n0}^{\ell}
\bigr)\bigr\rrvert \ge v_n^6 \bigr)
\\
&&\quad\qquad{} +\rP \bigl(\bigl
\llvert \gamma_{k,\ell}- \bigl((1-\tilde\ga)x_{n1}^{\ell
+1}+
\tilde\ga x_{n0}^{\ell+1} \bigr)\bigr\rrvert \ge
v_n^6 \bigr)
\\
&&\quad\qquad{} +\rP\bigl(\llvert \tilde\nu_{\ell+1,1}-x_{n1}\rrvert \ge2(
\ell+1) v_n^{6}\bigr)
\\
&&\qquad =o\bigl(n^{-t}\bigr).
\end{eqnarray*}
Substituting back to the recursive expression of $W_k$, we thus have
%
\begin{eqnarray}
\rP\biggl(\biggl\llvert W_k-\frac{a_n}{x_{n1}}\biggr\rrvert \ge
v_n^{6}\biggr)=o\bigl(n^{-t}\bigr).
\label{eqwwk7}
\end{eqnarray}
The proof of this lemma is complete.
\subsection{Proof of Lemma \texorpdfstring{\protect\ref{cross}}{3.8}(a)}\label{sec8.6}
When $\tau<k\le2\tau$, the lemma is obviously true because $\bgma
_{k-\tau}$ is independent of
$\bbA_k$. Similarly, the lemma is true when $T-\tau< k\le T$.

When $2\tau<k\le T/2$, similar to (\ref{eqwwk}), we have
\begin{eqnarray*}
&& \widetilde W_{k,\ldots,k+s\tau}
\\
&&\qquad := \bgma_{k-\tau}^*\bbA_{k,k+\tau,\ldots,k+(s-1)\tau}^{-1}
\bgma_{k+s\tau}
\\
&&\qquad = \frac{\bgma_{k-\tau}^*(\bbA_{k,k+\tau,\ldots,k+s\tau}+\bgma
_{k+(s+1)\tau}\bgma_{k+s\tau}^*)^{-1}\bgma_{k+s\tau}}{1+\bgma
_{k+(s+1)\tau}^*(\bbA_{k,k+\tau,\ldots,k+s\tau}
+\bgma_{k+(s+1)\tau}\bgma_{k+s\tau}^*)^{-1}\bgma_{k+s\tau}}
\\
&&\qquad = \biggl(\bgma_{k-\tau}^*\bbA_{k,k+\tau,\ldots,k+s\tau
}^{-1}\bgma_{k+s\tau}
\\
&&\hspace*{37pt}{}-\frac{\bgma_{k-\tau}^*\bbA_{k,k+\tau,\ldots,k+s\tau}^{-1}\bgma_{k+(s+1)\tau}\bgma_{k+s\tau}^*\bbA
_{k,k+\tau,\ldots,k+s\tau}^{-1}\bgma_{k+s\tau}}{1+\bgma_{k+s\tau
}^*\bbA_{k,k+\tau,\ldots,k+s\tau}^{-1}\bgma_{k+(s+1)\tau}}\biggr)
\\
&&\quad\qquad{}\Big/
\biggl(
1+\bgma_{k+(s+1)\tau}^*\bbA_{k,k+\tau,\ldots,k+s\tau}^{-1}\bgma
_{k+s\tau}
\\
&&\hspace*{47pt}{}-\frac{\bgma^*_{k+(s+1)\tau}\bbA_{k,k+\tau,\ldots,k+s\tau}^{-1}\bgma_{k+(s+1)\tau}\bgma_{k+s\tau}^*\bbA_{k,k+\tau,\ldots,k+s\tau}
^{-1}\bgma_{k+s\tau}}{1+\bgma_{k+s\tau}^*\bbA_{k,k+\tau,\ldots,k+s\tau}^{-1}\bgma_{k+(s+1)\tau}}\biggr)
\\
&&\qquad = \frac{\widetilde r_1(k+s\tau)-\widetilde W_{k,\ldots,k+(s+1)\tau
}a_n}{1+ r_2(k+s\tau)-a_nW_{k,\ldots,k+s\tau}},
\end{eqnarray*}
where
\begin{eqnarray*}
\widetilde r_1(k+s\tau)&=&\bgma_{k-\tau}^*
\bbA_{k,\ldots,k+s\tau
}^{-1}\bgma_{k+s\tau}\bigl(1+
\bgma_{k+s\tau}^*\bbA_{k,\ldots,k+s\tau
}^{-1}\bgma_{k+(s+1)\tau}
\bigr)
\\
&&{}-\widetilde W_{k,\ldots,k+(s+1)\tau}\bigl(\bgma_{k+s\tau}^*\bbA
_{k,\ldots,k+s\tau}^{-1}\bgma_{k+s\tau}-a_n\bigr).
\end{eqnarray*}
Similarly, one can show that
\begin{eqnarray*}
\rP\bigl(\bigl\llvert \widetilde r_1(k+s\tau)\bigr\rrvert \ge
v_n^{12}\bigr)=o\bigl(n^{-t}\bigr).
\end{eqnarray*}
When $\llvert  \widetilde r_1(t+s\tau)\rrvert  \le v_n^{12}$, $\llvert  r_2(k+s\tau)\rrvert  \le
v_n^{12}$, and $\llvert  W_{k,\ldots,k+s\tau}-\frac{a_n}{x_{n1}}\rrvert  \le v_n^6$,
we have
\begin{eqnarray*}
\llvert \widetilde W_{k,\ldots,k+s\tau}\rrvert &\le& \frac
{v_n^{12}}{\llvert  x_{n1}\rrvert  -v_n^{5}}+\llvert
\widetilde W_{k,\ldots,k+(s+1)\tau
}\rrvert \biggl\llvert \frac{\llvert  a_n\rrvert  }{\llvert  x_{n1}\rrvert  }+v_n^5
\biggr\rrvert
\\
&\le& 3v_n^{12}+\llvert \widetilde W_{k,\ldots,k+(s+1)\tau}
\rrvert \biggl(1-\frac{1}2\eta v_n^3+v_n^4
\biggr),
\end{eqnarray*}
where the second term follows from the fact that
\[
\frac{\llvert  a_n\rrvert  }{\llvert  x_{n1}\rrvert  }=\sqrt{\frac{\llvert  x_{n0}\rrvert  }{\llvert  x_{n1}\rrvert  }}\le1-\frac{1}2\eta
v_n^3. %
\]
Therefore, when $v_{n}^{-4}<\ell<v_n^{-5}$,
\[
\llvert \widetilde W_{k}\rrvert \le3\ell v_{n}^{12}+
\llvert \widetilde W_{k,\ldots,k+\ell
\tau}\rrvert \bigl\llvert 1-\tfrac{1}2\eta
v_n^3+v_n^4\bigr\rrvert
^\ell\le v_n^6.
\]
The lemma then follows by the fact that
\begin{eqnarray*}
&&\rP\bigl(\llvert \widetilde W_k\rrvert \ge v_n^6
\bigr)
\\
&&\qquad \le \sum_{s=1}^\ell \biggl(\rP \bigl(
\bigl\llvert \widetilde r_1(k+s\tau)\bigr\rrvert \ge
v_n^{12}\bigr)+\rP\bigl(\bigl\llvert r_2(k+s
\tau)\bigr\rrvert \ge v_n^{12}\bigr)
\\
&&\hspace*{118pt}{}+\rP\biggl(\biggl\llvert
W_{k,\ldots,k+s\tau}-\frac{a_n}{x_{n1}}\biggr\rrvert \ge v_n^6
\biggr) \biggr)
\\
&&\qquad =o\bigl(n^{-t}\bigr).
\end{eqnarray*}
The proof of the lemma is complete.
%
\subsection{Proof of Lemma \texorpdfstring{\protect\ref{x1}}{3.6}, part (\textup{b})}\label{sec8.7}
Let $ x_{1}$ and $x_{0}$ be the two roots of the quadratic equation
\[
x^2=x-\breve a^2, %
\]
where $\breve a=\breve a(z)=cm(z)/2$ and $m(z)$ satisfies (\ref{m}).
We claim that
%
\begin{eqnarray}
\sup_{u\in[a,b]}\frac{\llvert  x_0(z)\rrvert  }{\llvert  x_1(z)\rrvert  }\le1-\eta \label{eqb9}
\end{eqnarray}
for some $\eta\in(0,1)$. Otherwise, there will be a sequence $\{z_k\}
$ with $\Re(z_k)\in[a,b]$ and
\[
\frac{\llvert  x_0(z_k)\rrvert  }{\llvert  x_1(z_k)\rrvert  }\to1. %
\]
Then we can select a convergent subsequence $\{z_{k'}\}\to z_0$. If
$z_0=\infty$, then $\breve a(z_0)=0$ and hence $x_1=1$ and $x_0=0$. It
contradicts the fact that
\[
\frac{\llvert  x_0(z_0)\rrvert  }{\llvert  x_1(z_0)\rrvert  }=1. %
\]
The only case to make the equality above true is that $\breve a(z_0)$
is real and its absolute value is ${\ge}\frac{1}2$. That is, $z_0$ is
real and $\llvert  \breve a(z_0)\rrvert  \ge\frac{1}2$. Since $\breve a(\infty)=0$,
there is a real number $z'$ between $z_0$ and $\operatorname{sgn}(z_0)\infty$
such that $\llvert  \breve a(z')\rrvert  =\frac{1}2$ which contradicts the equation
(\ref{m}). Therefore, (\ref{eqb9}) is proved.

Since $m_n^0(z)\to m(z)$ uniformly for all $\Re(z)\in[a,b]$, we
conclude that there is a constant $\eta\in(0,1)$ such that
\[
\sup_{\Re(z)\in[a,b]}\frac{\llvert  \widetilde x_{n0}\rrvert  }{\llvert  \widetilde
x_{n1}\rrvert  }<1-\eta, %
\]
where $\widetilde x_{n1}$ and $\widetilde x_{n0}$ are the two roots
of the equation
\[
x^2=x-\tfrac{1}4c_n^2
\bigl(m_n^0(z)\bigr)^2. %
\]
By what has been proved in Section~\ref{sec3}, we have $\sup_{1>\Im
(z)\ge n^{-1/52}}\llvert  \rE m_n(z)-m_n^0(z)\rrvert  \to0$. Thus,
\[
\mathop{\sup_{\Re(z)\in[a,b]}}_{1> \Im(z)\ge n^{-1/52}}\frac
{\llvert  x_{n0}\rrvert  }{\llvert  x_{n1}\rrvert  }\le1-
\eta. %
\]
The conclusion (i)(b) follows.

We then prove the conclusion (v). In the proof of (i)(b), we actually
proved that there is a constant $\eta\in(0,\frac{1}2)$ such that for
all $u\in[a,b]$,
\[
\bigl\llvert \breve a(u)\bigr\rrvert <\tfrac{1}2 -\eta. %
\]
By the uniform continuity of $\breve a(z)$ for all $\Re(z)\in[a,b]$.
we have
\[
\sum_{u\in[a,b], v\in(0,\delta_n)}\bigl\llvert \breve a(u+iv)-\breve a(u)
\bigr\rrvert \to 0\qquad\mbox{as } \delta_n\to0. %
\]
Then conclusion (v) follows from the fact that
$\sup_{1>\Im(z)\ge n^{-1/52}}\llvert  \rE m_n(z)-m_n^0(z)\rrvert  \to0$.

The first conclusion of (ii)(b) is the same as (ii)(a) and the
second follows easily from the fact that $\llvert  a_n(z)\rrvert  \le\frac{1}2$ and the
argument that $\llvert  x_{n1}\rrvert  \le\frac{1}2(1+\sqrt{1+4\llvert  a_n^2\rrvert  })\le\frac{3}2$.

The conclusion (iii)(b) follows from the fact that
$\llvert  x_{n1}-x_{n0}\rrvert  =\llvert  \sqrt{1-4 a^2_n}\rrvert  \ge\sqrt{4\eta(1-\eta)}$.
The conclusion (iv)(b) follows from conclusions (ii)(b) and (iii)(b).
The goal of this section is reached.

\subsection{Proof of Lemma \texorpdfstring{\protect\ref{kkbar0}}{3.7}(b1)}\label{sec8.8}
When $k\le T-\log^2n$, noticing $\llvert  x_{{n0}}\rrvert  /\llvert  x_{{n1}}\rrvert  \le1-\eta$
established in part (b) of Lemma~\ref{x1}, so (\ref{eqwwk7}) remains
true, hence in turn implies the lemma. When $k>T-\log^2n$, we shall
recursively show the lemma by proving
%
\begin{eqnarray}
\rP\bigl(\llvert W_{k,\ldots,k+s\tau}\rrvert >1-\eta\bigr)=o\bigl(n^{-t}
\bigr), \label{eqb11}
\end{eqnarray}
for some $\eta\in(0,\frac{1}2)$. In fact, when $k+s\tau\ge
T>k+(s-1)\tau$, (\ref{eqb11}) follows easily by the fact that $\bgma
_{k+(s+1)\tau}$ is independent of $\bbA_{k,\ldots,k+s\tau}^{-1}$,
and hence $\rP(\llvert  W_{k,\ldots,k+s\tau}-a_{n}\rrvert  \ge v_{n}^{3})=o(n^{-t})$
and $\llvert  a_{n}\rrvert  \le1/2-\eta$.

By induction, assume that (\ref{eqb11}) is true for some $s\ge1$. By
(\ref{eqwwk1}) and Lemma~\ref{x1}(v), when $\llvert  r_{1}(k+s\tau)\rrvert  \le
v_{n}^{3}$ and $\llvert  r_{2}(k+s\tau)\rrvert  \le v_{n}^{3}$, we have
\begin{eqnarray*}
\llvert W_{k,\ldots,k+(s-1)\tau}\rrvert \le\frac{\sfrac{1}2-\eta+v_n^3}{1-(\sfrac{1}2-\eta)(1-\eta)-v_n^3}\le1-\eta\qquad\mbox{for
all large } n.
\end{eqnarray*}
Thus,
\begin{eqnarray*}
&&\rP\bigl(\llvert W_{k,\ldots,k+(s-1)\tau}\rrvert >1-\eta\bigr)
\\
&&\qquad \le\rP\bigl(\llvert W_{k,\ldots,k+s\tau}\rrvert >1-\eta\bigr)+\rP\bigl(\bigl
\llvert r_1(k+s\tau)\bigr\rrvert \ge v_n^3
\bigr)+\rP\bigl(\bigl\llvert r_2(k+s\tau)\bigr\rrvert \ge
v_n^3\bigr)
\\
&&\qquad =o\bigl(n^{-t}\bigr).
\end{eqnarray*}
The assertion (\ref{eqb11}) is proved, and thus the proof of the lemma
is complete.

\subsection{Proof of Lemma \texorpdfstring{\protect\ref{kkbar}}{3.9}}\label{sec8.9}
Define $\widetilde\bbA_k=\bbA_{k,k+\tau}+\bgma_{k+\tau}\bgma
_{k+2\tau}^*$. Recall $\bbA_k=\bbA_{k,k+\tau}+\bgma_{k+\tau}\bgma
_{k+2\tau}^*+\bgma_{k+2\tau}\bgma_{k+\tau}^*$, so we have
\begin{eqnarray*}
\bbA_k^{-1}=\bigl(\widetilde\bbA_k+
\bgma_{k+2\tau}\bgma_{k+\tau
}^*\bigr)^{-1}=\widetilde
\bbA_k^{-1}-\frac{\widetilde\bbA_k^{-1}\bgma
_{k+2\tau}\bgma_{k+\tau}^*\widetilde\bbA_k^{-1}}{1+\bgma_{k+\tau
}^*\widetilde\bbA_k^{-1}\bgma_{k+2\tau}}.
\end{eqnarray*}
Hence, we have
\begin{eqnarray*}
\bgma_{k+\tau}^*\bbA_k^{-1}=\bgma_{k+\tau}^*
\widetilde\bbA _k^{-1}-\frac{\bgma_{k+\tau}^*\widetilde\bbA_k^{-1}\bgma_{k+2\tau
}\bgma_{k+\tau}^*\widetilde\bbA_k^{-1}}{1+\bgma_{k+\tau
}^*\widetilde\bbA_k^{-1}\bgma_{k+2\tau}} =
\frac{\bgma_{k+\tau}^*\widetilde\bbA_k^{-1}}{1+\bgma_{k+\tau
}^*\widetilde\bbA_k^{-1}\bgma_{k+2\tau}}.
\end{eqnarray*}
Next, we have
\begin{eqnarray*}
\bgma_{k+\tau}^*\widetilde\bbA_k^{-1}&=&
\bgma_{k+\tau}^*\bbA _{k,k+\tau}^{-1}-\frac{\bgma_{k+\tau}^*\bbA_{k,k+\tau}^{-1}\bgma
_{k+\tau}\bgma_{k+2\tau}^*\bbA_{k,k+\tau}^{-1}}{1+\bgma_{k+2\tau
}^*\bbA_{k,k+\tau}^{-1}\bgma_{k+\tau}}
\\
&=&\bgma_{k+\tau}^*\bbA_{k,k+\tau}^{-1}-a_n
\bgma_{k+2\tau}^*\bbA _{k,k+\tau}^{-1}+R_{k1},
\end{eqnarray*}
where
\begin{eqnarray*}
R_{k1}&=&a_n\bgma_{k+2\tau}^*\bbA_{k,k+\tau}^{-1}-
\frac{\bgma
_{k+\tau}^*\bbA_{k,k+\tau}^{-1}\bgma_{k+\tau}\bgma_{k+2\tau
}^*\bbA_{k,k+\tau}^{-1}}{1+\bgma_{k+2\tau}^*\bbA_{k,k+\tau
}^{-1}\bgma_{k+\tau}}
\\
&=& \biggl(\frac{a_n-\bgma_{k+\tau}^*\bbA_{k,k+\tau}^{-1}\bgma
_{k+\tau}+a_n\bgma_{k+2\tau}^*\bbA_{k,k+\tau}^{-1}\bgma_{k+\tau
}}{1+\bgma_{k+2\tau}^*\bbA_{k,k+\tau}^{-1}\bgma_{k+\tau}} \biggr)\bgma_{k+2\tau}^*
\bbA_{k,k+\tau}^{-1}.
\end{eqnarray*}
Substituting back, we obtain
%
\begin{eqnarray}\label{eqbnk}
\bgma_{k+\tau}^*\bbA_k^{-1}
&=& \frac{\bgma_{k+\tau}^*\bbA
_{k,k+\tau}^{-1}-a_n\bgma_{k+2\tau}^*\bbA_{k,k+\tau
}^{-1}+R_{k1}}{1+\bgma_{k+\tau}^*\bbA_{k,k+\tau}^{-1}\bgma
_{k+2\tau}-a_n\bgma_{k+2\tau}^*\bbA_{k,k+\tau}^{-1}\bgma_{k+2\tau
}+R_{k1}\bgma_{k+2\tau}}\nonumber
\\
&=&
\bigl(\bgma_{k+\tau}^*\bbA_{k,k+\tau}^{-1}-a_n\bgma_{k+2\tau
}^*\bbA_{k,k+\tau}^{-1}+R_{k1}\bigr)
\nonumber\\[-8pt]\\[-8pt]\nonumber
&&{} /
\bigl(x_{n1}+\bgma_{k+\tau}^*\bbA
_{k,k+\tau}^{-1}\bgma_{k+2\tau}
\\
&&\hspace*{6pt}{}-a_n\bigl(\bgma_{k+2\tau}^*\bbA
_{k,k+\tau}^{-1}\bgma_{k+2\tau}-a_n/x_{n1}\bigr) +R_{k1}\bgma_{k+2\tau}\bigr). \nonumber
\end{eqnarray}
When $\llvert  \bgma_{k+2\tau}^*\bbA_{k,k+\tau}^{-1}\bgma_{k+\tau}\rrvert  \le
v_n^3$, $\llvert  a_n-\bgma_{k+\tau}^*\bbA_{k,k+\tau}^{-1}\bgma_{k+\tau
}\rrvert  \le v_n^3$, we have
\begin{eqnarray*}
\llVert R_{k1}\rrVert \le Kv_n^2.
\end{eqnarray*}
Using similar approach of the proof of Lemma~\ref{Ek+tau}(a), one can
prove that when $k\le T-\log^2n$,
$\llvert  \bgma_{k+l\tau}^*\bbA^{-1}_{k,\ldots,k+l\tau}\bgma_{k+(l+1)\tau
}\rrvert  \le v_n^3$, $\llvert  \bgma_{k+(l+1)\tau}^*\times\break \bbA^{-1}_{k,\ldots,k+l\tau
}\bgma_{k+l\tau}\rrvert  \le v_n^3$, and $\llvert  
\bgma_{k+l\tau}^*\bbA_{k,\ldots,k+l\tau}^{-1}\bgma_{k+l\tau}-a_n\rrvert
  \le v_n^3$, for $l=1,\ldots,\break [\log
^2n]$, we have
\begin{eqnarray*}
\rP\bigl(\bigl\llvert \bgma_{k+2\tau}^*\bbA_{k,k+\tau}^{-1}
\bgma_{k+2\tau
}-a_n/x_{n1}\bigr\rrvert \ge
v_n^3\bigr)=o\bigl(n^{-t}\bigr).
\end{eqnarray*}

Therefore, by (\ref{eqbnk}), we have
%
\begin{eqnarray}
\qquad \bigl\llVert \bgma_{k+\tau}^*\bbA_{k}^{-1}\bigr
\rrVert \le2\bigl\llVert \bgma_{k+\tau}^*\bbA _{k,k+\tau}^{-1}
\bigr\rrVert +\bigl(1-\eta'\bigr)\bigl\llVert \bgma_{k+2\tau}^*
\bbA_{k,k+\tau
}^{-1}\bigr\rrVert +Kv_n. \label{eqbnk2}
\end{eqnarray}
Similarly, one can prove that
%
\begin{eqnarray}\label{eqbnk3}
\qquad &&\bigl\llVert \bgma_{k+2\tau}^*\bbA_{k,k+\tau}^{-1}\bigr
\rrVert
\nonumber\\[-8pt]\\[-8pt]\nonumber
&&\qquad \le 2\bigl\llVert \bgma_{k+2\tau}^*\bbA_{k,k+\tau,k+2\tau}^{-1}
\bigr\rrVert +\bigl(1-\eta '\bigr)\bigl\llVert \bgma_{k+3\tau}^*
\bbA_{k,k+\tau,k+2\tau}^{-1}\bigr\rrVert +Kv_n.
\end{eqnarray}
By induction, for any $k\le T-[\log^2n]$ and $\ell\le[\log^2n]$,
one obtains
%
\begin{eqnarray}\label{eqbnk4}
&&\bigl\llVert \bgma_{k+\tau}\bbA_{k}^{-1}\bigr
\rrVert
\nonumber
\\
&&\qquad \le2\sum_{l=1}^\ell\bigl(1-
\eta'\bigr)^{l-1} \bigl\llVert \bgma_{k+l\tau}^*\bbA
_{k,\ldots,k+l\tau}^{-1}\bigr\rrVert
\\
&&\quad\qquad{}  +\bigl(1-\eta'
\bigr)^\ell\bigl\llVert \bgma_{k+(\ell+1)\tau
}^*\bbA_{k,\ldots,k+\ell\tau}^{-1}
\bigr\rrVert +K\ell v_n, \nonumber
\end{eqnarray}
where $\eta'\in(0,\eta)$ is a constant.
Since
\[
\bigl\llVert \bgma_{k+l\tau}^*\bbA^{-1}_{k,\ldots,k+l\tau}\bigr
\rrVert ^2\to\frac
{c}2\int\frac{1}{(x-u)^2}\,dF_c(x)=:K_1
\]
uniformly for $k\le T+\tau-[\log^2n]$ and $l\le[\log^2n]$, then for
any $K>\frac{2\sqrt{K_1+\ep}}{\eta'}$, when $n$ is large, we have
%
\begin{eqnarray}\label{eqbnk215}
&& \rP\bigl(\bigl\llVert \bgma_{k+\tau}^*\bbA_{k}^{-1}
\bigr\rrVert \ge K\bigr)\nonumber
\\
&&\qquad \le \sum_{l=1}^{[\log
^2n]}\bigl[ \rP\bigl(\bigl\llvert \bgma_{k+(l+1)\tau}^*\bbA^{-1}_{k,\ldots,k+l\tau
}\bgma_{k+l\tau}\bigr\rrvert \ge v_n^3\bigr)\nonumber
\\
&&\hspace*{60pt}{}+\rP\bigl(\bigl\llvert \bgma_{k+l\tau}^*\bbA^{-1}_{k,\ldots,k+l\tau}
\bgma _{k+(l+1)\tau}\bigr\rrvert \ge v_n^3\bigr)
\\
&&\hspace*{60pt}{}+\rP\bigl(\bigl\llvert \bgma_{k+l\tau}^*\bbA^{-1}_{k,\ldots,k+l\tau}
\bgma _{k+l\tau}-a_n\bigr\rrvert \ge v_n^3
\bigr) \bigr]
\nonumber
\\
&&\qquad =o\bigl(n^{-t}\bigr).\nonumber
\end{eqnarray}
This proves the lemma for $k\le T+\tau-[\log^2n]$.

When $k>T+\tau-[\log^2n]$, by the first equality of (\ref{eqbnk})
and Lemma~\ref{x1}(v), when $\llvert  \bgma_{k+2\tau}^*\bbA^{-1}_{k,k+\tau
}\bgma_{k+2\tau}\rrvert  \le1$ [which, by (\ref{eqb11}), occurs with
probability $1-o(n^{-t})$], we have
\begin{eqnarray*}
&&\bigl\llvert 1+\bgma_{k+\tau}^*\bbA_{k,k+\tau}^{-1}
\bgma_{k+2\tau}-a_n\bgma _{k+2\tau}^*\bbA_{k,k+\tau}^{-1}
\bgma_{k+2\tau}+R_{k1}\bgma _{k+2\tau}\bigr\rrvert
\\
&&\qquad \ge 1-v_n^3-\bigl(\tfrac{1}2-\eta
\bigr)-Kv_n^2\ge\tfrac{1}2+\eta',
\end{eqnarray*}
for some constant $\eta'>0$. Therefore,
\begin{eqnarray*}
\bigl\llVert \bgma_{k+\tau}^*\bbA_k^{-1}\bigr
\rrVert \le2\bigl\llVert \bgma_{k+\tau}^*\bbA _{k,k+\tau}^{-1}
\bigr\rrVert +\bigl(1-\eta'\bigr)\bigl\llVert \bgma_{k+2\tau}^*
\bbA_{k,k+\tau
}^{-1}\bigr\rrVert +Kv_n.
\end{eqnarray*}
Again, by using induction, the lemma can be proved for the case where
$k>T-\log^2n$.

Therefore, the proof of the lemma is complete.
\subsection{Proof of Lemma \texorpdfstring{\protect\ref{kkbar2}}{3.10}}\label{sec8.10}
As in last subsection, we first consider the case $k\le T+\tau-[\log
^2n]$. Note that
\begin{eqnarray*}
\bbA_k^{-1}&=&\bigl(\widetilde\bbA_k+
\bgma_{k+2\tau}\bgma_{k+\tau
}^*\bigr)^{-1}=\widetilde
\bbA_k^{-1}-\frac{\widetilde\bbA_k^{-1}\bgma
_{k+2\tau}\bgma_{k+\tau}^*\widetilde\bbA_k^{-1}}{1+\bgma_{k+\tau
}^*\widetilde\bbA_k^{-1}\bgma_{k+2\tau}},
\\
\widetilde\bbA_k^{-1}&=&\bbA_{k,k+\tau}^{-1}-
\frac{\bbA_{k,k+\tau
}^{-1}\bgma_{k+\tau}\bgma_{k+2\tau}^*\bbA_{k,k+\tau
}^{-1}}{1+\bgma_{k+2\tau}^*\bbA_{k,k+\tau}^{-1}\bgma_{k+\tau}},
\end{eqnarray*}
and
\begin{eqnarray*}
\bgma_{k+\tau}^*\bbA_k^{-1}=\bgma_{k+\tau}^*
\widetilde\bbA _k^{-1}-\frac{\bgma_{k+\tau}^*\widetilde\bbA_k^{-1}\bgma_{k+2\tau
}\bgma_{k+\tau}^*\widetilde\bbA_k^{-1}}{1+\bgma_{k+\tau
}^*\widetilde\bbA_k^{-1}\bgma_{k+2\tau}} =
\frac{\bgma_{k+\tau}^*\widetilde\bbA_k^{-1}}{1+\bgma_{k+\tau
}^*\widetilde\bbA_k^{-1}\bgma_{k+2\tau}}.
\end{eqnarray*}
By similar approach to prove Lemmas~\ref{Ek+tau} and~\ref{kkbar}, we
have
\begin{eqnarray*}
\bigl\llvert \bgma_{k+\tau}^*\bbA_{k,k+\tau}^{-1}\bgma
_{k+2\tau}\bigr\rrvert &\le&v_n^3\qquad\mbox{with
probability } 1-o\bigl(n^{-t}\bigr),
\\
\bigl\llvert \bgma_{k+\tau}^*\bbA_{k,k+\tau}^{-2}
\bgma_{k+2\tau}\bigr\rrvert &\le&v_n^3\qquad\mbox{with
probability } 1-o\bigl(n^{-t}\bigr),
\\
\bigl\llvert \bgma_{k+2\tau}^*\bbA_{k,k+\tau}^{-1}
\bgma_{k+2\tau
}-a_n/x_{n1}\bigr\rrvert &
\le&v_n^3\qquad\mbox{with probability } 1-o
\bigl(n^{-t}\bigr),
\\
\bigl\llvert \bgma_{k+\tau}^*\bbA_{k,k+\tau}^{-1}
\bgma_{k+\tau}-a_n\bigr\rrvert &\le &v_n^3\qquad\mbox{with probability } 1-o\bigl(n^{-t}\bigr).
\end{eqnarray*}
By Remark~\ref{lkj2},
\begin{eqnarray*}
\bgma_{k+\tau}^*\bbA_{k,k+\tau}^{-2}\bgma_{k+\tau
}=
\frac{1}{2T}\rtr\bbA^{-2}+o\bigl(v_n^3
\bigr)\le K\qquad\mbox{with probability } 1-o\bigl(n^{-t}\bigr).
\end{eqnarray*}
By Lemma~\ref{kkbar},
\begin{eqnarray}
\bigl\llVert \bgma_{k+2\tau}^*\bbA_{k,k+\tau}^{-1}\bigr
\rrVert ^2&=&\bigl\llvert \bgma_{k+2\tau
}^*\bbA_{k,k+\tau}^{-1}
\bigl(\bbA^*_{k,k+\tau}\bigr)^{-1}\bgma_{k+2\tau
}\bigr\rrvert
\le K\nonumber
\\
\eqntext{\mbox{with probability } 1-o\bigl(n^{-t}\bigr),}
\\
\bigl\llvert \bgma_{k+2\tau}^*\bbA_{k,k+\tau}^{-2}
\bgma_{k+2\tau}\bigr\rrvert &\le &\bigl\llvert \bgma_{k+2\tau}^*
\bbA_{k,k+\tau}^{-1}\bigl(\bbA^*_{k,k+\tau
}
\bigr)^{-1}\bgma_{k+2\tau}\bigr\rrvert \le K
\nonumber
\\
\eqntext{\mbox{with
probability } 1-o\bigl(n^{-t}\bigr).}
\end{eqnarray}
By Lemma~\ref{lkj},
\begin{eqnarray}
\bigl\llVert \bgma_{k+\tau}^*\bbA_{k,k+\tau}^{-1}\bigr
\rrVert ^2&=&\frac{1}{2T}\rtr\bbA^{-1}_{k,k+\tau}
\bigl(\bbA^*_{k,k+\tau}\bigr)^{-1}+o\bigl(v_n^3
\bigr)\le K\nonumber
\\
\eqntext{\mbox{with probability } 1-o\bigl(n^{-t}\bigr).}
\end{eqnarray}
Also, we have
\begin{eqnarray*}
\bgma_{k+\tau}^*\widetilde\bbA_k^{-1}
\bgma_{k+2\tau}&=&\bgma _{k+\tau}^*\bbA_{k,k+\tau}^{-1}
\bgma_{k+2\tau}-\frac{\bgma
_{k+\tau}^*\bbA_{k,k+\tau}^{-1}\bgma_{k+\tau}\bgma_{k+2\tau
}^*\bbA_{k,k+\tau}^{-1}\bgma_{k+2\tau}}{1+\bgma_{k+2\tau}^*\bbA
_{k,k+\tau}^{-1}\bgma_{k+\tau}}
\\
&=&-x_{n0}+o\bigl(v_n^3\bigr)\qquad\mbox{with
probability } 1-o\bigl(n^{-t}\bigr).
\end{eqnarray*}
Therefore, with probability $1-o(n^{-t})$, we have
\begin{eqnarray*}
&&\bigl\llVert \bgma_{k+\tau}^*\bbA_k^{-1}
\widetilde\bbA_k^{-1}\bgma _{k+2\tau}
\bgma_{k+\tau}^*\widetilde\bbA_k^{-1}\bigr\rrVert
\\
&&\qquad =\biggl\llVert \frac{\bgma_{k+\tau}^*\widetilde\bbA_k^{-1}}{1+\bgma
_{k+\tau}^*\widetilde\bbA_k^{-1}\bgma_{k+2\tau}}
\\
&&\hspace*{36pt}{}\times  \biggl(\bbA _{k,k+\tau}^{-1}-
\frac{\bbA_{k,k+\tau}^{-1}\bgma_{k+\tau}\bgma
_{k+2\tau}^*\bbA_{k,k+\tau}^{-1}}{1+\bgma_{k+2\tau}^*\bbA
_{k,k+\tau}^{-1}\bgma_{k+\tau}} \biggr)\bgma_{k+2\tau}\bgma_{k+\tau
}^*\widetilde
\bbA_k^{-1}\biggr\rrVert
\\
&&\qquad =\biggl\llvert \frac{1}{1+\bgma_{k+\tau}^*\widetilde\bbA_k^{-1}\bgma
_{k+2\tau}}\biggr\rrvert
\\
&&\quad\qquad{}\times  \biggl\llvert
\bgma_{k+\tau}^* \biggl(\bbA_{k,k+\tau
}^{-1}-
\frac{\bbA_{k,k+\tau}^{-1}\bgma_{k+\tau}\bgma_{k+2\tau
}^*\bbA_{k,k+\tau}^{-1}}{1+\bgma_{k+2\tau}^*\bbA_{k,k+\tau
}^{-1}\bgma_{k+\tau}} \biggr)^2\bgma_{k+2\tau}\biggr\rrvert
\\
&&\quad\qquad{}\times \biggl\llVert \bgma_{k+\tau}^* \biggl(\bbA_{k,k+\tau}^{-1}-
\frac{\bbA
_{k,k+\tau}^{-1}\bgma_{k+\tau}\bgma_{k+2\tau}^*\bbA_{k,k+\tau
}^{-1}}{1+\bgma_{k+2\tau}^*\bbA_{k,k+\tau}^{-1}\bgma_{k+\tau
}} \biggr)\biggr\rrVert
\\
&&\qquad \le M_1
\end{eqnarray*}
for some $M_1>0$.
By Remark~\ref{lkj1},
\begin{eqnarray}
\bigl\llVert \bgma_{k+\tau}^*\bbA_{k,k+\tau}^{-2}\bigr
\rrVert ^2&=&\frac{1}{2T}\rtr\bbA^{-2}\bigl(\bbA^*
\bigr)^{-2}+o\bigl(v_n^3\bigr)\le K\nonumber
\\
\eqntext{\mbox{with probability } 1-o\bigl(n^{-t}\bigr).}
\end{eqnarray}
This implies, with probability $1-o(n^{-t})$
\begin{eqnarray*}
&&\bigl\llVert \bgma_{k+\tau}^*\bbA_k^{-1}
\widetilde\bbA_k^{-1}\bigr\rrVert
\\
&&\qquad =\biggl\llVert \frac{\bgma_{k+\tau}^*}{1+\bgma_{k+\tau}^*\widetilde
\bbA_k^{-1}\bgma_{k+2\tau}} \biggl(\bbA_{k,k+\tau}^{-1}-
\frac{\bbA
_{k,k+\tau}^{-1}\bgma_{k+\tau}\bgma_{k+2\tau}^*\bbA_{k,k+\tau
}^{-1}}{1+\bgma_{k+2\tau}^*\bbA_{k,k+\tau}^{-1}\bgma_{k+\tau
}} \biggr)^2\biggr\rrVert
\\
&&\qquad \le M_2+\llvert b_n\rrvert \bigl\llVert
\bgma_{k+2\tau}^*\bbA_{k,k+\tau}^{-2}\bigr\rrVert
\end{eqnarray*}
for some $M_2>0$ and
\[
b_n=-\frac{\sfrac{c_n\rE m_n}2}{1-\sklfrac{c_n\rE
m_n}2\sklfrac{c_n\rE m_n}{2x_{n1}}}=-\frac{a_n}{x_{n1}}
\] with
\[
\bigg\llvert  \frac{a_n}{x_{n1}}\bigg\rrvert  \le\sqrt{\llvert  x_{n0}\rrvert  /\llvert  x_{n1}\rrvert  }\le\sqrt
{1-\eta}.\vadjust{\goodbreak}
\]
Therefore, we have
\begin{eqnarray*}
&&\bigl\llVert \bgma_{k+\tau}^*\bbA_k^{-2}\bigr
\rrVert
\\
&&\qquad =\biggl\llVert \bgma_{k+\tau}^*\bbA_k^{-1}
\biggl(\widetilde\bbA _k^{-1}-\frac{\widetilde\bbA_k^{-1}\bgma_{k+2\tau}\bgma_{k+\tau
}^*\widetilde\bbA_k^{-1}}{1+\bgma_{k+\tau}^*\widetilde\bbA
_k^{-1}\bgma_{k+2\tau}} \biggr)
\biggr\rrVert
\\
&&\qquad \le\bigl\llVert \bgma_{k+\tau}^*\bbA_k^{-1}
\widetilde\bbA _k^{-1}\bigr\rrVert +\biggl\llvert
\frac{1}{1+\bgma_{k+\tau}^*\widetilde\bbA
_k^{-1}\bgma_{k+2\tau}}\biggr\rrvert \bigl\llVert \bgma_{k+\tau}^*\bbA
_k^{-1}\widetilde\bbA_k^{-1}
\bgma_{k+2\tau}\bgma_{k+\tau
}^*\widetilde\bbA_k^{-1}
\bigr\rrVert
\\
&&\qquad \le (2+\ep)M_1+M_2+\sqrt{1-\eta}\bigl\llVert
\bgma_{k+2\tau}^*\bbA _{k,k+\tau}^{-2}\bigr\rrVert ,
\end{eqnarray*}
where $\ep>0$ is a constant.
Then similar to the proof of Lemma~\ref{kkbar}, using the recursion
above we have
\begin{eqnarray*}
\rP\bigl(\bigl\llvert \bgma_{k+\tau}^*\bbA_k^{-2}
\bigl(\bbA^*_k\bigr)^{-2}\bgma_{k+\tau
}\bigr\rrvert
\ge K\bigr)=o\bigl(n^{-t}\bigr)
\end{eqnarray*}
for some $K>0$. When $k>T-\log^2n$, one can similarly prove the
inequality above. The proof of the lemma is complete.

\subsection{Proof of Lemma \texorpdfstring{\protect\ref{ep40}}{3.11}}\label{sec8.11}

We first consider the case where $\log^2n<k<T-\log^2n$. Note that
$\bbA=\bbA_k+\bgma_k\bbbeta_k^*+\bbbeta_k\bgma_k^*$, where
$\bbbeta_k=\bgma_{k-\tau}+\bgma_{k+\tau}$.
We have
%
\begin{eqnarray} \label{eqexpand}
&& \rtr\bbA_k^{-1}-\rtr\bbA^{-1}\nonumber
\\
&&\qquad =
\frac{d}{dz}\log \bigl((1+\ep _1) (1+\ep_2)-
\bgma_k^*\bbA_k^{-1}\bgma_k
\bbbeta_k^*\bbA _k^{-1}\bbbeta_k
\bigr)
\nonumber\\[-8pt]\\[-8pt]\nonumber
&&\qquad =\frac{d}{dz}\log \biggl((1+\ep_1) (1+\ep_2)-(
\ep_3+\ep _4+a_n) \biggl(\ep_5+
\frac{2a_n}{x_{n1}}\biggr) \biggr)
\nonumber
\\
&&\qquad  =\frac{d}{dz}\log \biggl(x_{n1}-x_{n0}+
\ep_1+\ep_2+\ep_1\ep _2-a_n
\ep_5-\biggl(\frac{2a_n}{x_{n1}}+\ep_5\biggr) (
\ep_3+\ep_4) \biggr),\nonumber\hspace*{-10pt}
\end{eqnarray}
where $\ep_i$'s are defined in (\ref{defep}). Note that
\[
\rE(\ep_i\mid \bgma_j,j\ne k)=0\qquad\mbox{for
}i=1,2,3. %
\]
Therefore, by Taylor's expansion, Cauchy integral and Lemma~\ref{x1}
part~(b), we have
%
\begin{eqnarray}\label{ineqep}
&&\biggl\llvert \rE\bigl(\rtr\bbA_k^{-1}-\rtr
\bbA^{-1}\bigr)-\frac{d}{dz}\log (x_{n1}-x_{n0})
\biggr\rrvert
\nonumber
\\
&&\qquad \le\biggl\llvert \frac{d}{dz}\rE \biggl[\log \biggl(1+
\frac{\ep_1+\ep
_2+\ep_1\ep_2-a_n\ep_5-(\sklfrac{2a_n}{x_{n1}}+\ep_5)(\ep_3+\ep
_4)}{x_{n1}-x_{n0}} \biggr)
\nonumber\\[-8pt]\\[-8pt]\nonumber
&&\hspace*{197pt}{} -\frac{\ep_1+\ep_2}{x_{n1}-x_{n0}}-\frac{2\ep_3
a_n}{x_{n1}(x_{n1}-x_{n0})} \biggr]\biggr
\rrvert
\nonumber
\\
&&\qquad \le Kv^{-1}_n\sup_{\llvert  \xi-z\rrvert  =v_n/2} \Biggl[\sum
_{i=1}^5 \bigl(\rE \bigl\llvert
\ep_i^2(\xi)\bigr\rrvert \bigr)+\bigl\llvert \rE
\ep_4(\xi)\bigr\rrvert +\bigl\llvert \rE\ep_5(\xi)\bigr
\rrvert \Biggr].\nonumber
\end{eqnarray}
By applying Lemmas~\ref{kkbar} and~\ref{kkbar2}, one can easily
verify that
%
\begin{equation}
\rE\bigl\llvert \ep_i^2(\xi)\bigr\rrvert =O
\bigl(n^{-1}\bigr)\qquad\mbox{for } i=1,2,3. \label{ineq2}
\end{equation}
Also, by (\ref{eqlacing1}),
%
\begin{equation}
\bigl\llvert \rE\ep_4(\xi)\bigr\rrvert =\biggl\llvert
\frac{1}{2T}\rE\bigl(\rtr\bbA_k^{-1}(\xi)-\rtr\bbA
^{-1}(\xi)\bigr)\biggr\rrvert \le\frac{K}{Tv_n}, \label{ineq3}
\end{equation}
and similar to the proof of (\ref{eqw2})
%
\begin{equation}
\qquad \bigl\llvert \rE\ep_4^2(\xi)\bigr\rrvert \le
\frac{1}{4T^2}\rE\bigl\llvert \rtr\bbA_k^{-1}(\xi )-\rE
\rtr\bbA_k^{-1}(\xi)\bigr\rrvert ^2+\bigl\llvert
\rE\ep_4(\xi)\bigr\rrvert ^2=O\biggl(\frac{1}{n}
\biggr). \label{ineq4}
\end{equation}
By the proof of Lemma~\ref{Ek+tau}(a) with noticing
$\llvert  x_{n0}/x_{n1}\rrvert  \le1-\eta$, when $\log^2n\le k\le T-\log^2n$, for
$i=1,2$, one can prove that
%
\begin{eqnarray}\label{ineq5}
\rE\biggl\llvert \bgma_{k+\tau}^*\bbA_k^{-1}
\bgma_{k+\tau}-\frac
{a_n}{x_{n1}}\biggr\rrvert ^i&=&o(1),
\nonumber\\[-8pt]\\[-8pt]\nonumber
\rE\biggl\llvert \bgma_{k-\tau}^*\bbA_k^{-1}\bgma
_{k-\tau}-\frac{a_n}{x_{n1}}\biggr\rrvert ^i&=&o(1),
\end{eqnarray}
and by the proof of Lemma~\ref{cross}(a),
%
\begin{equation}
\rE\bigl\llvert \bgma_{k-\tau}^*\bbA_{k}^{-1}
\bgma_k\bigr\rrvert ^i=o(1),\qquad\bigl\llvert \rE\bgma
_{k+\tau}^*\bbA_k^{-1}\bgma_{k-\tau}\bigr
\rrvert ^i=o(1). \label{ineq6}
\end{equation}
inequalities (\ref{ineq5}) and (\ref{ineq6}) imply that
%
\begin{eqnarray}
\rE\bigl\llvert \ep_5(\xi)\bigr\rrvert ^i=o(1).
\label{ineq7}
\end{eqnarray}
Combining (\ref{ineqep}), (\ref{ineq2}), (\ref{ineq3}), (\ref
{ineq4}) and (\ref{ineq7}), the first conclusion of Lemma~\ref{ep40}
is proved when $\log^2n\le k\le T-\log^2n$. If $k>T-\log^2n$, by
Lemmas~\ref{kkbar0}(b1)
and~\ref{cross}(a), one may modify the right-hand sides of (\ref
{ineq5})--(\ref{ineq6}) as $O(1)$. This also proves the lemma. The
conclusion for $k<\log^2n$ can be proved similarly.

The second conclusion of the lemma can be proved similarly. The proof
of the lemma is complete.

\subsection{Proof of Lemma \texorpdfstring{\protect\ref{ep41}}{3.7}(b2)}\label{sec8.12}
We assume that $k< T-\log^{2}n$ and prove the first statement only, as
the second follows by symmetry.
As in the proof of Lemma~\ref{Ek+tau}(a), write $W_{k}=\bgma_{k+\tau
}^*\bbA_k^{-1}\bgma_{k+\tau}$ and $W_{k,k+\tau,\ldots,k+s\tau
}=\bgma_{k+(s+1)\tau}^*\*\bbA_{k,k+\tau,\ldots,k+s\tau}^{-1}\bgma
_{k+(s+1)\tau}$. Then by (\ref{eqwwk1}), we have
\[
W_{k,\ldots,k+(s-1)\tau}=\frac{a_{n}+r_1(k+s\tau
)}{1-a_{n}W_{k,\ldots,k+s\tau}+r_2(k+s\tau)}, %
\]
where
\begin{eqnarray*}
r_1(k+s\tau) &=&\bgma_{k+s\tau}^*\bbA_{k,\ldots,k+s\tau}^{-1}
\bgma_{k+s\tau
}-a_{n},
\\
r_2(k+s\tau)&=&-\bigl(
\bgma_{k+s\tau}^*\bbA_{k,\ldots,k+s\tau
}^{-1}\bgma_{k+s\tau}-a_{n}
\bigr)\bgma_{k+(s+1)\tau}^*\bbA_{k,\ldots,k+s\tau}^{-1}
\bgma_{k+(s+1)\tau}
\\
&&{}+\bgma_{k+(s+1)\tau}^*\bbA_{k,\ldots,k+s\tau}^{-1}\bgma
_{k+s\tau}+\bgma_{k+s\tau}^*\bbA_{k,\ldots,k+(s+1)\tau}^{-1}\bgma
_{k+s\tau}
\\
&&{} +\bgma_{k+(s+1)\tau}^*\bbA_{k,\ldots,k+s\tau}^{-1}\bgma
_{k+s\tau}\bgma_{k+s\tau}^*\bbA_{k,\ldots,k+(s+1)\tau}^{-1}\bgma
_{k+s\tau}.
\end{eqnarray*}
Therefore, we have
%
\begin{eqnarray}\label{eqbb1}
&&W_k-\frac{a_{n}}{x_{n1}}
\nonumber
\\
&&\qquad =\frac{a_{n}+r_1(k+\tau)}{1-a_{n}W_{k,k+\tau}+r_2(k+\tau)}-\frac
{a_{n}}{x_{n1}}
\nonumber\\[-8pt]\\[-8pt]\nonumber
&&\qquad =\frac{r_1(k+\tau)}{1-a_{n}W_{k,k+\tau}+r_2(k+\tau)}-\frac
{a_{n}r_2(k+\tau)}{x_{n1}(1-a_{n}W_{k,k+\tau}+r_2(k+\tau))}
\nonumber
\\
&&\quad\qquad{}  + \frac{a_{n}^2(W_{k,k+\tau}-\sklfrac
{a_{n}}{x_{n1}})}{x_{n1}(1-a_{n}W_{k,k+\tau}+r_2(k+\tau))}.\nonumber
\end{eqnarray}
By Lemma~\ref{ep40}, when $k+s\tau\le T$,
\begin{eqnarray*}
\bigl\llvert \rE r_1(k+s\tau)\bigr\rrvert &=&\biggl\llvert
\frac{1}{2T}\rE\rtr\bbA_{k,\ldots,k+s\tau
}^{-1}-a_n\biggr
\rrvert =O \biggl(\frac{s}n \biggr)=O \biggl(\frac{\log
^2n}{n}
\biggr).
\end{eqnarray*}

Using this estimate together with Lemmas~\ref{rs} and~\ref{kkbar},
one can prove that
%
\begin{eqnarray}
&& \rE\bigl(\bigl\llvert r_1(k+s\tau)\bigr\rrvert ^p
\bigr)\nonumber
\\
&&\qquad \le K \bigl(\bigl\llvert \rE r_1(k+s\tau)\bigr\rrvert
^p+\rE \bigl\llvert r_1(k+s\tau)-\rE
r_1(k+s\tau)\bigr\rrvert ^p \bigr)
\nonumber\\[-8pt]\\[-8pt]\nonumber
&&\qquad \le K \bigl(n^{-p}\log^{2p}n+n^{-p}\rE \bigl(
\rtr\bbA_{k,\ldots,k+s\tau}^{-1}\bigl(\bbA_{k,\ldots,k+s\tau}^*
\bigr)^{-1} \bigr)^{p/2} \bigr)
\nonumber
\\
&&\qquad \le Kn^{-p/2},\nonumber
\end{eqnarray}
which implies that for any fixed $\gd>0$,
%
\begin{eqnarray}
\rP\bigl(\bigl\llvert r_1(k+s\tau)\bigr\rrvert \ge
n^{-0.5+\gd}\bigr)=o\bigl(n^{-t}\bigr). \label{eqbb2}
\end{eqnarray}
By this and Lemmas~\ref{kkbar0}(b1) and~\ref{rs}, one can prove that
%
\begin{eqnarray}
\rP\bigl(\bigl\llvert r_2(k+s\tau)\bigr\rrvert \ge n^{-0.5+\gd}
\bigr)=o\bigl(n^{-t}\bigr). \label{eqbb3}
\end{eqnarray}
In Section~\ref{sec3}, we have proved that with probability
$1-o(n^{-t})$, $\llvert  W_{k,k+\tau}-\frac{a_n}{x_{n1}}\rrvert  \le v_n^6$. Also by
Lemma~\ref{x1}(ii)(b), we have $\llvert  x_{n1}\rrvert  \ge\frac{1}2$
which implies that $\llvert  \frac{1}{1-a_nW_{k,k+\tau}+r_2(k+\tau
)}\rrvert  $ is bounded by 3 with probability $1-o(n^{-t})$.

Moreover, by the fact that $\llvert  \frac{a_n}{x_{n1}}\rrvert  =\sqrt{\llvert  \frac
{x_{n0}}{x_{n1}}\rrvert  }\le\sqrt{1-\eta}<1-\frac{1}2\eta$, we have, with
probability $1-o(n^{-t})$,
\begin{eqnarray*}
\biggl\llvert \frac{a_n}{1-a_nW_{k,k+\tau}+r_2(k+\tau)}\biggr\rrvert &\le&\frac
{\llvert  a_n\rrvert  }{\llvert  x_{n1}\rrvert  -v_n^4}\le
\frac{(1-\sklfrac{1}2\eta)\llvert  x_{n1}\rrvert  }{\llvert  x_{n1}\rrvert  -v_n^4}
\\
&\le&\frac{1-\sklfrac{1}2\eta}{1-2v_n^4}\le1-\eta',
\end{eqnarray*}
for some $0<\eta'<\frac{1}2\eta$.
In (\ref{eqbb1}), split the first term as
\begin{eqnarray*}
&& \frac{r_1(k+\tau)}{1-a_{n}W_{k,k+\tau}+r_2(k+\tau)}
\\
&&\qquad =\frac
{r_1(k+\tau)}{1-a_{n}W_{k,k+\tau}}\\
&&\qquad\quad{}-\frac{r_1(k+\tau)r_2(k+\tau
)}{(1-a_nW_{k,k+\tau})
(1-a_{n}W_{k,k+\tau}+r_2(k+\tau))} %
\end{eqnarray*}
and the second term as
\begin{eqnarray*}
&& \frac{a_{n}r_2(k+\tau)}{x_{n1}(1-a_{n}W_{k,k+\tau}+r_2(k+\tau
))}
\\
&&\qquad =\frac{a_{n}r_2(k+\tau)}{x_{n1}(1-a_{n}W_{k,k+\tau})}\\
&&\qquad\quad{}-\frac
{a_{n}r_2^2(k+\tau)}{x_{n1}(1-a_{n}W_{k,k+\tau})(1-a_{n}W_{k,k+\tau
}+r_2(k+\tau))}.
\end{eqnarray*}
Noting that $\llvert  W_k\rrvert  \le Kv^{-1}_n$, we have
%
\begin{eqnarray}
\label{eq1} &&\biggl\llvert \rE W_k-\frac{a_n}{x_{n1}}\biggr
\rrvert
\nonumber
\\
&&\qquad \le Kn^{-1+2\gd}+K\bigl\llvert \rE r_1(k+\tau)\bigr\rrvert
+K\bigl\llvert \rE r_2(k+\tau)\bigr\rrvert \nonumber
\\
&&\qquad\quad{} +\bigl(1-\eta
'\bigr)^2\biggl\llvert \rE W_{k,k+\tau}-
\frac{a_n}{x_{n1}}\biggr\rrvert
\nonumber\\[-8pt]\\[-8pt]\nonumber
&&\qquad \vdots
\nonumber
\\
&&\qquad \le K\ell n^{-1+2\gd}+K\sum_{s=1}^{\ell}
\bigl\llvert \rE r_1(k+s\tau)\bigr\rrvert +K\sum
_{s=1}^{\ell}\bigl\llvert \rE r_2(k+s\tau)
\bigr\rrvert
\nonumber
\\
&&\quad\qquad{}+\bigl(1-\eta'\bigr)^{2\ell}\biggl\llvert \rE
W_{k,\ldots,k+\ell\tau}-\frac
{a_n}{x_{n1}}\biggr\rrvert .\nonumber
\end{eqnarray}
%
By choosing $\ell=[\log^2n]$ and $\gd<1/106$, we can show that $\sum_{s=1}^{\ell}\llvert  \rE r_i(k+s\tau)\rrvert  =o(1/(nv_n))$, $i=1,2$
and that $(1-\eta')^{2\ell}\llvert  \rE W_{k,\ldots,k+\ell\tau
}-\frac{a_n}{x_{n1}}\rrvert  =o(1/(nv_n))$.
Substituting all the above into (\ref{eq1}), we have $\llvert  \rE W_k-\frac
{a_n}{x_{n1}}\rrvert  =o(1/(nv_n))$.

\subsection{Proof of Lemma \texorpdfstring{\protect\ref{ep43}}{3.7}(b3)}\label{sec8.13}
Again, we assume that $k< T-\log^{2}n$ and prove the first statement
only, as the second follows by symmetry. As in the proof of Lemma~\ref
{ep41}(b2), we have
%
\begin{eqnarray}
\label{eq2} &&\rE\biggl\llvert W_k-\frac{a_n}{x_{n1}}\biggr
\rrvert ^2
\nonumber
\\
&&\qquad \le K\rE\bigl\llvert r_1(k+\tau)\bigr\rrvert ^2+K\rE
\bigl\llvert r_2(k+\tau)\bigr\rrvert ^2+\bigl(1-
\eta'\bigr)^4\rE \biggl\llvert W_{k,k+\tau}-
\frac{a_n}{x_{n1}}\biggr\rrvert ^2
\nonumber
\\
&&\qquad \vdots
\nonumber\\[-8pt]\\[-8pt]\nonumber
&&\qquad \le K\sum_{s=1}^{\ell}\rE\bigl\llvert
r_1(k+s\tau)\bigr\rrvert ^2+K\sum
_{s=1}^{\ell}\rE \bigl\llvert r_2(k+s\tau)
\bigr\rrvert ^2\nonumber
\\
&&\quad\qquad{} +\bigl(1-\eta'\bigr)^{4\ell}\rE
\biggl\llvert W_{k,\ldots,k+\ell
\tau}-\frac{a_n}{x_{n1}}\biggr\rrvert ^2
\nonumber
\\
&&\qquad \le K\ell n^{-1+2\gd}=o\bigl(1/(nv_n)\bigr).\nonumber
\end{eqnarray}
The proof of the lemma is complete.
\subsection{Proof of Lemma \texorpdfstring{\protect\ref{ep42}}{3.8}(b1)}\label{sec8.14}
By symmetry, we only consider the case $k\le T/2$. As in the proof of
Lemma~\ref{cross}(a), write
\[
\widetilde W_{k,\ldots,k+s\tau}:=\bgma_{k-\tau}^*\bbA_{k,k+\tau,\ldots,k+(s-1)\tau}^{-1}
\bgma_{k+s\tau}.
\]
Then we have
%
\begin{equation}
\widetilde W_{k,\ldots,k+s\tau}=\frac{\widetilde r_1(k+s\tau
)-\widetilde W_{k,\ldots,k+(s+1)\tau}(a_n+\widetilde r_2(k+s\tau))}{1+
r_2(k+s\tau)-a_nW_{k,\ldots,k+s\tau}}, \label{eqbb430}
\end{equation}
where
\begin{eqnarray*}
\widetilde r_1(k+s\tau)&=&\bgma_{k-\tau}^*
\bbA_{k,\ldots,k+s\tau
}^{-1}\bgma_{k+s\tau}\bigl(1+
\bgma_{k+s\tau}^*\bbA_{k,\ldots,k+s\tau
}^{-1}\bgma_{k+(s+1)\tau}
\bigr),
\\
\widetilde r_2(k+s\tau)&=&\bgma_{k+s\tau}^*
\bbA_{k,\ldots,k+s\tau
}^{-1}\bgma_{k+s\tau}-a_n.
\end{eqnarray*}
Similar to the proof of (\ref{eqbb3}), one has
%
\begin{equation}
\rP\bigl(\bigl\llvert \widetilde r_i(k+\tau)\bigr\rrvert \ge
n^{-0.5+\gd}\bigr)=o\bigl(n^{-t}\bigr),\qquad i=1,2. \label{eqbb531}
\end{equation}

Similar to the proof of (\ref{eq1}), one can prove that for some $\eta'>0$,
\begin{eqnarray*}
&&\llvert \rE\widetilde W_{k,\ldots,k+s\tau}\rrvert \\
&&\qquad\le Kn^{-1+2\gd}+K\bigl
\llvert \rE \widetilde r_1(k+s\tau)\bigr\rrvert +\bigl(1-
\eta'\bigr)\llvert \rE\widetilde W_{k,\ldots
,k+(s+1)\tau}\rrvert .
\end{eqnarray*}
Therefore, when $k\le T/2$,
\begin{eqnarray*}
\llvert \rE\widetilde W_{k}\rrvert &\le& K\ell n^{-1+2\gd}+K\sum
_{s=1}^\ell\bigl\llvert \rE \widetilde
r_1(k+s\tau)\bigr\rrvert +\bigl(1-\eta'
\bigr)^\ell\llvert \rE\widetilde W_{k,\ldots
,k+\ell\tau}\rrvert
\\
& =&o
\bigl(1/(nv_n)\bigr).
\end{eqnarray*}
The proof of the lemma is complete.

\subsection{Proof of Lemma \texorpdfstring{\protect\ref{ep44}}{3.8}(b2)}\label{sec8.15}
Using the notation of Lemma~\ref{ep42}(b1), by triangle inequality,
we have
\begin{eqnarray*}
&&\bigl(\rE\llvert \widetilde W_{k+s\tau}\rrvert ^2
\bigr)^{1/2}\\
&&\qquad\le K\bigl(\rE\bigl\llvert \widetilde r_1(k+s
\tau)\bigr\rrvert ^2\bigr)^{1/2}+\bigl(\bigl(1-
\eta'\bigr)\rE\llvert \widetilde W_{k,\ldots
,k+(s+1)\tau}\rrvert
^2\bigr)^{1/2}.
\end{eqnarray*}
Therefore, when $k\le T/2$ and $\ell=[\log^2n]$,
\begin{eqnarray*}
&&\bigl(\rE\llvert \widetilde W_{k}\rrvert ^2
\bigr)^{1/2}\\
&&\qquad\le K\sum_{s=1}^\ell
\bigl(\rE \bigl\llvert \widetilde r_1(k+s\tau)\bigr\rrvert
^2\bigr)^{1/2}+\bigl(1-\eta'
\bigr)^{\ell/2}\bigl(\rE \llvert \widetilde W_{k,\ldots,k+\ell\tau}\rrvert
^2\bigr)^{1/2}
\\
&&\qquad\le K\log^2nn^{-1/2+\gd}.
\end{eqnarray*}
Therefore, when $2\delta<1/212$,
\[
\rE\llvert \widetilde W_{k}\rrvert ^2\le K
\log^4nn^{-1+2\gd}=o\bigl(1/(nv_n)\bigr) %
\]
and the proof of the lemma is complete.
\end{appendix}

\section*{Acknowledgements}
The authors would like to thank the referees for their careful reading
and invaluable comments which greatly improved the quality of the paper.



\printaddresses
\end{document}